\documentclass[3p, 12pt, sort&compress]{elsarticle}
\usepackage{float} 
\usepackage{graphicx}
\usepackage{amsmath}
\usepackage{mathtools}
\usepackage{amssymb}
\usepackage{bm}
\usepackage{framed}
\usepackage{subfig}
\usepackage{hyperref}
\usepackage{nomencl}
\usepackage{color}
\usepackage{textcomp}
\usepackage{graphics}
\usepackage{epsfig}
\usepackage{bm}
\usepackage{epstopdf}
\usepackage{amsthm}
\usepackage{float}
\usepackage{latexsym,amsfonts}
\usepackage{url}
\usepackage{longtable}
\usepackage[figuresright]{rotating}
\usepackage{listings}
\usepackage{algorithm}
\usepackage{algpseudocode}
\usepackage{footnote}
\usepackage{xcolor}
\usepackage[normalem]{ ulem }
\usepackage{soul}

\usepackage{makeidx}
\makeindex

\newtheorem{proposition}{Proposition}

\newdefinition{rmk}{Remark}
\newproof{pf}{Proof}
\newproof{pot}{Proof of Theorem \ref{thm2}}

\newcommand\myeq{\mathrel{\stackrel{\makebox[0pt]{\mbox{\normalfont\tiny def}}}{=}}}

\newcommand{\vecbold}[1]{\bm{#1}}
\renewcommand{\vec}[1]{\vecbold{#1}}

\DeclareMathOperator*{\argmin}{arg\,min}

\definecolor{tocorrect}{rgb}{0.97, 0.04, 0.56}

\newcounter{lnote}

\begin{document}

\begin{frontmatter}
\title{Fast Bayesian experimental design: Laplace-based importance sampling for the expected information gain}

\author[rvt]{Joakim Beck}
\ead{joakim.beck@kaust.edu.sa}

\author[fct]{Ben Mansour Dia}
\ead{mansourben2002@yahoo.fr}

\author[rvt]{Luis FR Espath}
\ead{espath@gmail.com}

\author[utrc]{Quan Long}
\ead{longq@utrc.utc.com}

\author[rvt]{Raul Tempone}
\ead{raul.tempone@kaust.edu.sa}

\address[rvt]{Computer, Electrical and Mathematical Science and Engineering Division (CEMSE), King Abdullah University of Science and Technology (KAUST), Thuwal, 23955-6900, Saudi Arabia}
\address[fct]{CIPR, College of Petroleum Engineering and Geosciences, King Fahd  University of Petroleum and Minerals, Dhahran 31261, Saudi Arabia}
\address[utrc]{United Technologies Research Center, East Hartford, CT, 06108, United States}

\begin{abstract}
In calculating expected information gain in optimal Bayesian experimental design, the computation of the inner loop in the classical double-loop Monte Carlo requires a large number of samples and suffers from underflow if the number of samples is small. These drawbacks can be avoided by using an importance sampling approach. We present a computationally efficient method for optimal Bayesian experimental design that introduces importance sampling based on the Laplace method to the inner loop. We derive the optimal values for the method parameters in which the average computational cost is minimized according to the desired error tolerance. We use three numerical examples to demonstrate the computational efficiency of our method compared with the classical double-loop Monte Carlo, and a more recent single-loop Monte Carlo method that uses the Laplace method as an approximation of the return value of the inner loop. The first example is a scalar problem that is linear in the uncertain parameter. The second example is a nonlinear scalar problem. The third example deals with the optimal sensor placement for an electrical impedance tomography experiment to recover the fiber orientation in laminate composites.
\end{abstract}

\begin{keyword}
Bayesian experimental design, Expected information gain, Monte Carlo, Laplace approximation, Importance sampling, Composite materials.
\end{keyword}

\end{frontmatter}
AMS 2010 subject classification: 62K05, 65N21, 65C60, 65C05




\section{Introduction}
This work proposes an efficient method for the computation of expected information gain \cite{kullback1951, kullback1959, gin} in optimal Bayesian experimental design. The expected information gain, also known as the expected Kullback-Leibler divergence, is an information metric commonly used to estimate the information provided by the proposed experiment. 
Its straightforward estimator, the double-loop Monte Carlo (DLMC), typically requires a large number of samples \cite{ryan, huan} because it embeds a nested sample-average structure. The computation of the inner loop can suffer from arithmetic underflow, especially for small sample sizes.

Another approach is to employ the Laplace method (MCLA) \cite{Stigler1986} to approximate the inner integral of the expected information gain with a second-order Taylor expansion around the mode and then compute the inner integral analytically \cite{tierney1986,tierney1989,kass1990}, which leads to a single integral. Other studies have introduced the Laplace method to the design of experiments for partial differential equation (PDE) models and achieved a substantial improvement in the efficiency compared to DLMC \cite{key15,key51}. The posterior distribution was expanded at the maximum a posteriori (MAP) estimate; therefore, the associated integrals of the inner loop could be approximated analytically by Gaussian integration. However, this leads to an extra bias, except when the Laplace approximation is exact. Furthermore, in \cite{papadimitriou}, the Laplace method was used to approximate the entropy of the posterior distribution for optimal sensor placement, and in \cite{key52}, the method was extended to handle under-determined experiments. A truncated Gaussian has also been used to approximate the posterior distribution, see \cite{long2016}. Related work in optimal experimental design are the Bayesian A-optimality for PDE models \cite{alex2016}, and the use of Gaussian process models for the approximation of the expected information gain in the context of sequential design of computer experiments \cite{BG2016}.

A self-normalized importance sampling approach for the computation of the inner loop of DLMC has been proposed by Feng \cite{F2015}, where a multivariate normal distribution is used as the inner sampling distribution with weighted sample mean and weighted sample covariance based on the outer samples. This leads to a substantial reduction in the number of inner samples. However, the approximation of the normalization constant can suffer from underflow when the posterior distributions are concentrated, in particular when the measurement error is small or the number of repetitive experiments is large.

In this work, we instead propose a Laplace-based importance sampling where the inner sampling distribution is the Laplace approximation on the MAP estimates. These MAP estimates are obtained by solving an optimization problem. The method does not introduce any extra bias, in contrast to the Laplace method, and shows a reduction of orders of magnitude in the number of inner samples. This approach also mitigates the risk of underflow.

We also devise a strategy to estimate optimal values for the method parameters for the desired error tolerance. The method parameters are the outer and inner number of samples, as well as the discretization parameter in the context of PDEs.

The outline of the paper is as follows: Section \ref{section2} deals with optimal Bayesian experimental design where the data model is composed of a deterministic computer model and an additive observational noise component. A brief introduction of the expected information gain criteria is also presented. In Section \ref{sec:DLMC}, we provide an error analysis of the DLMC estimator, and derive optimal values for the method parameters for a given error tolerance. In Section \ref{sec:la}, we perform the minimization of the average computational work of the MCLA estimator. In Section \ref{sec:DLMCIS}, we propose our double-loop Monte Carlo importance sampling (DLMCIS) method, which is proven to substantially reduce the average computational work, mitigate the risk of underflow, and in contrast to the Laplace method it does not introduce any additional bias. In Section \ref{sec:num}, three numerical examples are used to compare the methods in terms of robustness and computational efficiency. The first example is a linear model, the second example is a nonlinear model, and the third example is a sensor placement design problem, where the goal is to maximize signal information during electrical impedance tomography (EIT) in order to inversely obtain the parameters of the inter-ply delaminations.

\section{Optimal Bayesian experimental design}\label{section2}

\subsection{Problem setting}
We consider the data model given by
\begin{equation}
\bm{y}_i(\bm{\xi}) = \bm{g}(\bm{\theta}_t,\bm{\xi}) + \bm{\epsilon}_i, \hspace{0.1cm} i=1,\dots,N_e,
\label{eq_datamodel}
\end{equation}
where $\bm{y}_i \in \mathbb{R}^{q}$ is a vector of $q$ observed responses, $\bm{g}(\bm{\theta}_t,\bm{\xi}) \in \mathbb{R}^q$ is the deterministic model responses, $\bm{\theta}_t \in \mathbb{R}^{d}$ is the true parameter vector, $\bm{\xi} \in \Xi$ is the design parameter vector, $\Xi$ is the experimental design space, $\bm{\epsilon}_i$ are independent and identically distributed (i.i.d.) zero-mean Gaussian measurement errors with covariance matrix $\bm{\Sigma_{\epsilon}}$, and $N_e$ is the number of repetitive experiments. The observed dataset is denoted by $\bm{Y}=\{\bm{y}_i\}^{N_e}_{i=1}$. In our problem setting, we assume the true value of $\bm{\theta}_t$ is unknown. We let $(\Omega,\mathcal{F},\mathbb{P})$ be a complete probability space where $\mathcal{F}$ is the $\sigma$-field of events, $\mathbb{P} : \mathcal{F} \rightarrow [0,1]$ is a probability measure, and $\Omega$ is the set of outcomes. We consider a vector of random variables, $\bm{\theta}: \Theta$ $\mapsto$ $\mathbb{R}^{d}$, with prior space $\Theta$ and prior distribution $\pi(\bm{\theta})$ in lieu of the unknown vector $\bm{\theta}_t$, i.e.,
\begin{equation}\label{model}
\bm{y}_i(\omega_1, \omega_2, \bm{\xi}) = \bm{g}(\bm{\theta}(\omega_1),\bm{\xi}) + \bm{\epsilon}_i(\omega_2), \hspace{0.1cm} i=1,\dots,N_e.
\end{equation}
We note that the data model accounts for parametric but not structural uncertainty; for example, we do not account for the model error in $\bm{g}$ that cannot be eliminated by evaluating $\bm{g}$ at $\bm{\theta}_t$.

We denote the resulting approximation of the forward model $\bm{g}$ using mesh discretization parameter $h$ by $\bm{g}_h$. As $h \rightarrow 0$ asymptotically, the convergence order of $\bm{g}_h$ is given by
$$\mathbb{E} \left[ \left\| \bm{g}(\bm{\theta}) - \bm{g}_h(\bm{\theta})\right\|_2 \right] = \mathcal{O}\left(h^{\eta}\right),$$
where $\eta>0$ is the $h$-convergence rate. The work of $\bm{g}_h$ is assumed to be $\mathcal{O}\left(h^{-\gamma}\right)$, for some $\gamma>0$. We also assume that $\bm{g}$ is twice differentiable with respect to $\bm{\theta}$. The matrix norm given by $ \|\bm{x}\|^{2}_{\bm{\Sigma}^{-1}} = \bm{x}^{T} \bm{\Sigma}^{-1} \bm{x}$ for a vector $\bm{x}$ and covariance matrix $\bm{\Sigma}$ is used throughout.

The objective of optimal Bayesian experimental design is to determine the most informative experimental design setup about $\bm{\theta}_t$, denoted by $\bm{\xi}^* \in \Xi$. The utility function employed is the expectation of the Kullback-Leibler divergence, see, e.g., \cite{huan,ryan,key15}. For the sake of conciseness, $\bm{\xi}$ is omitted until Section \ref{sec:num}, as the expected information gain for each $\bm{\xi}$ is computed separately.

\subsection{Expected information gain}
The Kullback-Leibler divergence, $D_{kl}$, also known as the information gain \cite{key13, key12}, is an entropic function that can quantify our uncertainty about $\bm{\theta}_t$ through the distance between the prior $\pi(\bm{\theta})$ and the posterior $\pi(\bm{\theta} \vert \bm{Y})$ as
\begin{equation}
D_{kl}(\bm{Y}) = \int_{\Theta}{\pi(\bm{\theta} \vert \bm{Y})  \log \left( \frac{\pi(\bm{\theta} \vert \bm{Y})}{\pi(\bm{\theta})} \right)  d\bm{\theta}}, \label{eq_dkl1}
\end{equation}
where $\bm{Y}=(\bm{y}_1,\ldots,\bm{y}_{N_e})$ is the data, and $\pi(\bm{\theta})$ and $\pi(\bm{\theta} \vert \bm{Y})$ are the prior and posterior probability density functions (pdfs), respectively. The larger the value of $D_{kl}$, the more informative the dataset is about the unknown $\bm{\theta}_t$.

The value of $D_{kl}(\bm{Y})$ can be derived exactly when the prior and posterior pdfs are both Gaussian. For the one-dimensional case (d=1), the Gaussian integral formula yields
\begin{eqnarray}
D_{kl}(\bm{Y}) = \log\left(\frac{\sigma_{prior}}{\sigma_{post}(\bm{Y})}\right)+\frac{1}{2}\left[\frac{\sigma_{post}^2(\bm{Y})}{\sigma_{prior}^2}-1+\frac{(\mu_{post}(\bm{Y})-\mu_{prior})^2}{\sigma_{prior}^2} \right],
\end{eqnarray}
where the first and second moments of the prior distribution $\pi(\theta)$ are denoted by $\mu_{prior}$ and $\sigma_{prior}$, and the moments of the posterior pdf $\pi(\theta \vert \bm{Y})$ by $\mu_{post}(\bm{Y})$ and $\sigma_{post}(\bm{Y})$.

In the first stage of experimental design, there are no available observations. Thus, we take the expectation of $D_{kl}$, denoted by $I$, over the sample space, $\mathcal{Y} \subseteq \mathbb{R}^q$, i.e.,
\begin{align}\label{expecinfgain}
I =& \int_{\mathcal{Y}}{\int_{\Theta}{\log \left( \frac{\pi(\bm{\theta} \vert \bm{Y})}{\pi(\bm{\theta})} \right) \pi(\bm{\theta} \vert \bm{Y}) d\bm{\theta}} p(\bm{Y}) d\bm{Y}} \nonumber\\
=& \int_{\Theta} \int_{\mathcal{Y}} \log \left( \frac{p(\bm{Y} \vert \bm{\theta})}{p(\bm{Y})} \right)  p(\bm{Y} \vert \bm{\theta}) d\bm{Y} \pi(\bm{\theta}) d\bm{\theta}.
\end{align}
The latter equality in \eqref{expecinfgain} follows from Bayes' rule. For clarity, we adopt the notation $\pi(\cdot)$ for the pdf of the parameters and $p(\cdot)$ for the pdf of the data sample. In the rest of the paper, we use the approximation of the likelihood with respect to $\bm{g}_h$, given by
\begin{eqnarray}
 p(\bm{Y} \vert \bm{\theta}) =  \left(2\pi \vert \bm{\Sigma_\epsilon} \vert\right)^{-\frac{N_e}{2}} \exp \left( -\frac{1}{2}  \sum_{i=1}^{N_e} \left\| \bm{g}_h(\bm{\theta}_t) + \bm{\epsilon}_i - \bm{g}_h(\bm{\theta})\right\|^2_{\bm{\Sigma_\epsilon}^{-1}} \right).
\end{eqnarray}

\subsection{Fast numerical estimators for expected information gain}
Here we propose a strategy for designing a computationally-efficient numerical estimator $\mathcal{I}$ of $I=\mathbb{E}[D_{kl}]$ that satisfies the tolerance $\hbox{TOL}>0$ at a confidence level given by $0<\alpha \ll 1$:
\begin{eqnarray}
\mathbb{P}\left(\vert I-\mathcal{I} \vert \leq \hbox{TOL} \right) \geq 1-\alpha.
\end{eqnarray}
So the absolute difference $\vert \mathcal{I} - I \vert$ should be less or equal to $\hbox{TOL}$, with a probability of $1-\alpha$.

Similar to \cite{key14,key56}, we optimize the numerical estimator by minimizing its average computational cost based on the $h$-convergence rate, $\eta$, and the work rate, $\gamma$, of the underlying forward problem. Before defining the cost minimization problem, we introduce some notation. The total error is split into a bias component and a statistical error,
\begin{eqnarray}
  \lvert I - \mathcal{I} \rvert \leq \lvert I - \mathbb{E}\left[\mathcal{I}\right] \rvert + \lvert \mathbb{E} \left[\mathcal{I}\right] - \mathcal{I} \rvert,
\end{eqnarray}
and we introduce a balance parameter $\kappa \in ]0,1[$ such that
\begin{eqnarray}
\lvert I - \mathbb{E}\left[\mathcal{I}\right] \rvert \leq (1-\kappa) \hbox{TOL}, \text{and} \\
\nonumber\\
\lvert \mathbb{E} \left[\mathcal{I}\right] - \mathcal{I} \rvert \leq \kappa \hbox{TOL},
\end{eqnarray}
where $(1-\kappa) \hbox{TOL}$ is the bias tolerance and $\kappa \hbox{TOL}$ is the statistical error tolerance. As in a previous study \cite{key56}, we recast the statistical error constraint by using the central theorem limit (CTL) as follows:
\begin{eqnarray}
\mathbb{V}\left[\mathcal{I}\right] \leq \left( \frac{\kappa \hbox{TOL}}{C_\alpha}\right)^2,
\end{eqnarray}
where $C_\alpha = \Phi^{-1}(1-\frac{\alpha}{2})$ and $\Phi^{-1}(\cdot)$ is the inverse cumulative distribution function of the standard normal distribution.  Furthermore, $W$ denotes the average computational work of a single evaluation of $\mathcal{I}$, and the method parameters of the estimator are denoted by $\bm{\zeta}$. For instance, the classical Monte Carlo (MC) estimator has $\bm{\zeta}=\{ N,h \}$, where $N$ is the number of samples and $h$ is the discretization parameter for $\bm{g}_h$.

Within this framework, the optimal setting for the estimator $\mathcal{I}_{\bm{\zeta}}$ is the solution of the cost minimization problem,
\begin{eqnarray}
 \label{costmin}
 (\bm{\zeta}^*,\kappa^*)=\argmin_{(\bm{\zeta},\kappa)} W(\bm{\zeta})\;\;\;\hbox{subject to}\;\;\; 
 \left\{
\begin{array}{lll}
\mathbb{V}\left[\mathcal{I}_{\bm{\zeta}}\right] \leq \left(\kappa \hbox{TOL}/C_\alpha\right)^2\\\\
\lvert I - \mathbb{E}\left[\mathcal{I}_{\bm{\zeta}}\right] \rvert  \leq (1-\kappa ) \hbox{TOL},
\end{array}
\right.
\end{eqnarray}
for the specified tolerance $\hbox{TOL}>0$ at a confidence level given by $1-\alpha$. We note that the balance parameter $\kappa$ is chosen in conjuction with the method parameters. In practice, we need to provide estimates for $W$, $\mathbb{V}\left[\mathcal{I}\right]$, and $\lvert I - \mathbb{E}\left[\mathcal{I}\right] \rvert$, which depend on $\bm{\zeta}$.  Two optimal settings can be derived, with and without a mesh discretization, which lead to different optimal $\boldsymbol{\zeta^*}$. In this work, we consider both settings, but only the derivation with the mesh discretization is presented.

Below we review the DLMC, MCLA, and DLMCIS estimators, and derive their optimal parameter settings.

\section{Double-loop Monte Carlo} 
\label{sec:DLMC}

\subsection{Double-loop Monte Carlo (DLMC) estimator}
The DLMC estimator, $\mathcal{I}_{dl}$, of the expected information gain \eqref{expecinfgain} is given by
\begin{eqnarray}\label{eig_mc}
\mathcal{I}_{dl} \myeq \frac{1}{N} \sum_{n=1}^{N}{\log\left( \frac{p(\bm{Y}_n \vert \bm{\theta}_n)}{\hat{p}_M(\bm{Y}_n)} \right)},
\end{eqnarray}
where $\hat{p}_M(\bm{Y}_n)$ denotes the inner averaging,
\begin{equation}
\hat{p}_M(\bm{Y}_{n}) \myeq \frac{1}{M} \sum_{m=1}^{M}{p(\bm{Y}_n \vert \bm{\tilde{\theta}}_{n,m})},
\end{equation}
and $N$ and $M$ are the respective numbers of samples. The tilde ($\bm{\sim}$) on $\bm{\theta}$ differentiates the inner loop samples from the outer loop samples, as they are independent of each other. The average computational work of $\mathcal{I}_{dl}$ is assumed to follow
\begin{equation}
\label{eq:costmodeldl}
W_{dl} \propto NMh^{-\gamma},
\end{equation}
where $h^{-\gamma}$ is proportional to the average work of a single evaluation of $\bm{g}_h$. The next section shows the optimal setting (i.e., $N$, $M$, and $h$) with respect to the choice of $\hbox{TOL}$.

\subsection{Optimal setting for the DLMC estimator}
The optimal setting for the DLMC estimator is derived by exploiting Proposition \ref{prop1}:
\begin{proposition}
\label{prop1}
The bias and variance of DLMC estimator $\mathcal{I}_{dl}$ can be estimated by:
\begin{eqnarray}
 \lvert I - \mathbb{E}\left[\mathcal{I}_{dl}\right] \rvert \leq C_{dl,3}h^\eta + \frac{C_{dl,4}}{M}+o(h^{\eta})+\mathcal{O}\left(\frac{1}{M^2}\right),\label{biasprop1} \\
\nonumber\\
 \mathbb{V}\left[ \mathcal{I}_{dl} \right] = \frac{C_{dl,1}}{N} + \frac{C_{dl,2}}{NM}+\mathcal{O}\left(\frac{1}{NM^2}\right),\label{statprop1}
\end{eqnarray}
respectively, where
\begin{eqnarray*}
C_{dl,1} & = &  \mathbb{V}\left[ \log \left(\frac{p(\bm{Y} \vert \bm{\theta})}{p(\bm{Y})}\right) \right],\\
C_{dl,2} & = & \left(1+\mathbb{E}\left[ \log \left(\frac{p(\bm{Y} \vert \bm{\theta})}{p(\bm{Y})}\right) \right] \right) \mathbb{E}\left[\mathbb{V} \left[ \frac{ p(\bm{Y} \vert \bm{\theta})}{p(\bm{Y})} \vert \bm{Y}\right]\right]-\mathbb{E}\left[ \log \left(\frac{p(\bm{Y} \vert \bm{\theta})}{p(\bm{Y})}\right)\mathbb{V} \left[ \frac{ p(\bm{Y} \vert \bm{\theta})}{p(\bm{Y})} \vert \bm{Y}\right] \right],  \\
C_{dl,4} & = &  \frac{1}{2}\mathbb{E}\left[\mathbb{V} \left[ \frac{ p(\bm{Y} \vert \bm{\theta})}{p(\bm{Y})} \vert \bm{Y}\right] \right],
\end{eqnarray*}
and $C_{dl,3}$ is the constant of the $h$-convergence of $\mathcal{I}_{dl}$.
\end{proposition}

The proof of Proposition \ref{prop1} is given in \ref{ap:secA}. We remark that a second-order Taylor expansion is used in the proof, instead of a first order as was used elsewhere \cite{key59}, which results in a better estimate of $C_{dl,2}$. The terms $C_{dl,3}h^\eta$ and $\frac{C_{dl,4}}{M}$ respectively account for the biases due to the numerical discretization of the forward problem and due to the error of the inner averaging that arise through the nonlinearity of the logarithmic function. We insert \eqref{biasprop1}, \eqref{statprop1}, and \eqref{eq:costmodeldl} into \eqref{costmin}, and the optimal parameter setting of the estimator is obtained by
\begin{eqnarray*}
 (N^*,M^*,h^*,\kappa^*)=\argmin_{(N,M,h,\kappa)} W_{dl}\;\;\;\hbox{subject to}\;\;\;
 \left\{
 \begin{array}{ll}
\frac{C_{dl,1}}{N} + \frac{C_{dl,2}}{NM} \leq \left(\kappa \hbox{TOL}/C_\alpha\right)^2\\\\
 C_{dl,3} h^{\eta} + \frac{C_{dl,4}}{M} \leq (1-\kappa) \hbox{TOL},
 \end{array}
 \right.
\end{eqnarray*}
where superscript $*$ is used to denote the optimal solution for the method parameters. The solution to the above problem using the Pontryagin's principle of minimization is given by
\begin{eqnarray*}
 \frac{1}{2 C_{dl,1}} \left( 1 + \frac{\gamma}{2\eta}\right)^{2} {\kappa^*}^2 \hbox{TOL} - \left[\frac{1}{2} + \left(1- \frac{1}{C_{dl,1}} \hbox{TOL} \right)\left( 1 + \frac{\gamma}{2\eta}\right)\right] \kappa^* + \left[1 + \frac{1}{2C_{dl,1}} \hbox{TOL}\right] = 0,
\end{eqnarray*}
for $\kappa^* \in ]0,1[$, and
\begin{eqnarray*}
 N^* &=& \frac{C_\alpha^2}{2\kappa^*} \frac{C_{dl,1}}{1 - \kappa^*\left(1+\frac{\gamma}{2\eta}\right)} \hbox{TOL}^{-2},\\
 M^* &=& \frac{C_{dl,2}}{2\left[1 - \kappa^*\left(1+\frac{\gamma}{2\eta}\right)\right]} \hbox{TOL}^{-1},\\
 h^* &=& \left( \frac{\gamma}{\eta} \frac{\kappa^*}{2C_{dl,3}}\right)^{1/\eta}\hbox{TOL}^{1/\eta}.
\end{eqnarray*}
In practice, we take the ceilings of $N$ and $M$ to obtain the optimal parameter values, $N^*$ and $M^*$. The minimal average work is then given by
\begin{equation}
W_{dl}^* \propto \hbox{TOL}^{-\left(3+\frac{\gamma}{\eta}\right)}.
\end{equation}
We note that the optimal setting without discretization, which is not presented herein but is applied in the numerical section (Section \ref{sec:num}), is different from the one presented above with $\gamma/\eta=0$; however, the same rates are attained with respect to $\hbox{TOL}$. The constants $C_{dl}$ are independent of $h$, but for the numerical estimation we approximate these constants using $\bm{g}_h$.

It has been proposed that the same samples can be used for the inner and outer loops as a means to reduce the number of evaluations of $\bm{g}_h$ \cite{huan}; more specifically, the computational cost is reduced from $\mathcal{O}(NM)$ to $\mathcal{O}(N)$ whenever the cost is dominated by the forward problem solver. Generally, such an approach leads to an extra bias \cite{key52}.

\subsection{Arithmetic underflow}\label{sec:underflow}
Arithmetic underflow can cause the return value of $\hat{p}_M(\bm{Y}_n)$, which is the denominator of its logarithm term in the DLMC estimator \eqref{eig_mc}, to be numerically zero. This occurs if all of the likelihoods, $p(\bm{Y}_n \vert \bm{\tilde{\theta}}_{n,m})$, $m=1,2,\ldots,M$, return values below the machine precision. For example, this is likely to happen if the prior distribution $\pi$ is not sufficiently concentrated on the posterior of $\bm{\theta}_n$. We also note that the posterior distribution itself becomes more concentrated for lower measurement noise, and the likelihood decreases exponentially with the size of $\bm{Y}$ ($=qN_e$). This can be seen from looking at the following leading-order terms of the log-likelihood:

\begin{eqnarray}\label{eq:like-inner1}
  \log{\left(p(\bm{Y}_n \vert \bm{\tilde{\theta}}_{n,m})\right)} &=&  -\frac{N_e}{2}\log{\left(2\pi \vert \bm{\Sigma_\epsilon} \vert \right)}-\frac{N_e}{2}\left\|\bm{g}(\bm{\theta}_{n}) - \bm{g}(\bm{\tilde{\theta}}_{n,m})\right\|_{\bm{\Sigma_\epsilon}^{-1}}^2 \nonumber \\
 && - \frac{1}{2} \bm{v_\epsilon}^{T}\left(\bm{g}(\bm{\theta}_n) - \bm{g}(\bm{\tilde{\theta}}_{n,m})\right)\mathcal{O}_{\mathbb{P}}\left(\sqrt{N_e}\right) - \frac{1}{2}\sum_{j=1}^{q} \sigma_{\epsilon_{j}}^{-1}\left(N_e + \mathcal{O}_{\mathbb{P}}\left(\sqrt{N_e}\right) \right), \nonumber \\
\end{eqnarray}
where $\bm{v_\epsilon}$ is the vector of the diagonal elements of $\bm{\Sigma_{\epsilon}}^{-1/2}$. The derivation of the above expression is given in \ref{proof_underflow}. We remark that the term $ -\frac{N_e}{2}\log{\left(2\pi \vert \bm{\Sigma_\epsilon} \vert \right)}$ is canceled out by the ratio in the DLMC estimator. The notation $X_{M}=\mathcal{O}_\mathbb{P}(a_{M})$ for a sequence of random variables $X_{M}$ with corresponding constants $a_{M}$, indexed by $M$, is defined as follows: for any $\epsilon>0$, there exists a finite $K>0$ and a finite $M_0>0$, such that for all $M \ge M_0$, $\mathbb{P}\left(\lvert X_{M} \rvert > K \lvert a_{M} \rvert \right) < \epsilon$.

\begin{figure}[H]
\centering
\includegraphics[width=.5\textwidth] {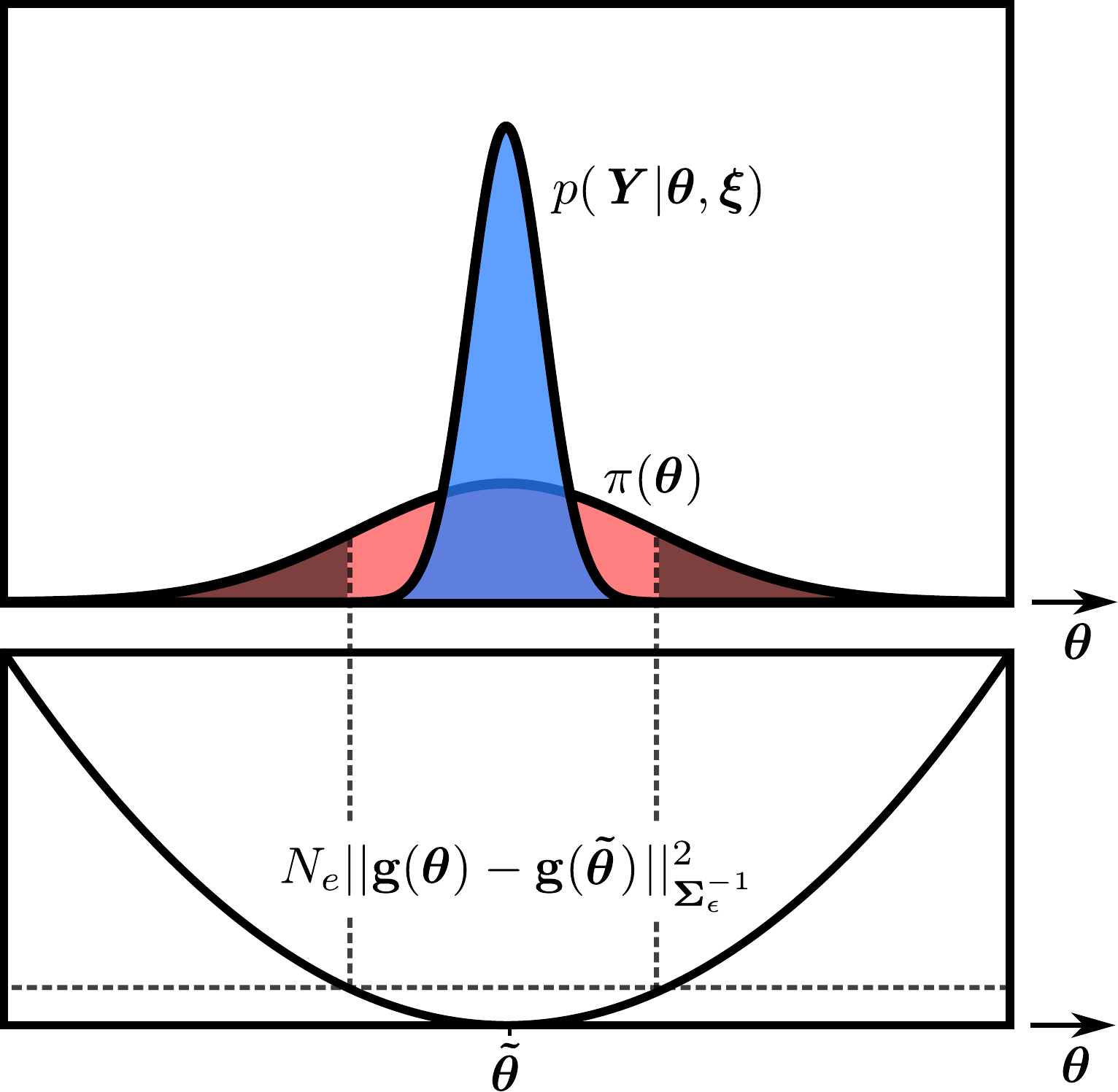} 
\caption{Illustration of the underflow that can occur when sampling the inner integral using the conventional DLMC.}
\label{fig:underflow}
\end{figure}

The top of Figure \ref{fig:underflow} presents a normal prior $\pi(\bm{\theta})$ (depicted in red) and a particular likelihood $p(\bm{Y} \vert \bm{\theta},\bm{\xi})$ (depicted in blue), both centered at $\boldsymbol{\tilde{\theta}}$. The shaded region over $\pi(\bm{\theta})$ represents the non-overlapping region between the two distributions. The lower illustration in Figure \ref{fig:underflow} presents the distance, multiplied by $N_e$, between the model $\bm{g}$ evaluated at $\bm{\theta}$ and $\bm{\tilde{\theta}}$.

\section{Monte Carlo with the Laplace method}\label{sec:la}
In this section, we consider MCLA for the estimation of the expected information gain \cite{key15}. MCLA approximates the inner part of the nested integral analytically as a normal approximation of the posterior pdf. Thus, the estimation of the expected information gain \eqref{expecinfgain} is reduced to a single integral over the prior space $\Theta$, typically at the price of some additional bias. The size of the bias decreases with the increasing number of measurements, i.e., the number of independent and repetitive experiments, $N_e$. MCLA has been extended previously to the approximation of posterior pdfs characterized by a non-informative manifold \cite{key52}.

In this section, we first show the main steps of the construction of MCLA, which can also be found elsewhere \cite{key15,key51}, and then present our optimized MCLA estimator.

\subsection{Laplace approximation}
The posterior pdf of $\bm{\theta}$ is given by Bayes' rule,
\begin{eqnarray}\label{eq31}
\pi(\bm{\theta} \vert \bm{Y}) =  \frac{\prod_{i=1}^{N_e}{\exp \left(-\frac{1}{2} \bm{r}_i^T(\bm{\theta}) \bm{\Sigma_\epsilon}^{-1} \bm{r}_i(\bm{\theta}) \right)} \pi(\bm{\theta})}{p(\bm{Y})}, 
\end{eqnarray}
where $\bm{r}_i(\bm{\theta}) = \bm{g}(\bm{\theta}_t) + \bm{\epsilon}_i - \bm{g}(\bm{\theta})$ is the residual of the $i$-th experimental data.

The Gaussian approximation of the posterior pdf \eqref{eq31} can be written as
\begin{equation}\label{eq32}
\pi_{g}(\bm{\theta} \vert \bm{Y}) = (2 \pi)^{-\frac{d}{2}} \vert \bm{\hat{\Sigma}} \vert ^{-\frac{1}{2}} \exp\left(-\frac{1}{2} \|\bm{\theta} - \bm{\hat{\theta}}\|^{2}_{\bm{\hat{\Sigma}}^{-1}}\right),   
\end{equation}
where $\bm{\hat{\theta}}$ is the MAP estimate, i.e.,
\begin{eqnarray}
  \bm{\hat{\theta}} \myeq \underset{\bm{\theta} \in \Theta}{\arg\min} \left[ \sum_{i=1}^{N_e} \left\| \bm{y}_i - \bm{g}(\bm{\theta}) \right\|^2_{\bm{\Sigma_{\epsilon}}^{-1}} + h(\bm{\theta}) \right], \text{and}
\end{eqnarray}
\begin{eqnarray}\label{eq_sigma}
\bm{\hat{\Sigma}}^{-1} = N_e \bm{J}(\bm{\hat{\theta}})^T\bm{\Sigma_\epsilon}^{-1} \bm{J}(\bm{\hat{\theta}}) -  \nabla_{\bm{\theta}}  \nabla_{\bm{\theta}}  h(\bm{\hat{\theta}}) +  \mathcal{O}_\mathbb{P}\left(\sqrt{N_e}\right)
\end{eqnarray}
is the inverse Hessian matrix of the negative logarithm of the posterior pdf evaluated at $\bm{\hat{\theta}}$, and $h(\bm{\theta}) = \log (\pi(\bm{\theta}))$. It has been shown \cite{key15} that 
\begin{eqnarray}\label{eq_theta}
\bm{\hat{\theta}} = \bm{\theta}_t - \left(N_e \bm{J}^T(\bm{\theta}_t) \bm{\Sigma_\epsilon}^{-1} \bm{J}(\bm{\theta}_t) + \bm{H}^T \bm{\Sigma_\epsilon}^{-1} \bm{E_\epsilon} - \nabla \nabla h(\bm{\theta}_t)\right)^{-1} \bm{J}^T \bm{\Sigma_\epsilon}^{-1} \bm{E_\epsilon} + \mathcal{O}_\mathbb{P}\left(\frac{1}{N_e}\right)\,,
\end{eqnarray} 
where
\begin{eqnarray*}
  \bm{E_\epsilon}(\bm{\theta})  = \sum_{i=1}^{N_e} \bm{r}_i^T(\bm{\theta}), \;\;\;
 \bm{J}(\bm{\theta}) =- \nabla_{\bm{\theta}} \bm{g}(\bm{\theta})\;\;\;\hbox{and}\;\;\; 
  \bm{H}(\bm{\theta}) =-\nabla_{\bm{\theta}}  \nabla_{\bm{\theta}} \bm{g}(\bm{\theta}).
\end{eqnarray*}

For a sufficiently large $N_e$, we have
\begin{eqnarray}
\bm{\hat{\theta}} &= \bm{\theta}_t  + \mathcal{O}_\mathbb{P}\left(\frac{1}{\sqrt{N_e}}\right),
\end{eqnarray}
which is the approximation used in the MCLA estimator, i.e., $\bm{\hat{\theta}} \approx \bm{\theta}_t$.

\subsection{Monte Carlo with the Laplace approximation (MCLA) estimator}
The Gaussian approximation \eqref{eq32} with $\bm{\hat{\theta}}$ and $\bm{\hat{\Sigma}}$, given by \eqref{eq_theta} and \eqref{eq_sigma}, respectively, leads to an analytical expression of the Kullback-Leibler divergence, which subsequently yields the following approximation of the expected information gain:
\begin{eqnarray}
\label{eq38}
I = \int_{\Theta}{\left[-\frac{1}{2} \log((2 \pi)^{d} \vert \bm{\hat{\Sigma}}(\hat{\bm{\theta}}) \vert) -\frac{d}{2} - h(\hat{\bm{{\theta}}}) \right] p(\bm{\theta}) d\bm{\theta}} + \mathcal{O}\left(\frac{1}{N_e}\right).
\end{eqnarray}
The detailed proof of (\ref{eq38}) is given in \ref{ap:secB}. Then we can use MC sampling on the integral in \eqref{eq38}. The additional bias (Laplace error) is of the order $\mathcal{O}\left(\frac{1}{N_e}\right)$. Many other schemes could also be applied, e.g., Gaussian quadrature \cite{key15,key51}. The covariance matrix is proportional to $N_e^{-1}$, which can be seen by expanding \eqref{eq_sigma} using the Sherman-Morrison formula; hence, $\vert \hat{\vec\Sigma} \vert = \mathcal{O}\left(N_e^{-d}\right)$ and $I = \mathcal{O}\left(\log(N_e)\right)$, where $d = \dim(\bm{\theta})$.

The MCLA estimator of $I$ is defined as
\begin{eqnarray}
 \mathcal{I}_{_{la}} \myeq \frac{1}{N} \sum_{n=1}^{N}  \left(-\frac{1}{2} \log((2 \pi)^{d} \vert \bm{\hat{\Sigma}}(\bm{\theta}_n) \vert) -\frac{d}{2} - h(\bm{{\theta}_n}) \right),
\end{eqnarray} 
where $N$ is the number of MC samples. The average computational work of MCLA can be estimated by
$$W_{la} \propto NN_{jac} h^{-\gamma},$$ 
where $N_{_{jac}}$ is the number of forward model evaluations required in the Jacobian matrix $\bm{J}(\bm{\theta})$ in \eqref{eq_theta} and \eqref{eq_sigma}. For the finite difference,
\begin{eqnarray*}
N_{_{jac}} = \left\{
\begin{array}{lll}
d + 1\;\;\;\;\hbox{if forward or backward differentiation,}\\\\
2 d \;\;\;\;\;\; \;\; \hbox{if central differentiation.}
\end{array}
\right.
\end{eqnarray*}

To derive the optimal setting for the MCLA estimator $\mathcal{I}_{la}$, we use Proposition \ref{prop2}.
\begin{proposition}
\label{prop2}
The bias and variance of the MCLA estimator $\mathcal{I}_{la}$ can be estimated by:
\begin{eqnarray}
 \lvert I - \mathbb{E}\left[\mathcal{I}_{la}\right] \rvert&\leq&\frac{C_{la,2}}{N_e} + C_{la,3} h^{\eta}+o(h^{\eta}), \\
\nonumber\\
 \mathbb{V}\left[ \mathcal{I}_{la} \right] &=&\frac{C_{la,1}}{N},
\end{eqnarray}
where $C_{la,1}$, $C_{la,2}$ and $C_{la,3}$ are constants.
\end{proposition}
The total bias is composed of the bias $C_{la,2}/N_e$ introduced by the Laplace approximation in \eqref{eq38}, and the bias $C_{la,3} h^{\eta}$ from the numerical discretization. The constant $C_{la,1}$ is given by $C_{la,1}=\mathbb{V}[D_{kl}]$, since
\begin{eqnarray}
 \mathbb{V}\left[\mathcal{I}_{_{la}}\right] = \mathbb{V}\left[ \frac{1}{N} \sum_{n=1}^{N} D_{kl}(\bm{\theta}_n)\right] = \frac{\mathbb{V}[D_{kl}]}{N}.
\end{eqnarray}

\subsection{Optimal setting for the MCLA estimator}
The cost minimization problem for the MCLA estimator of a given $\hbox{TOL}>0$ is
\begin{eqnarray}
 (N^*,h^*,\kappa^*)=\argmin_{(N,h,\kappa)} W_{la}\;\;\;\hbox{subject to}\;\;\;
 \left\{
\begin{array}{lll}
\frac{C_{la,1}}{N} \leq \left(\kappa \hbox{TOL}/C_\alpha\right)^2,\\\\
\frac{C_{la,2}}{N_e} + C_{la,3}h^{\eta}  \leq (1-\kappa) \hbox{TOL}.
\end{array}
\right.
\end{eqnarray}

We note that the Laplace bias, $C_{la,2}/N_e$, is related to the number of repetitive experiments, $N_e$, and is therefore not a method parameter. Also, there exists no solution to the above problem if the constraint,
\begin{eqnarray}\label{constraint}
(1-\kappa) \hbox{TOL} \geq \frac{C_{la,2}}{N_e},
\end{eqnarray}
does not hold. Thus, for $\hbox{TOL} \geq \frac{C_{la,2}}{N_e}$, the solution to the Pontryagin's principle of minimization is given by
\begin{eqnarray*}
\kappa^* &=& \frac{1 - \frac{ C_{la,2}}{N_e} \hbox{TOL}^{-1} }{1 + \frac{ \gamma}{2\eta}},\\
h^* &=& \left( \frac{\gamma}{\eta} \frac{\kappa^*}{2 C_{la,3}}\right)^{1/\eta} \hbox{TOL} ^{1/\eta},\\
N^* &=& \frac{2\eta}{\gamma}  C_{la,1} \frac{1-\kappa^*}{{\kappa^*}^3} C_\alpha \hbox{TOL}^{-2} - \frac{2\eta}{\gamma} \frac{ C_{la,1}  C_{la,2}}{N_e{\kappa^*}^3} C_\alpha^2 \hbox{TOL}^{-3},
\end{eqnarray*}
with an average work of order $\mathcal{O}\left(\hbox{TOL}^{-\left(2+\frac{\gamma}{\eta}\right)}\right)$. This means that the work rate of MCLA is better than that of DLMC, although an accuracy lower than $C_{la,2}/N_e$ cannot be achieved.

\section{Double-loop Monte Carlo with Laplace-based importance sampling}\label{sec:DLMCIS}
In this section, we devise a new approach that entails a measure change based on the Laplace approximation for the inner loop of DLMC. The approach does not introduce any extra bias, avoids the occurence of underflow, and substantially reduces the number of samples of the inner loop. In other words, we use the Laplace approximation to achieve a more efficient sampling on the posterior of $\bm{\theta}_n$.

\subsection{Double-loop Monte Carlo with importance sampling (DLMCIS) estimator}
To compute the evidence term in \eqref{eig_mc}, we change the sampling pdf, $\pi$, to a new pdf, $\tilde{\pi}_n$, given by $\tilde{\pi}_n(\bm{\theta}) \sim \mathcal{N} \left(\bm{\hat{\theta}}_n, \bm{\hat{\Sigma}}(\bm{\hat{\theta}}_n)\right)$, where 
\begin{eqnarray}
 \label{fittheta}
 \bm{\hat{\theta}}_n  &\myeq & \underset{\bm{\theta} \in \Theta}{\arg\min} \left[ \frac{1}{2} \sum_{i=1}^{N_e} \left\| \bm{y}_{n,i} - \bm{g}_h(\bm{\theta}) \right\|^2_{\bm{\Sigma_\epsilon}^{-1}} + \log (\pi(\bm{\theta})) \right],
\end{eqnarray}
and
\begin{eqnarray}
 \bm{\hat{\Sigma}} (\bm{\hat{\theta}}_n)&=& \Big( N_e \bm{J}(\bm{\hat{\theta}}_n)^T\bm{\Sigma_\epsilon}^{-1} \bm{J}(\bm{\hat{\theta}}_n)-\nabla_{\bm{\theta}}  \nabla_{\bm{\theta}}  h(\bm{\hat{\theta}}_n) \Big)^{-1} +  \mathcal{O}_\mathbb{P}\left(\frac{1}{\sqrt{N_e}}\right).
\end{eqnarray}
The importance sampling based on the Laplace approximation leads to the DLMCIS estimator,
\begin{eqnarray}
\label{eq:dlis}
\mathcal{I}_{dlis} \myeq \frac{1}{N} \sum_{n=1}^{N} \log \left( \frac{p(\bm{Y}_n \vert \bm{\theta}_n)}{\frac{1}{M} \sum _{m=1} ^{M}L(\bm{Y}_n,\bm{\tilde{\theta}}_{n,m})}\right),
\end{eqnarray}
where
\begin{eqnarray}
L(\bm{Y}_n,\bm{\tilde{\theta}}_{n,m}) =  \dfrac{p (\bm{Y}_n \vert \bm{\tilde{\theta}}_{n,m}) \pi(\bm{\tilde{\theta}}_{n,m})}{\tilde{\pi}_n(\bm{\tilde{\theta}}_{n,m})}.
\end{eqnarray}  
Here $\bm{\theta}_n$ and $\bm{\tilde{\theta}}_{n,m}$ are sampled from the prior $\pi$ and the Laplace-based pdf $\tilde{\pi}_n$ respectively. 

\begin{figure}[ht!]
\centering
\subfloat[Uniform prior with $\left(\bm{\theta}_n,  \bm{\hat{\Sigma}} (\bm{\theta} _n)\right)$]{
\includegraphics[width=.475\textwidth] {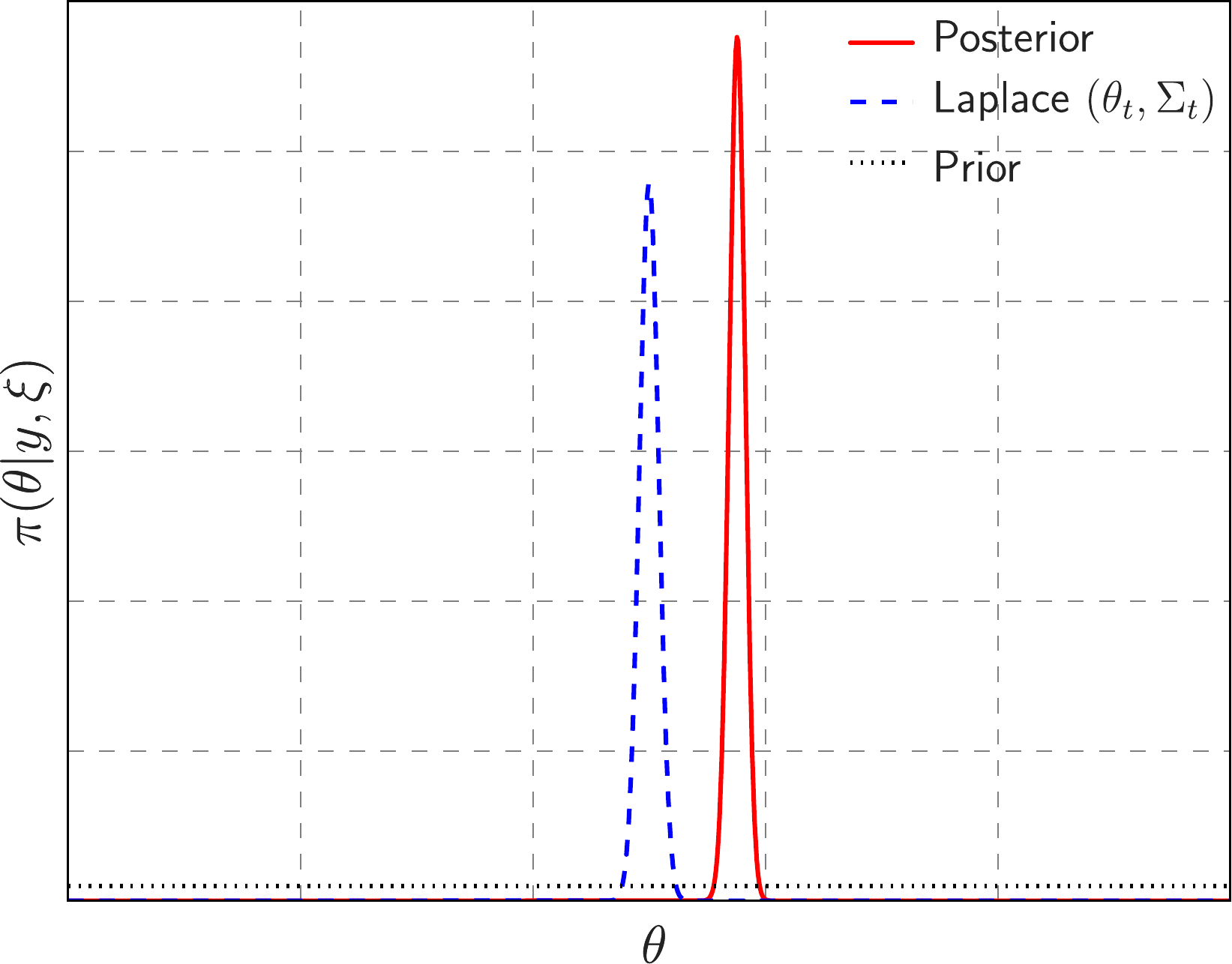}}
\hfill
\subfloat[Normal prior with $\left(\bm{\theta}_n,  \bm{\hat{\Sigma}} (\bm{\theta} _n)\right)$]{
\includegraphics[width=.475\textwidth] {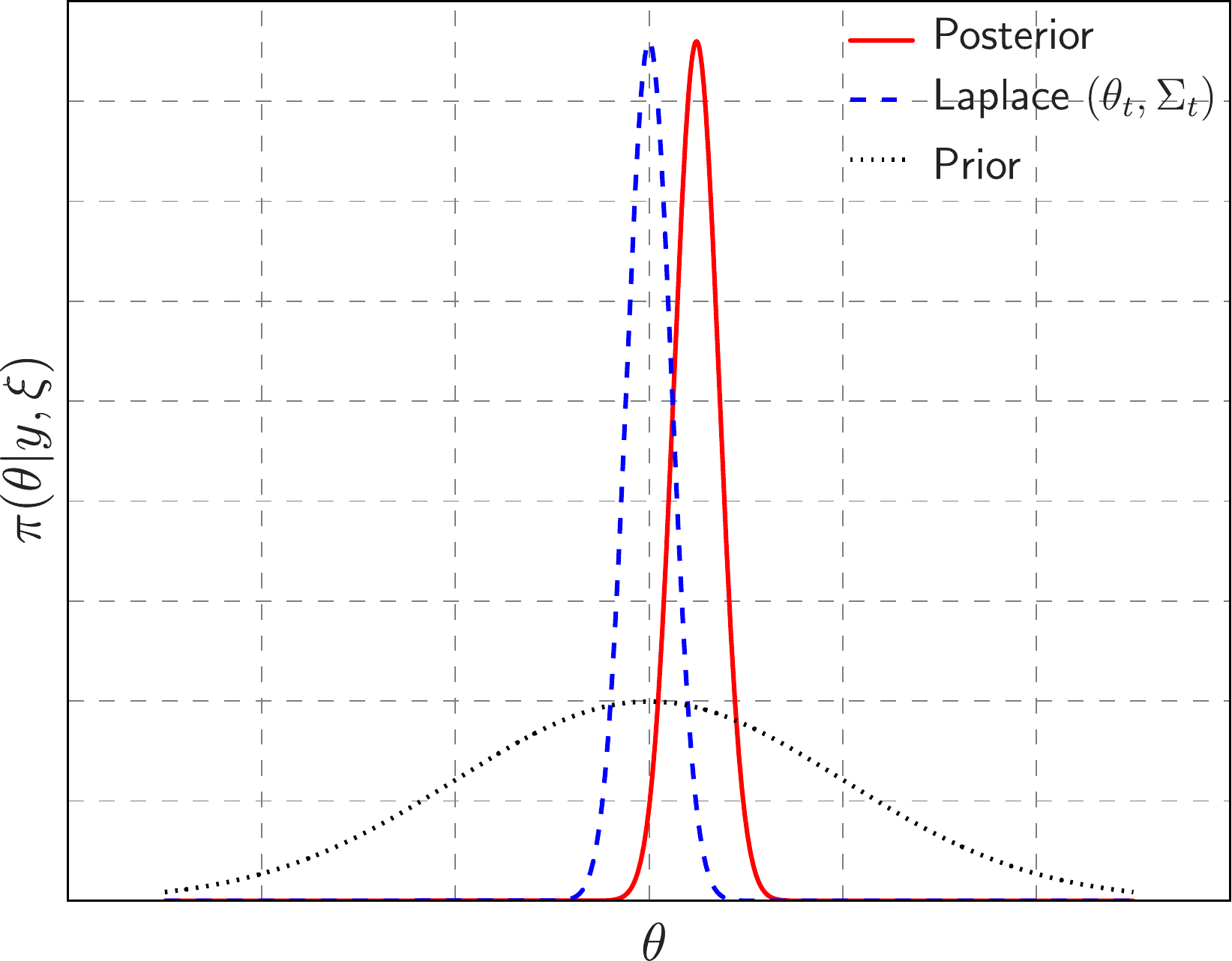} }%
\hfill

\subfloat[Uniform prior with $\left(\hat{\bm{\theta}}_n,  \bm{\hat{\Sigma}} (\hat{\bm{\theta}} _n)\right)$]{
\includegraphics[width=.475\textwidth] {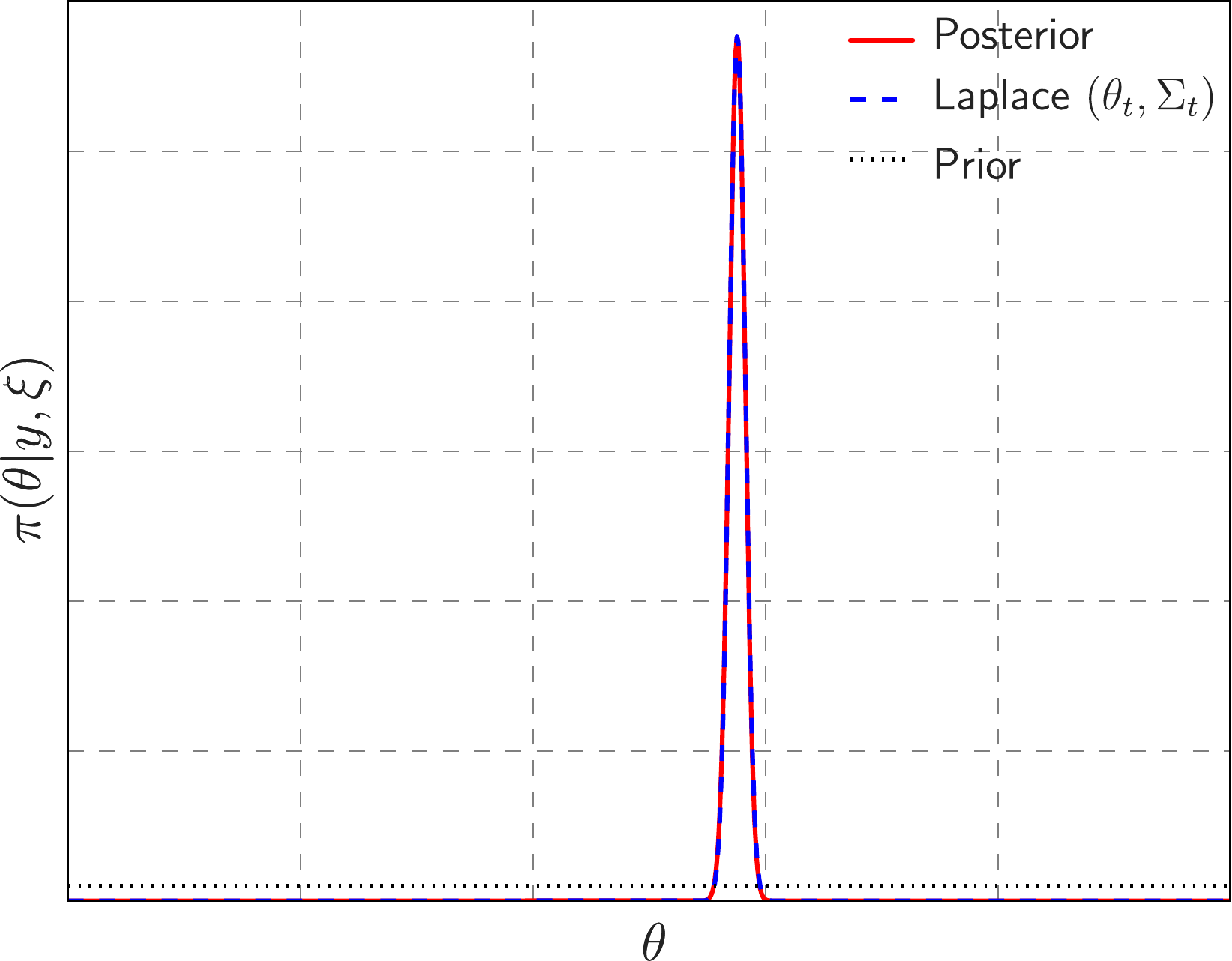}   }
\hfill
\subfloat[Normal prior with $\left(\hat{\bm{\theta}}_n,  \bm{\hat{\Sigma}} (\hat{\bm{\theta}} _n)\right)$]{
\includegraphics[width=.475\textwidth] {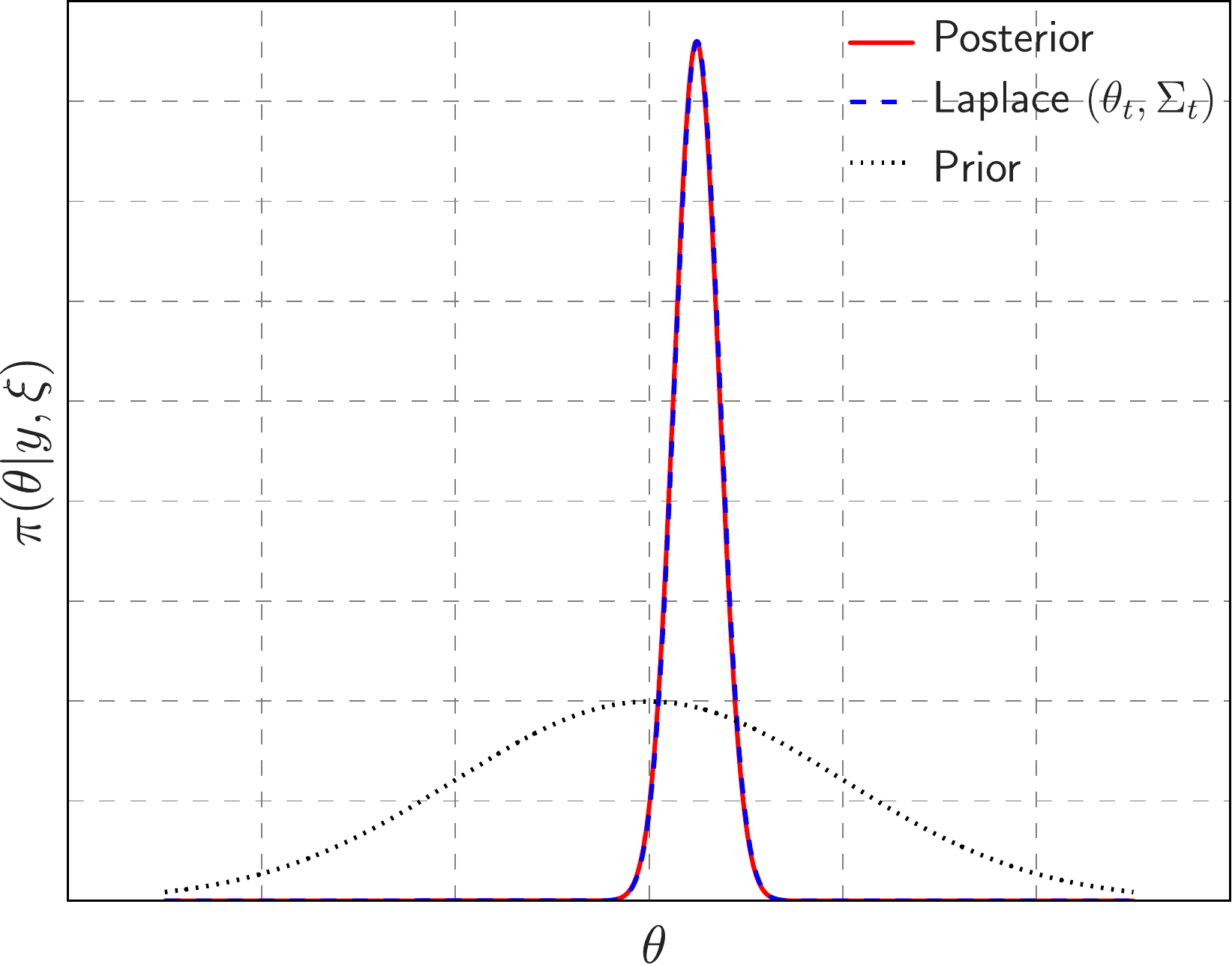}
\label{meanplot} }%
\hfill                        
\caption{Posterior pdfs by DLMCIS: top with $\bm{\theta}_n$ and bottom with $\bm{\hat{\theta}}_n$.}
\label{fig:sampling}   
\end{figure}

It would be computationally less costly to construct the Laplace approximation around $\bm{\theta}_n$ rather than $\bm{\hat{\theta}}_n$. Unfortunately, we have verified numerically that it would be inefficient as the discrepancy between $\bm{\theta}_n$ and $\bm{\hat{\theta}}_n$ can be too large, and thus not enough concentration on the posterior; see Figure \ref{fig:sampling}.

\subsection{Optimal setting for the DLMCIS estimator}
The average computational work of the DLMCIS estimator \eqref{eq:dlis} is assumed to be
\begin{equation}
W_{dlis} = C_1NMh^{-\gamma} + C_2 N N_{jac} h^{-\gamma},
\end{equation}
with constants $C_1,C_2>0$. The cost minimization problem can be written as
\begin{eqnarray*}
(N^*,M^*,h^*,\kappa^*)=\argmin_{(N,M,h,\kappa)} W_{dlis}\;\;\;\hbox{subject to}\;\;\;
\left\{
\begin{array}{ll}
\frac{C_{dlis,1}}{N} + \frac{C_{dlis,2}}{NM} \leq  \left(\kappa \hbox{TOL}/C_\alpha\right)^2,\\\\
C_{dlis,3} h^{\eta} + \frac{C_{dlis,4}}{M} \leq (1-\kappa) \hbox{TOL},
\end{array}
\right.
\end{eqnarray*}
with constants $C_{dlis,k}$, $k=1,\cdots,4$. The optimal $\kappa$, denoted by $\kappa^*$, is the solution of a second order polynomial equation that belongs to $]0,1[$, and the optimal $N$, $M$, and $h$ are
\begin{eqnarray*}
 N^* &=& \frac{C_{dlis,1} C_\alpha^2}{{\kappa^*}^2} \hbox{TOL}^{-2} + \left(1 - \kappa^* \left(1+ \frac{\gamma}{2\eta}\right)\right) \frac{C_{dlis,1} C_\alpha^2}{{\kappa^*}^2} \hbox{TOL}^{-1}, \\
 M^* &=& \frac{C_{dlis,4}}{1 - \kappa^*\left(1+ \frac{\gamma}{2\eta}\right)} \hbox{TOL}^{-1},\\
 h^* &=& \left( \frac{\gamma}{\eta} \frac{\kappa^*}{2C_{dlis,3}}\right)^{1/\eta}\hbox{TOL}^{1/\eta},\text{ respectively}.
\end{eqnarray*}
The average computational work is $W_{dlis}^* \propto TOL^{-\left(3+\frac{\gamma}{\eta}\right)}$, i.e., the method preserves the order as DLMC, but the constants $C_{dlis,1}$, $C_{dlis,2}$, and $C_{dlis,3}$, are typically several magnitudes smaller than their corresponding counterparts $C_{dl}$, due to the variance reduction achieved by the importance sampling; $C_{dlis,3}$ is still equal to $C_{dl,3}$.

\subsection{Effect of the change of measure on the arithmetic underflow}
The change of measure mitigates the risk of underflow in DLMC as discussed in Section \ref{sec:underflow}. The advantage of changing measure in this regard can be seen by observing that the constants in front of the leading-order terms of the following expansion of the likelihood become much smaller:
\begin{eqnarray}
 \log\left(L(\bm{Y}_n;\bm{\tilde{\theta}}_{n,m})\right) &=&  \frac{N_e}{2}\left(-\log{\left(2\pi \vert \bm{\Sigma_\epsilon} \vert\right)}+\log{\left(2\pi \vert \bm{\hat{\Sigma}}(\bm{\hat{\theta}}_n) \vert\right)}\right)-\frac{N_e}{2}\left\|\bm{g}(\bm{\theta}_n) - \bm{g}(\bm{\tilde{\theta}}_{n,m})\right\|^2_{\bm{\Sigma_\epsilon}^{-1}}
\nonumber \\
 && +\frac{N_e}{2} \left\|\bm{J}(\bm{\hat{\theta}}_n)\left(\bm{\tilde{\theta}}_{n,m} - \bm{\hat{\theta}}_n\right)\right\|^2_{ \bm{\Sigma_\epsilon}^{-1}} 
   - \frac{1}{2} \left\| \bm{\tilde{\theta}}_{n,m} - \bm{\hat{\theta}}_n \right\|^2_{ \nabla_{\bm{\theta}} \nabla_{\bm{\theta}} h(\bm{\hat{\theta}}_n)} \nonumber\\
  &&\left.  - \frac{1}{2}\bm{v_\epsilon}^{T}\left(\bm{g}(\bm{\theta}_n) - \bm{g}(\bm{\tilde{\theta}}_{n,m})\right)\mathcal{O}_{\mathbb{P}}\left(\sqrt{N_e}\right) - \frac{1}{2}\sum_{j=1}^{q} \sigma_{\epsilon_{j}}^{-1}\left(N_e + \mathcal{O}_{\mathbb{P}}\left(\sqrt{N_e}\right) \right) \right). \nonumber\\
\end{eqnarray}
The detailed proof can be found in \ref{correct_underflow}.

\section{Numerical examples}\label{sec:num}

We present three examples to demonstrate the numerical performance of our proposed method, DLMCIS, in comparison with DLMC and MCLA. In Example 1, we consider a scalar model whose parameter follows a normal prior distribution. In Example 2, we apply the methods to a nonlinear benchmark test, which was studied in \cite{huan}. In the third example, we design the placement of electrode sensors in order to maximize signal information about unknown parameters, in this case the orientation of fibers, in a laminate composite material during EIT experiments.

\subsection{Example 1: Linear scalar model}
In this example, we consider an algebraic linear model given by
\begin{equation}
y_i(\theta,\xi) = \theta \left( 1 + \xi \right)^2 + \epsilon_i, \quad \text{for} \, \ i=1,2,\cdots,N_e,
\end{equation}
with prior $\pi(\theta) \sim \mathcal{N}\left(1, 0.01\right)$ and observational noise $\epsilon_i \sim \mathcal{N}\left(0, (2+(\xi-10)/10)^2\right)$. Two repetitive experiments ($N_e=2$) are conducted. 

\subsubsection{Optimality and error convergence}
We study the problem in two parts. First, we analyze the optimality, that is, $N^*$, $M^*$, and $\kappa^*$ for the tolerance range $[10^{-5},1]$. Second, we analyze the error, computational time, and work in relation to $\hbox{TOL}$; due to computational constraints, we shrink the tolerance range to $[10^{-2},1]$.

We compute the expected information gain for a fixed experimental setup, $\xi=10$, at a given confidence level of $97.5\%$, i.e., $\alpha=0.05$. For DLMCIS, the constants $C_{dlis,1}$ and $C_{dlis,2}$ are estimated using DLMC with $N=M=100$; for MCLA, the constant $C_{la,1}$ is estimated by MC with $N=100$. Because the model is linear on $\theta$ and the prior is normal, the Laplace method is exact in this case and does not induce any bias. We estimate $C_{la,2}$ by the absolute difference between MCLA and DLMCIS in $D_{kl}$ for a small, fixed number of outer samples with identical realizations.

\begin{figure}[!htb]
  \centering
  \includegraphics[width=0.75\linewidth]{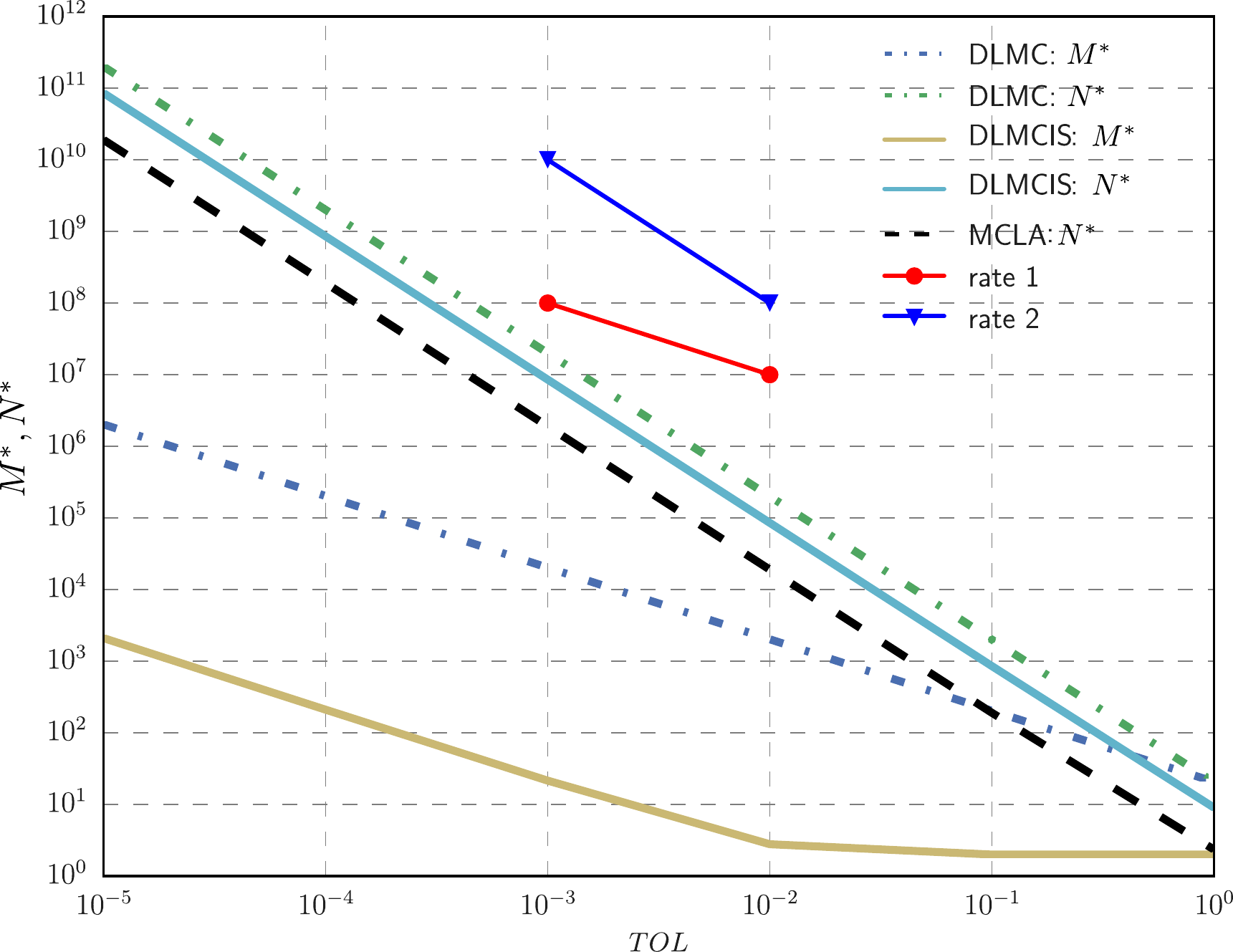}
  \caption{Optimal setting (outer $N^*$ and inner $M^*$ number of samples) vs. tolerance for linear scalar model with normal prior (Example 1).}
  \label{optimality_N}
\end{figure}

Figure \ref{optimality_N} shows $N^*$ and $M^*$ for the tolerance range. The optimal value $M^*$ for DLMCIS follows the asymptotic rate of the order $\mathcal{O}\left(\hbox{TOL}^{-1}\right)$, except for tolerances larger than $10^{-2}$, where it is small (less than 5). The number of samples needed with DLMCIS can be several magnitudes lower than the number required with DLMC, which shows the efficiency of the importance sampling with the Laplace approximation. The optimal value $\kappa^*$ is constant over the considered tolerance window: roughly 0.64, 0.67, and 1, for DLMC, DLMCIS, and MCLA, respectively. This means that the statistical error contributes more than the bias to the total error in this example, especially for MCLA.

\begin{figure}[!htb]
  \centering
  \subfloat[MCLA]{\label{errvstol_N.a}\includegraphics[width=0.5\linewidth]{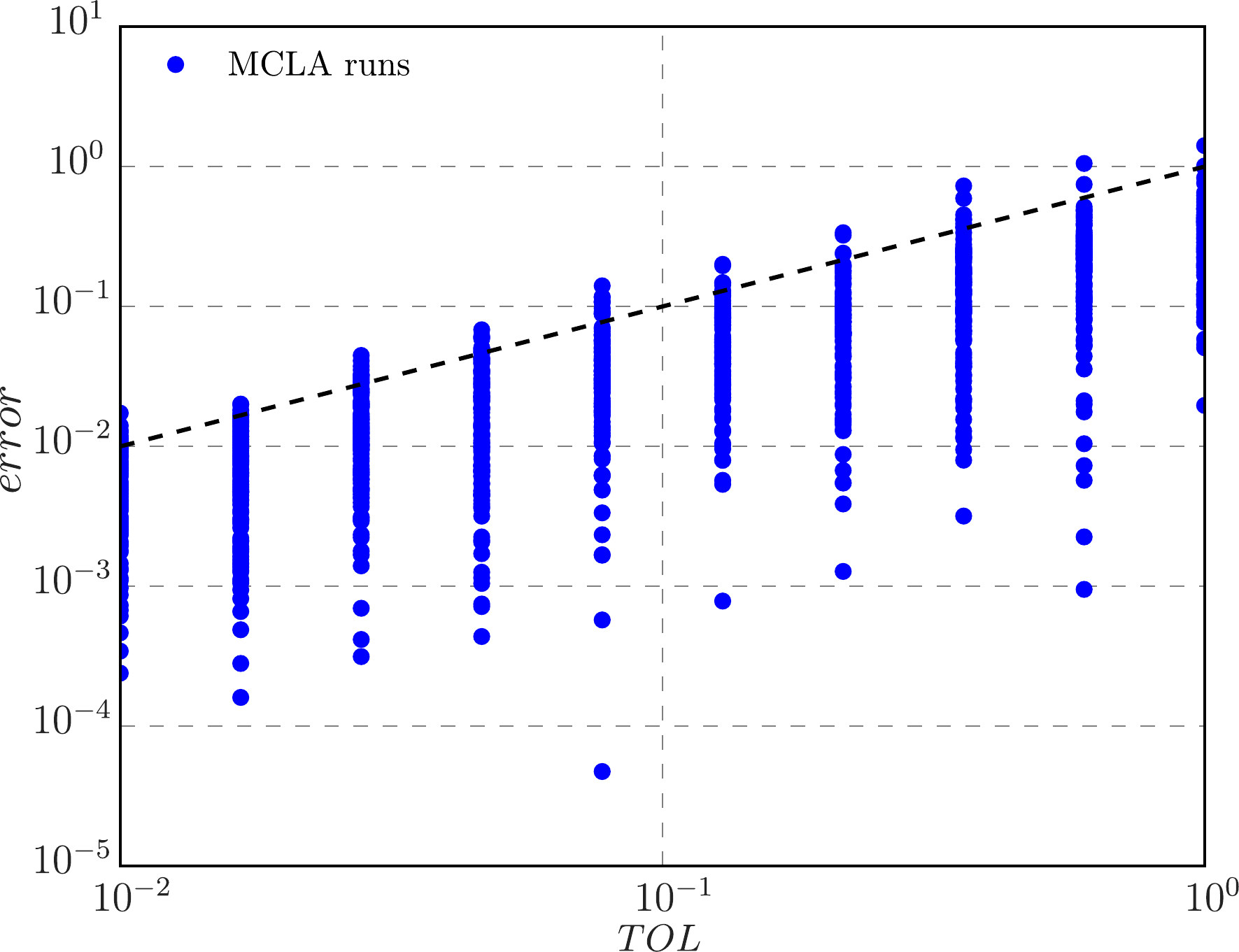}}
  \subfloat[DLMCIS]{\label{errvstol_N.b}\includegraphics[width=0.5\linewidth]{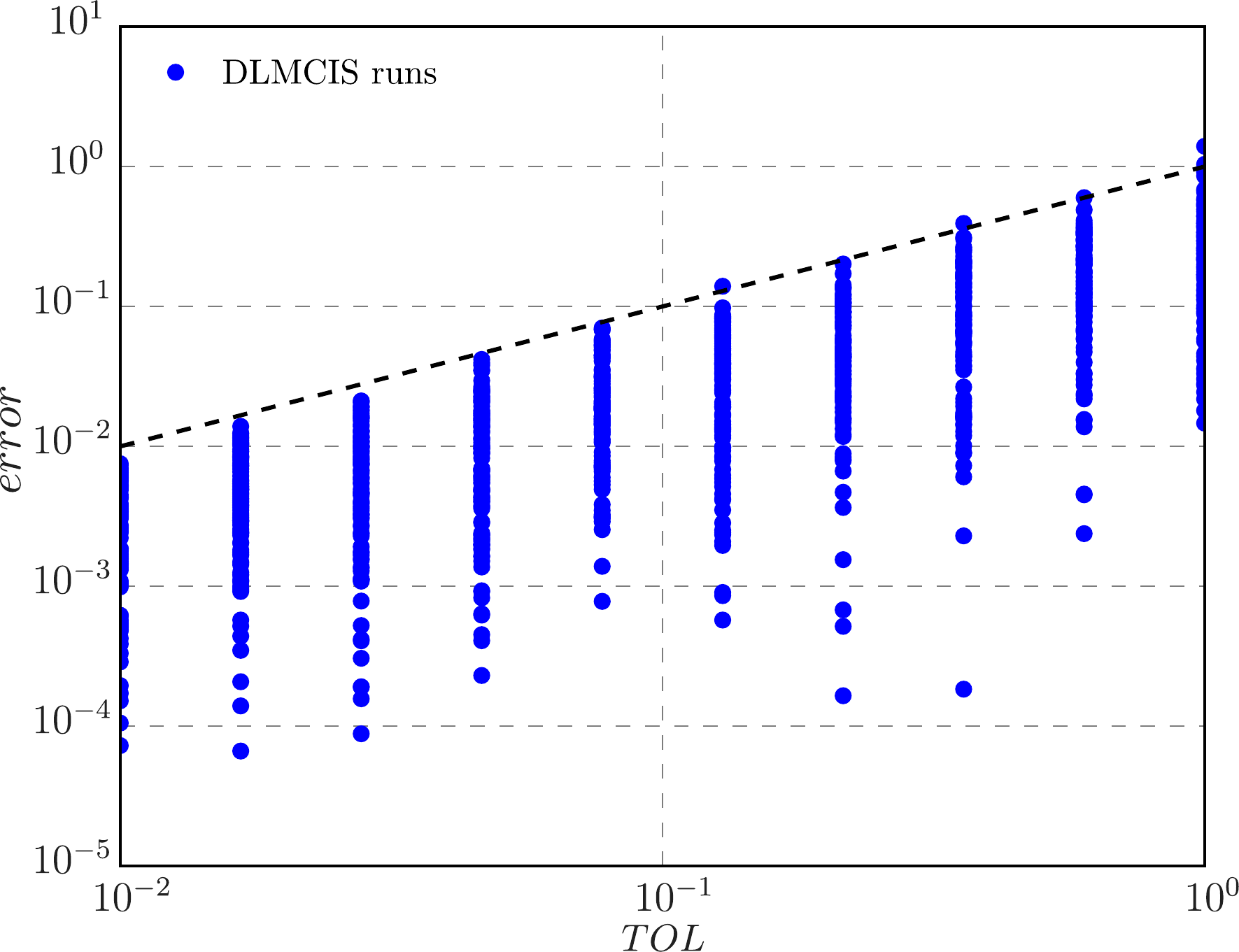}}
  \caption{Error vs. tolerance for linear model with normal prior (Example 1).}
  \label{errvstol_N}
\end{figure}

We now look at the accuracies achieved by MCLA and DLMCIS using their optimal settings with respect to the specified tolerance $\hbox{TOL}$. In Figure \ref{errvstol_N}, we provide a consistency test between the actual computed absolute error, $\vert \mathcal{I}-I \vert$, and $\hbox{TOL}$. The absolute error and $\hbox{TOL}$ are in agreement for both MCLA and DLMCIS, which numerically validates our bias estimation in the optimal setting derivation. The absolute error of MCLA is consistently a bit higher than the tolerance, but this is expected since the error is below $\hbox{TOL}$ with $97.5\%$ probability. Here DLMCIS is here a bit more conservative than expected.

\begin{figure}[!htb]
  \centering
  \subfloat[Work]{\label{time_work_N.b}\includegraphics[width=0.5\linewidth]{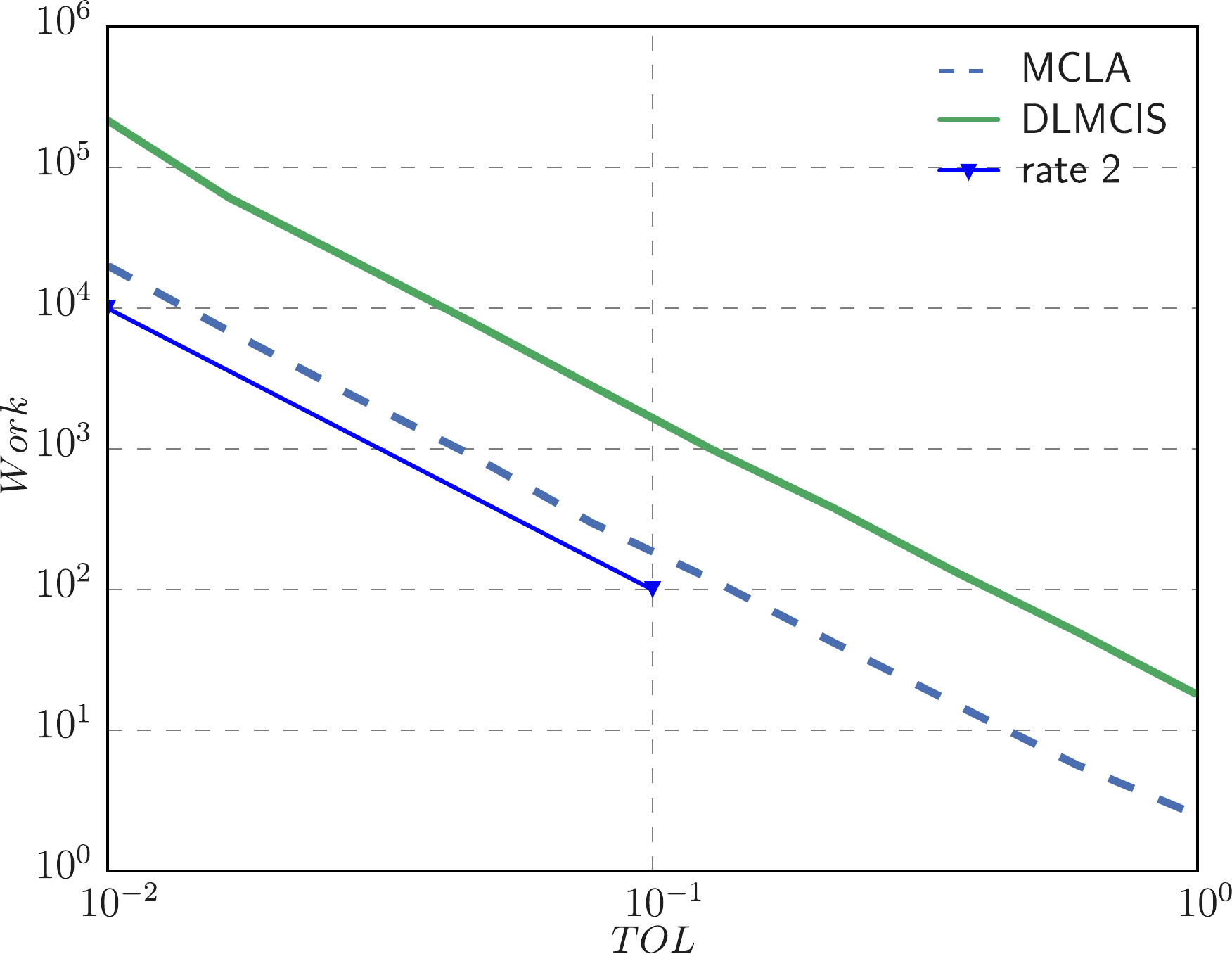}}
  \subfloat[Time]{\label{time_work_N.a}\includegraphics[width=0.5\linewidth]{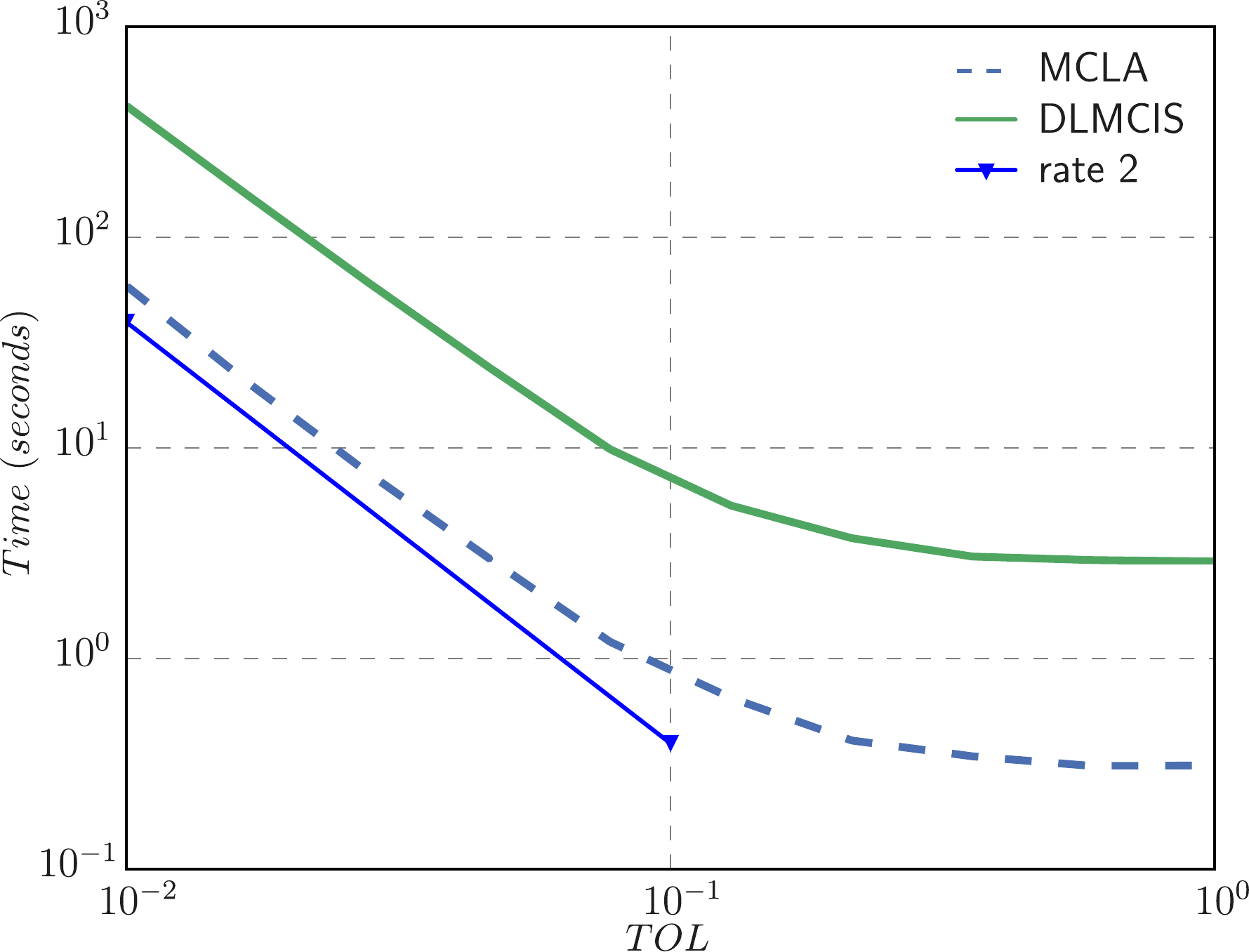}}
  \caption{Average computational work and running time vs. tolerance for linear scalar model with normal prior (Example 1).}
  \label{time_work_N}
\end{figure}

The average computational work required for MCLA is consistent with the theoretical asymptotic result $\mathcal{O}(\hbox{TOL}^{-2})$, whereas the work for DLMCIS, interestingly enough, has the same rate of $\mathcal{O}(\hbox{TOL}^{-2})$, which is better than the theoretical result, $\mathcal{O}(\hbox{TOL}^{-3})$, for this method; see Figures \ref{time_work_N.b} and \ref{time_work_N.a} for the running time and the average computational work. The DLMCIS work rate is a pre-asymptotic rate since the optimal value $M^*$ is unaffected over the tolerance range $[0.01,1]$ (see Figure \ref{optimality_N}). The computation of the inner loop using DLMCIS is very efficient in this example. We observe that MCLA has a lower error constant than DLMCIS; this is due to the fact that the Laplace approximation is exact owing to as the Gaussian response of the linear model with normal prior.

\subsubsection{Expected information gain}
\begin{figure}[!htp]
\centering
\includegraphics[width=0.75\linewidth]{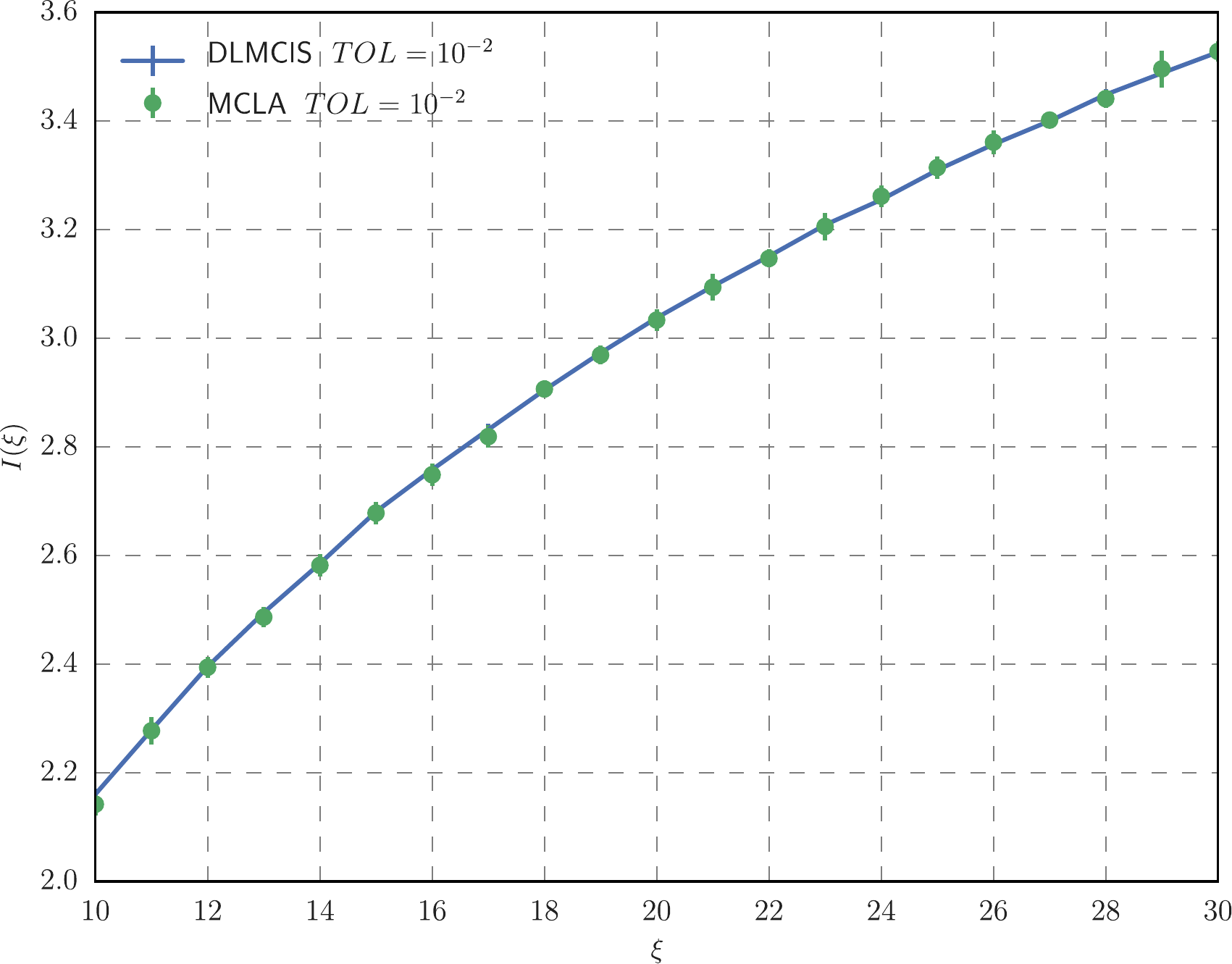}
\caption{Expected information gain for linear scalar model with normal prior (Example 1).}
\label{I_N}
\end{figure}

In this section, we show the expected information gain estimated by the methods at $\hbox{TOL}=0.01$ over a set of experiment setups, $\xi\in[10,30]$. As seen in Figure \ref{I_N}, DLMCIS and MCLA show agreement, though the variance of the MCLA estimator is larger than that of DLMCIS.

\subsection{Example 2: Nonlinear scalar model}
We consider the algebraic nonlinear model found in Huan and Marzouk \cite{huan}:
\begin{equation}
y_i = \theta^3 \xi^2 + \theta \exp \left( - |0.2 - \xi |\right) + \epsilon_i, \quad \text{for} \, \ i=1,2,\cdots,N_e,
\end{equation}
where $N_e=1$ and $N_e=10$. The prior is $\pi(\theta) \sim \mathcal{U}\left(0, 1\right)$, and the observational noise is $\epsilon_i \sim \mathcal{N}\left(0, 10^{-3}\right)$. 

\subsubsection{Optimality and error convergence}
Same as in Example 1, we provide a comparison between MCLA and DLMCIS; here for the experiment $\xi=1$. Example 2 is more challenging for MCLA than Example 1, since the use of the Laplace approximation leads to an inherent bias due to the uniform prior and the nonlinearity in the forward model with respect to the unknown parameter, $\theta$.

\begin{figure}[!htp]
  \centering
  \subfloat[$N_e=1$]{\label{optimality_MN_U.a}\includegraphics[width=0.5\linewidth]{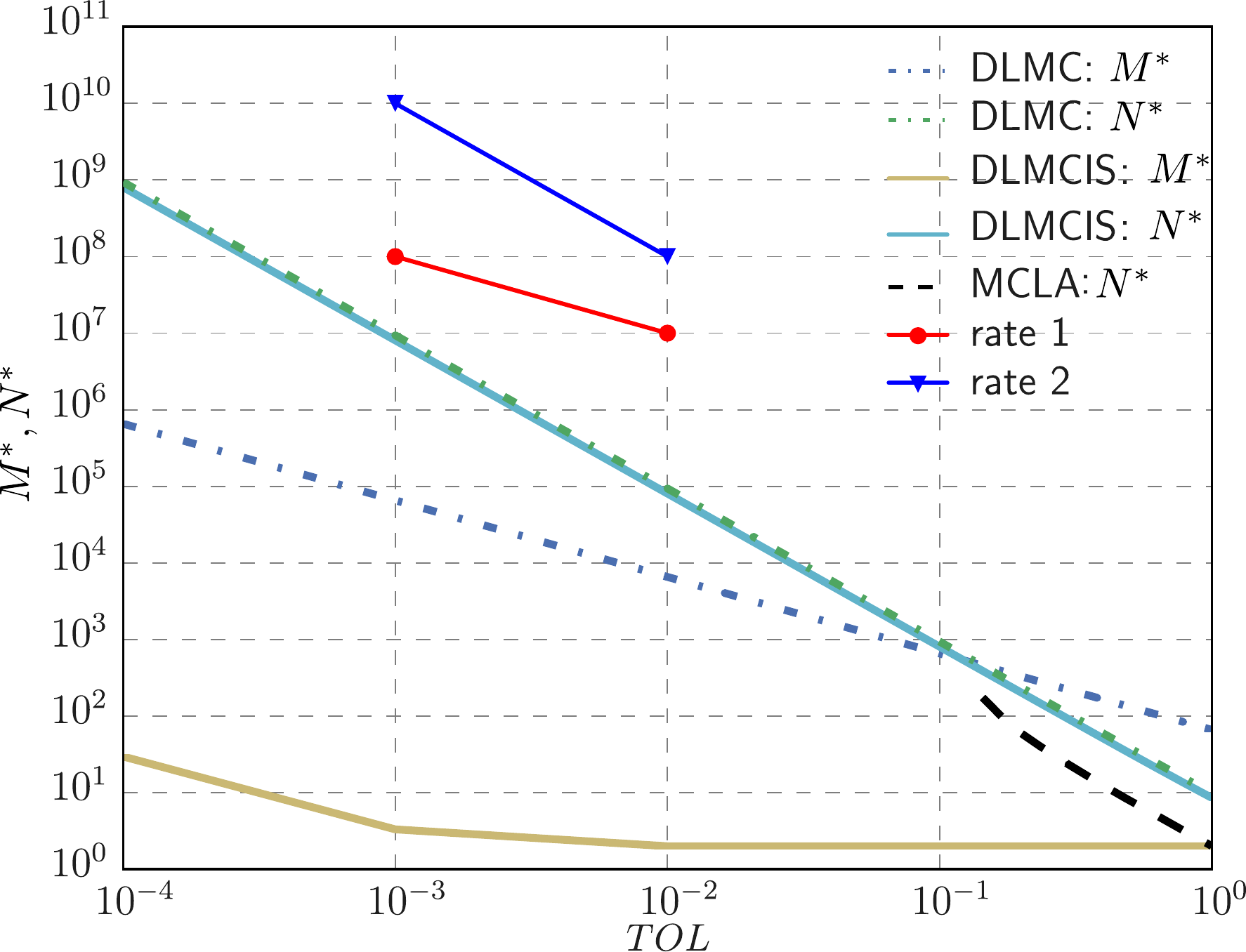}}
  \subfloat[$N_e=10$]{\label{optimality_MN_U.b}\includegraphics[width=0.5\linewidth]{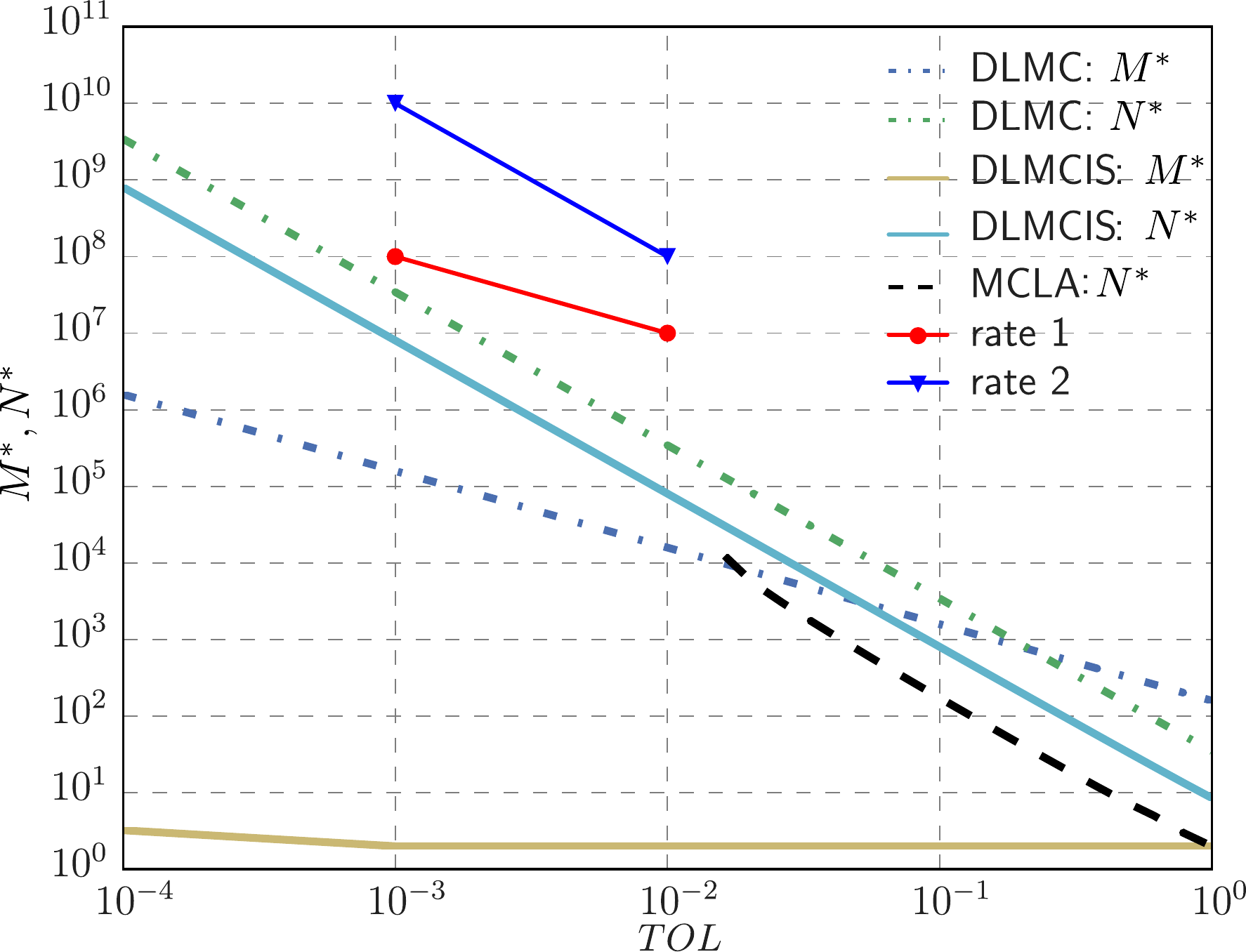}}
  \caption{Optimal setting (outer $N^*$ and inner $M^*$ number of samples) vs. tolerance for nonlinear scalar model with uniform prior (Example 2).}
  \label{optimality_MN_U}
\end{figure}

We start with a tolerance range $[10^{-4},1]$ to analyze the optimal method parameters, i.e., $N^*$, $M^*$, and $\kappa^*$. Figure \ref{optimality_MN_U} shows the optimal setting for $N_e=1$ (top) and $N_e=10$ (bottom). MCLA can reach lower tolerances with $N_e=10$ than with $N_e=1$, which is a constraint given by the Laplace bias $\mathcal{O}\left(1/N_e\right)$. DLMCIS can achieve an accuracy of $10^{-3}$ with $M^{*}=5$, compared to $M^*=10^5$ with DLMC. However, there is an additional cost associated with the change of measure in DLMCIS; namely, the method requires about $30$ extra forward model solves per outer sample $\theta_n$ in order to find $\hat{\theta}_n$ by solving an optimization problem \eqref{fittheta}.

The optimal balance parameter $\kappa^*$ is shown against $\hbox{TOL}$ in Figure \ref{optimality_nu_U}; in contrast to the previous example, here $\kappa^*$ is not constant. For MCLA, the optimal value $\kappa^*$ is quadratically decreasing to zero as $\hbox{TOL}$ approaches the size of the Laplace bias. MCLA is unable to achieve an accuracy better than the Laplace bias; accordingly, the Laplace bias decreases linearly with $N_e$.

\begin{figure}[!htp]
  \centering
  \subfloat[$N_e=1$]{\label{optimality_nu_U.a}\includegraphics[width=0.5\linewidth]{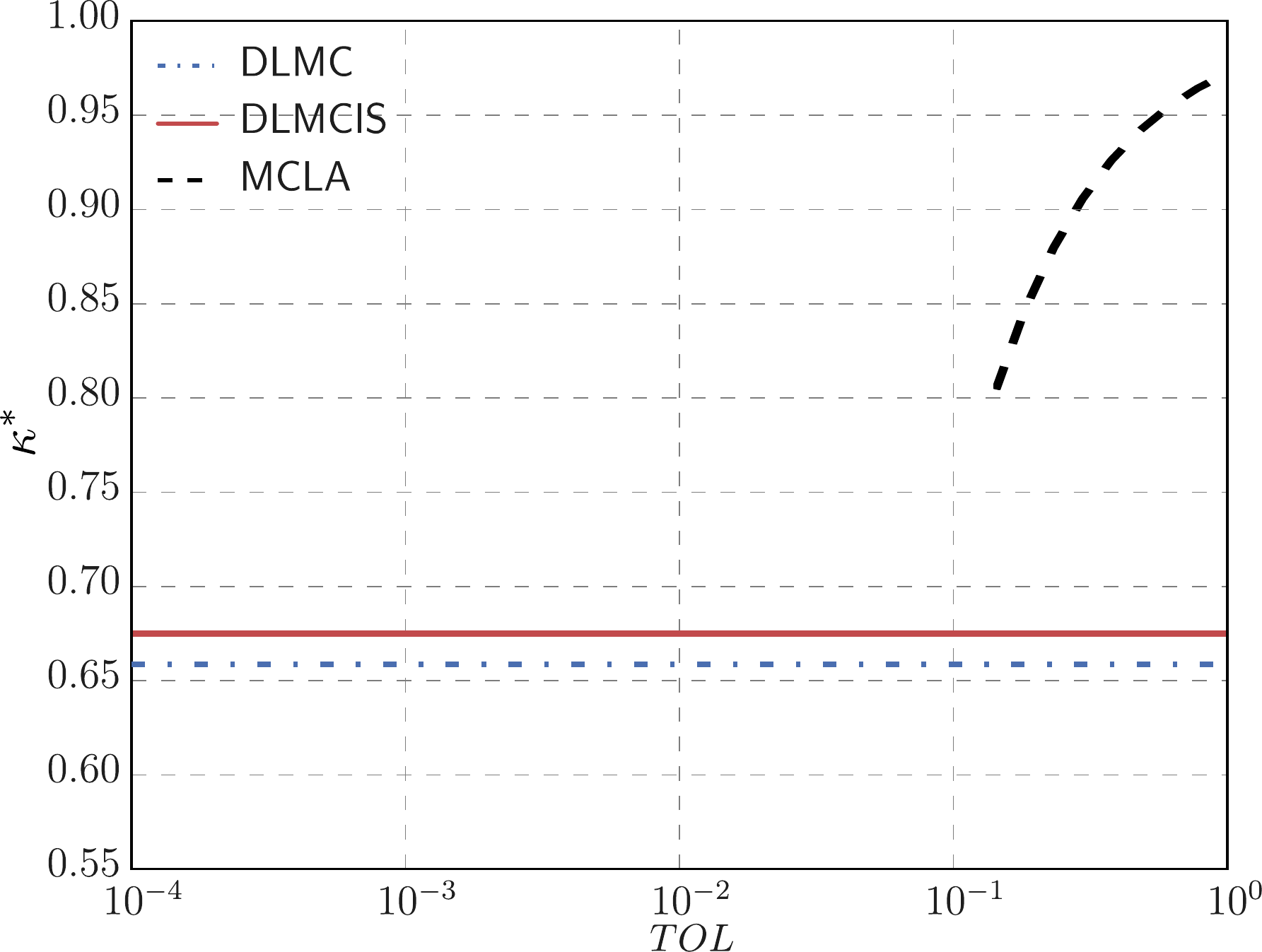}}
  \subfloat[$N_e=10$]{\label{optimality_nu_U.b}\includegraphics[width=0.5\linewidth]{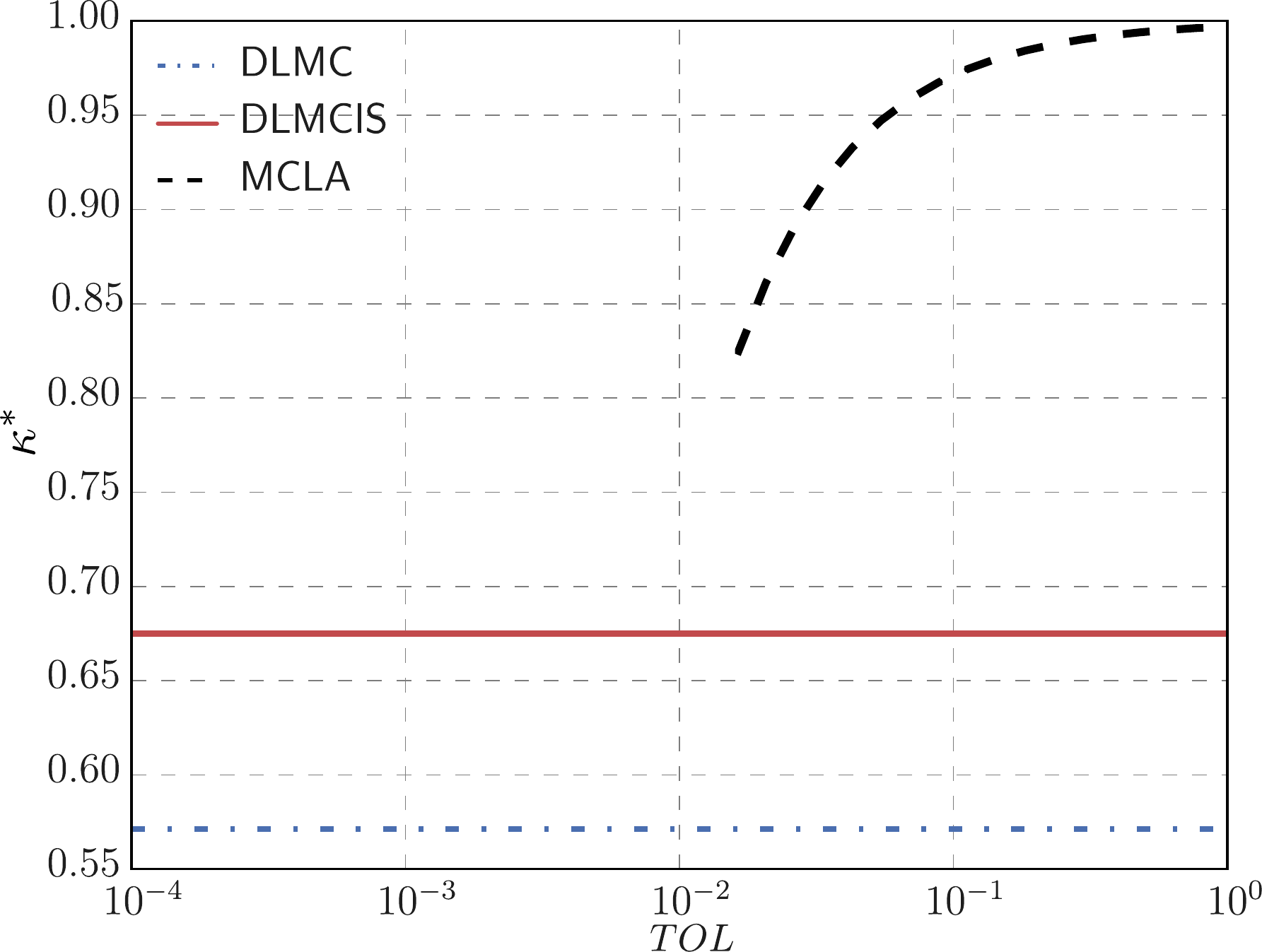}}
  \caption{Optimal split factor $\kappa^*$ vs. tolerance for nonlinear scalar model with uniform prior (Example 2).}
  \label{optimality_nu_U}
\end{figure}

\begin{figure}[!htp]
  \centering
  \subfloat[MCLA]{\label{errvstol_U.a}\includegraphics[width=0.5\linewidth]{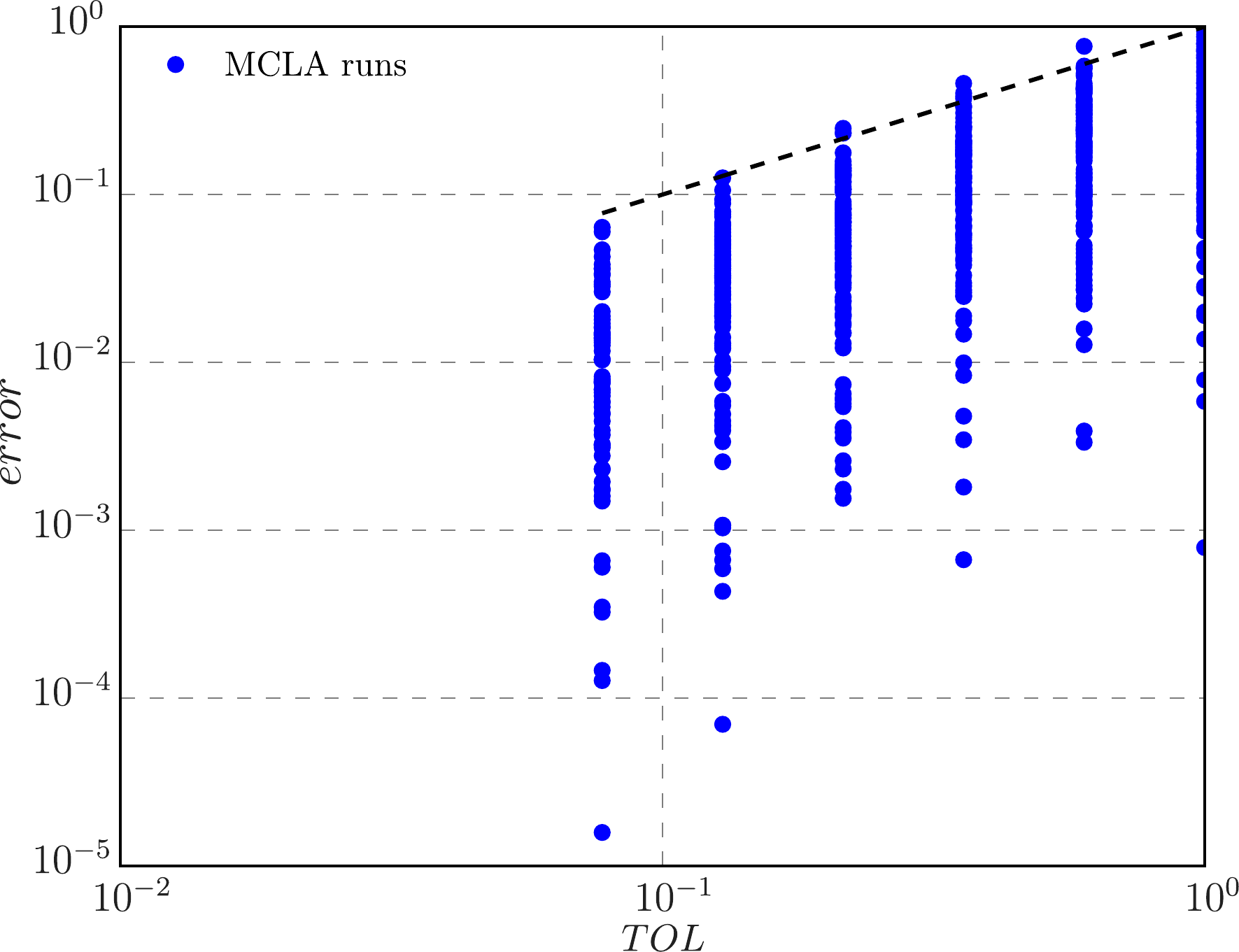}}
  \subfloat[DLMCIS]{\label{errvstol_U.b}\includegraphics[width=0.5\linewidth]{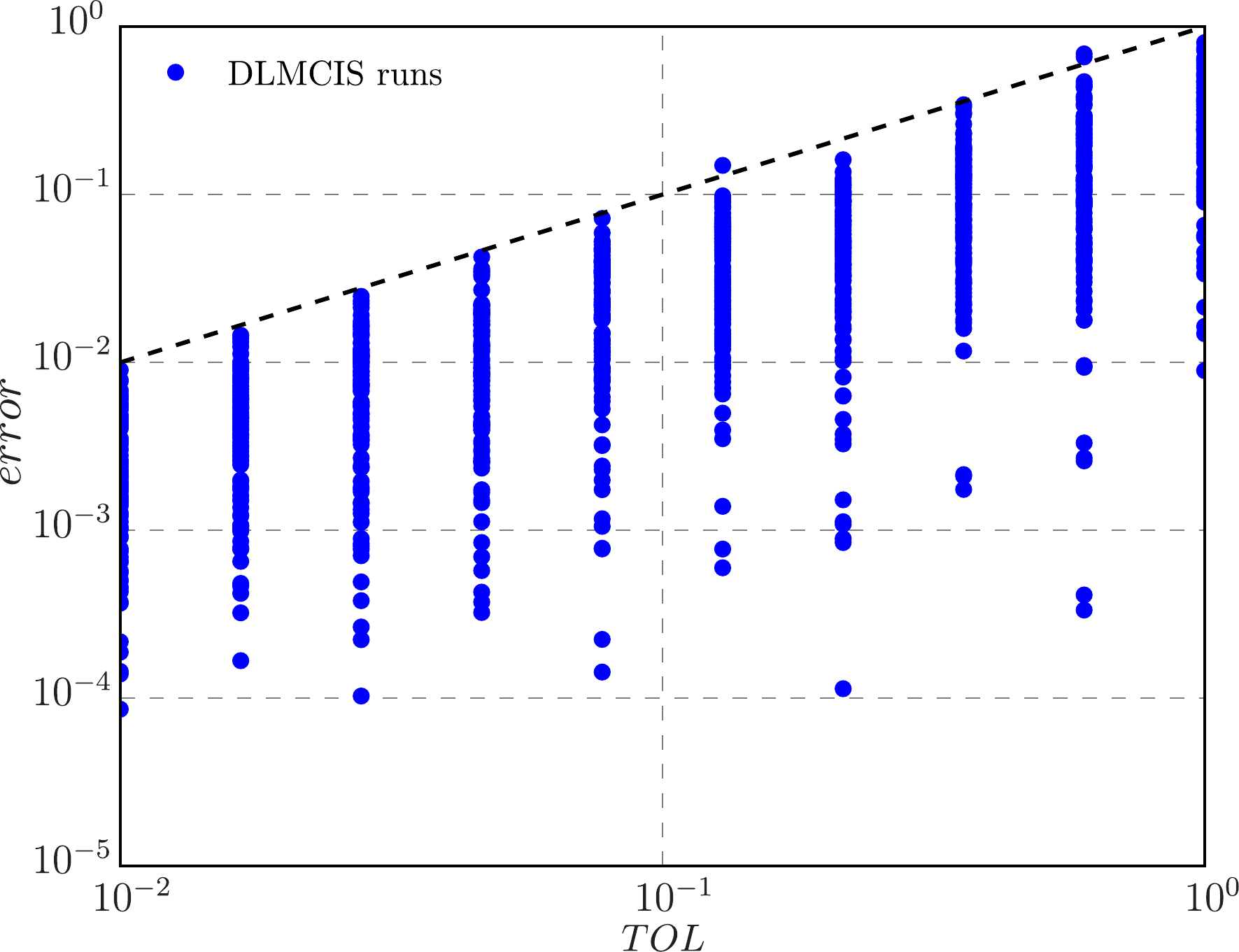}}
  \caption{Error vs. tolerance for nonlinear model with uniform prior (Example 2).}
  \label{errvstol_U}
\end{figure}

The absolute error for all three methods is below the specified tolerance $\hbox{TOL}$ in the range $[10^{-2},1]$ with at least $97.5\%$ probability (as specified); see Figure \ref{errvstol_U}. We follow the same procedure as in Example 1 for the estimation of the constants in the optimal parameter setting derivation.

\begin{figure}[!htp]
  \centering
  \subfloat[Work]{\label{time_work_U.b}\includegraphics[width=0.5\linewidth]{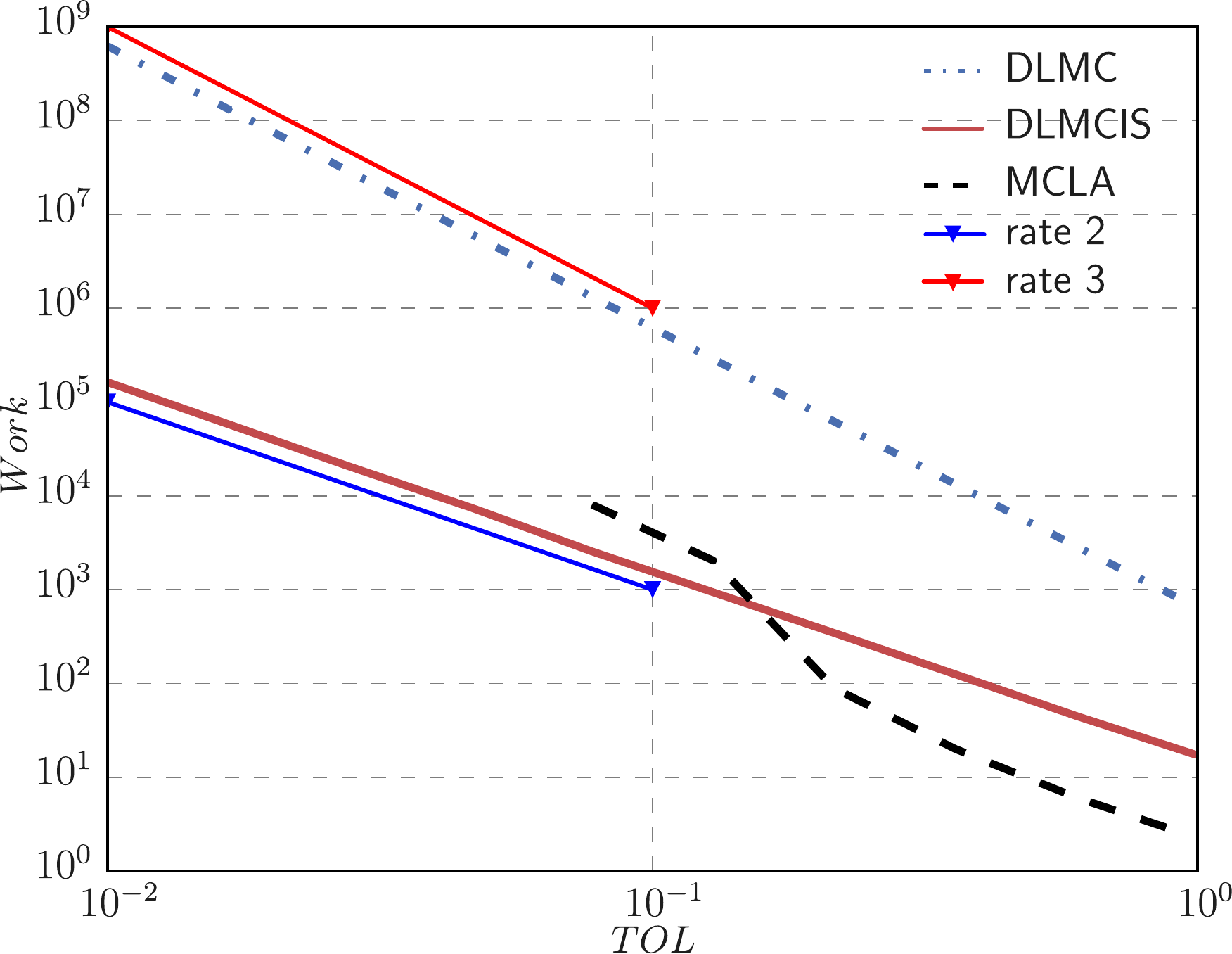}}
  \subfloat[Time]{\label{time_work_U.a}\includegraphics[width=0.5\linewidth]{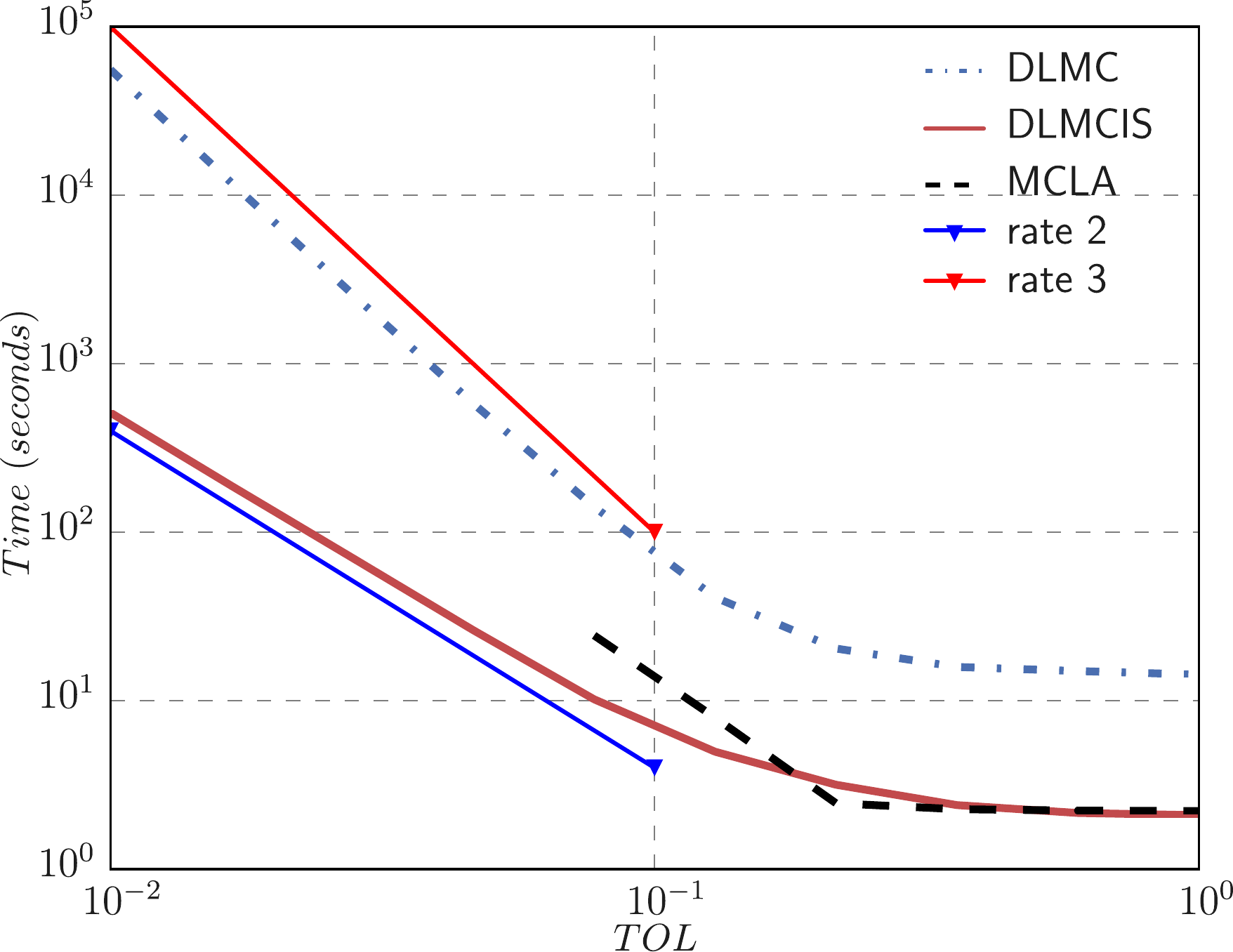}}
  \caption{Average computational work and running time vs. tolerance for nonlinear scalar model with uniform prior (Example 2).}
  \label{time_work_U}
\end{figure}

The running time against $\hbox{TOL}$ has rate 2 for MCLA and DLMCIS, and rate 3 for DLMC; see Figure \ref{time_work_U}. Again, as in Example 1, DLMCIS exhibits a pre-asymptotic rate that is one order less than the theoretical asymptotic rate, since $M^*$ is constant for tolerances within $[10^{-2},1]$. In this case, DLMCIS is about hundred times faster than DLMC for $\hbox{TOL}=0.01$. MCLA performs similarly to DLMCIS; however, MCLA could not be computed for tolerances lower than its bias.

\subsubsection{Expected information gain}
\begin{figure}[H]
  \centering
  \subfloat[MCLA for $\hbox{TOL}=3\times 10^{-2}$, and DLMCIS for $\hbox{TOL}=10^{-3}$.]{\label{I_U.a}\includegraphics[width=0.5\linewidth]{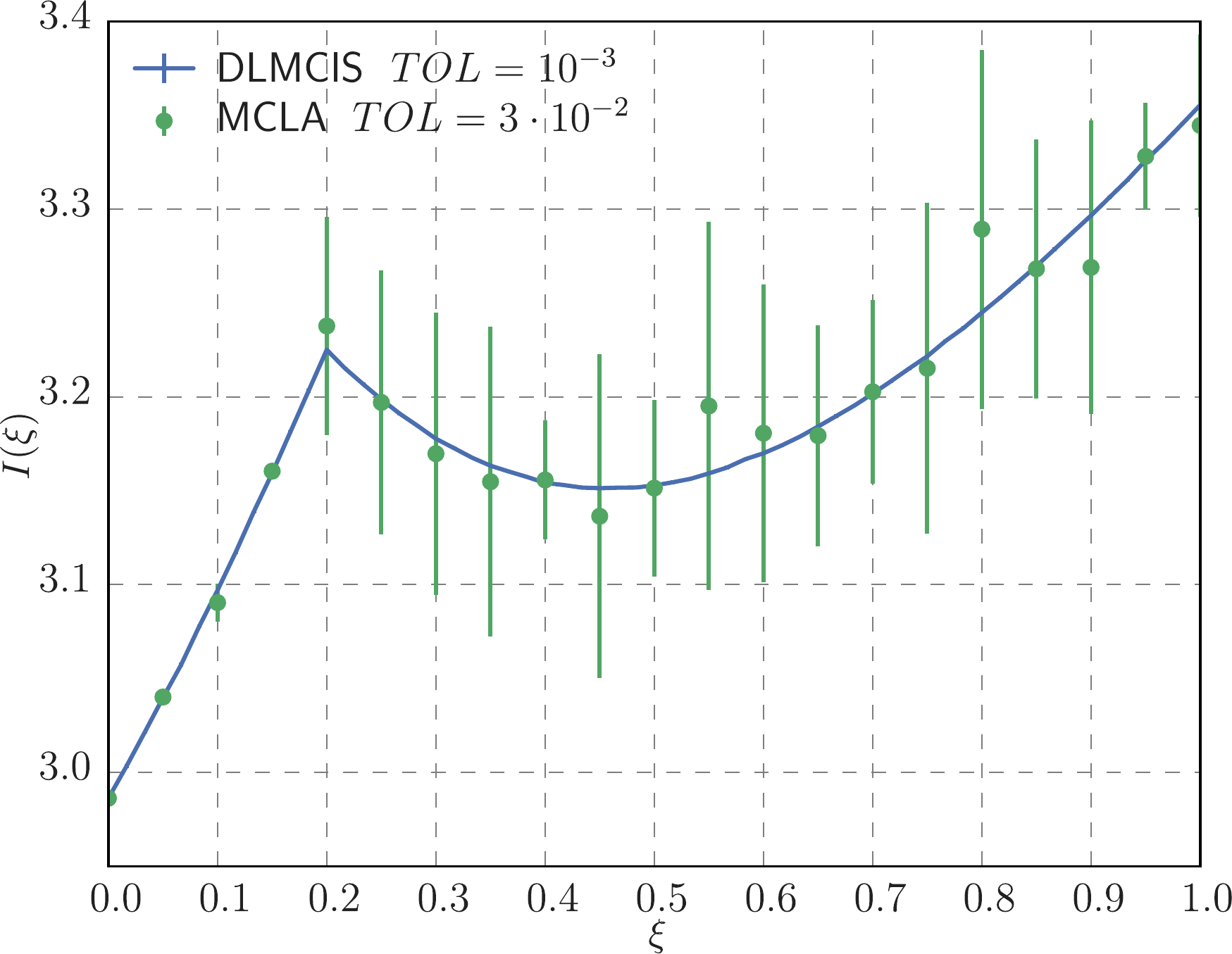}}
  \subfloat[MCLA when enforcing $\kappa^*=1$ for $\hbox{TOL}=10^{-3}$.]{\label{I_U.b}\includegraphics[width=0.5\linewidth]{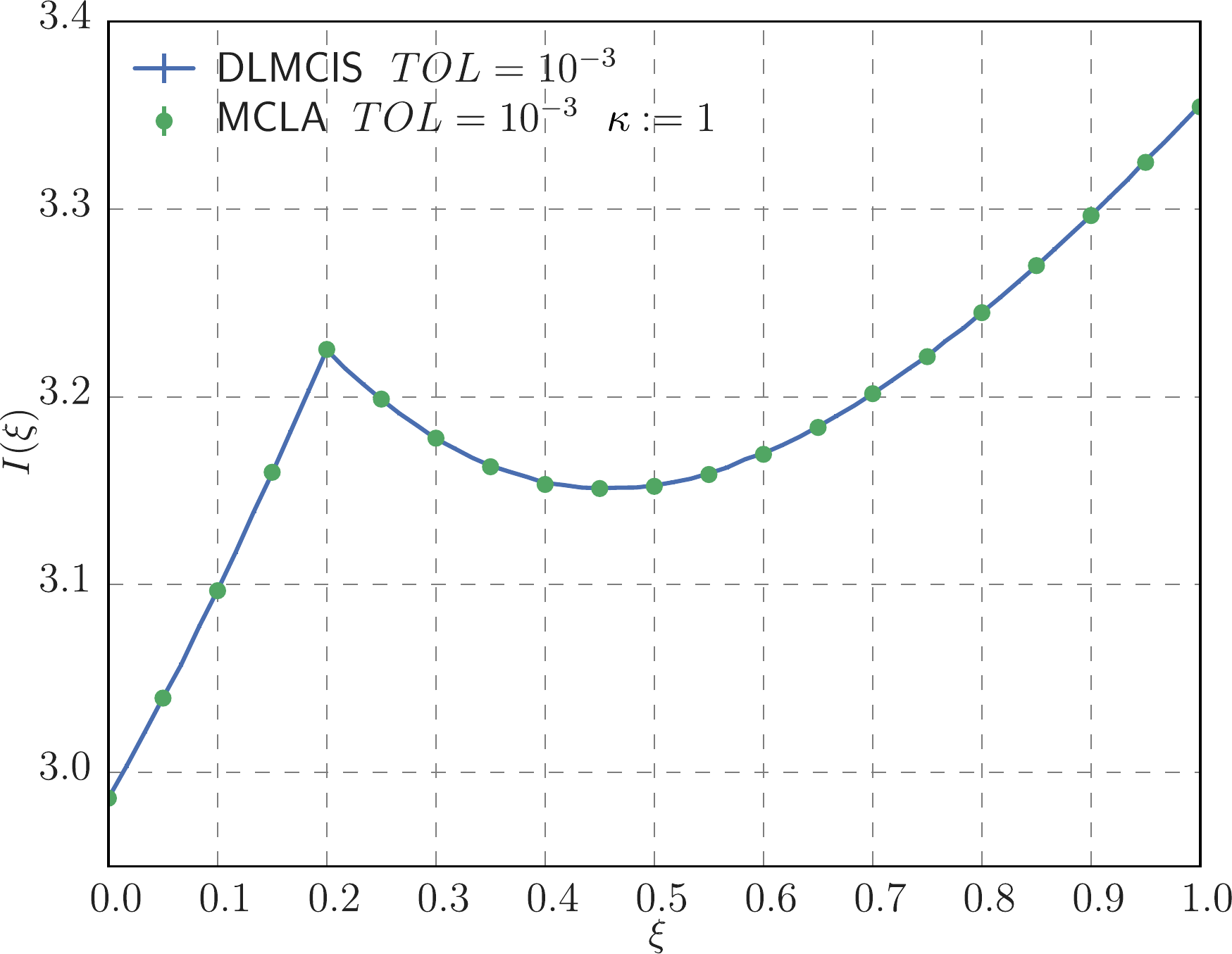}}
  \caption{Expected information gain for nonlinear scalar model with uniform prior (Example 2).}
  \label{I_U}
\end{figure}

In Figure \ref{I_U}, we present the estimation of the expected information gain using MCLA and DLMCIS for the experiment setups $\xi \in [0,1]$. In Figure \ref{I_U.a}, MCLA is applied for $\hbox{TOL}=3\times10^{-2}$, and DLMCIS for $\hbox{TOL}=10^{-3}$. The confidence bars of the MCLA curve show the $97.5\%$ confidence intervals. Different tolerances are specified for each of the two methods due to the Laplace bias constraint. However, in Figure \ref{I_U.b}, we omit the bias constraint by enforcing $\kappa^*=1$, and see that the MCLA curve matches well with the DLMCIS curve for $\hbox{TOL}=10^{-3}$. The resulting expected information gain curve is in agreement with the one reported previously, see \cite{huan}.

\subsection{Example 3: Electrical impedance tomography}
Here, we consider the optimal design of EIT experiments, an imaging technique in which the conductivity is inferred of a closed body from the measurements of electrodes placed on its boundary surface. In the forward problem, low-frequency electrical currents are injected through the electrodes attached to a composite laminate material, where each ply is orthotropic. The potential field in the body of the material is considered quasi-static for a given conductivity. The mathematical model is a second-order partial differential equation with an electrode boundary model.

\subsubsection{Mathematical model}
We use the complete electrode model (CEM) to formulate the EIT problem for a composite laminate material, see \citep{somersalo}. The composite body, denoted by $D$, is composed of $N_{_p}$ orthotropic laminated plies, which yields a macroscale, anisotropic material. The configuration is such that the plies overlap with their fibers at different orientation angles.

The equations governing the potential field $\bm{u}$ are
\begin{eqnarray}
\label{eq_eit1}
\nabla \cdot \bm{\jmath}(\omega, \bm{x}) \!\!\! & = & \!\!\! 0,\;\;\; \hbox{in}\;\;\; D, \text{and} \\
\bm{\jmath}(\omega, \bm{x}) \!\!\! & = & \!\!\! \bm{\bar{\sigma}}(\omega,\bm{x}) \cdot \nabla \bm{u}(\omega,\bm{x}),
\label{eq_eit2}
\end{eqnarray}
where $\bm{\jmath}$ is the flux of electric current, the conductivity $\bm{\bar{\sigma}}$ is given by
\begin{eqnarray*}
\bm{\bar{\sigma}}(\omega,\bm{x}) = \bm{Q}^T(\theta_k(\omega)) \cdot \bm{\sigma} \cdot \bm{Q}(\theta_k(\omega)),\;\;\;\hbox{for}\;\;\; \bm{x} \in D_k,\;\; k=1,\cdots,N_{_p},
\end{eqnarray*}
and the boundary conditions are specified in equations \eqref{eq_bcfree}-\eqref{eq_uzu}. The domain of a ply $k$ is denoted by $D_k$; thus, $D = \bigcup_{k=1}^{N_{_p}} D_k$. The orthogonal matrix $\bm{Q}(\theta_k)$ is the rotational matrix that defines the orientation of the fibres, in ply $k$ for a given angle $\theta_k$, while $\bm{\sigma}$ stands for the orthotropic conductivity:
\begin{eqnarray*}
\bm{Q}(\theta_k) =\begin{bmatrix}
\cos(\theta_k)  &     0      &      -\sin(\theta_k) \\ \\
0                       &     1     &        0 \\ \\
\sin(\theta_k)   &     0     &        \cos(\theta_k)
\end{bmatrix}
\;\;\;\hbox{and}\;\;\;
\bm{\sigma} =\begin{bmatrix}
\sigma_{1}      &       0               &       0 \\ \\
0                     &  \sigma_{2}     & 0 \\ \\
0                   &         0               & \sigma_{3}
\end{bmatrix}. 
\end{eqnarray*}  
The upper and lower surfaces of boundary $\partial D$ are equipped with $N_{_{\tiny{\hbox{el}}}}$ square-shaped electrodes $E_l$, $l=1,\cdots,N_{_{\tiny{\hbox{el}}}}$, with dimensions $e_{el}$. On the free surface of the boundary, $\partial D \backslash \left( \cup E_l\right)$, we assume a no-flux condition, i.e., no current flow in the out-of-surface direction:
\begin{eqnarray}
\bm{\jmath} \cdot \bm{n} = 0, \qquad \left(\bm{\bar{\sigma}}(\omega,\bm{x}) \cdot \nabla \bm{u}(\omega,\bm{x})\right) \cdot \bm{n} = 0,
\label{eq_bcfree}
\end{eqnarray}
where $\bm{n}$ represents the outward normal unit vector. CEM \cite{somersalo} is applied in the electrodes $E_l$. This means that we adopt \eqref{eq_bcfree} along with the assumption that the total injected current $I_l$ through each electrode is known and given by
\begin{eqnarray}
\label{eq_bcelec}
\int _{E_l} \bm{\jmath} \cdot \bm{n} \, d \bm{x} = I_l \;\;\;  \text{on} \;\;\;  E_l, \;\;\; l= 1,\cdots,N_{_{\tiny{\hbox{el}}}},
\end{eqnarray}
and that the shared interface of the electrode and the material has an infinitesimally thin layer with a surface impedance of $z_l$:
\begin{eqnarray}
\label{eq_uzu}
\frac{1}{E_l}\int_{E_l} \bm{u}d\bm{x} + z_l \int_{E_l}\bm{\jmath} \cdot \bm{n} d\bm{x}  = U_l \;\;\;\hbox{on} \;\;\;  E_l , \;\;\; l= 1,\cdots,N_{_{\tiny{\hbox{el}}}}.
\end{eqnarray}
For sake of well-posedness, the Kirchhoff law of charge conservation and the ground potential condition are set as constraints to complete the boundary model,
\begin{eqnarray}
\label{eq_constr}
\displaystyle{ \sum_{l=1}^{N_{_{\tiny{\hbox{el}}}}} I_l = 0} \;\;\; \hbox{and}\;\;\;
\displaystyle{ \sum_{l=1}^{N_{_{\tiny{\hbox{el}}}}} U_l = 0.}
\end{eqnarray}
In the rest of the paper, the EIT model refers to \eqref{eq_eit1}, \eqref{eq_eit2}, \eqref{eq_bcfree}, \eqref{eq_bcelec}, \eqref{eq_uzu},  and \eqref{eq_constr}. Due to the randomness of $\theta_k$, the conductivity field $\bm{\bar{\sigma}}$ is random and assumed to be a uniformly and strictly positive element of $L^\infty(\Omega\times D)$ in order to guarantee ellipticity. The vectors $\bm{I} = \left( I_1,I_2,\cdots,I_{N_{_{\tiny{\hbox{el}}}}} \right)^T$, and $\bm{U} = \left( U_1,U_2,\cdots,U_{N_{_{\tiny{\hbox{el}}}}} \right)^T$ respectively design the vector of the injected (deterministic) current and the vector measurement of the (random) potential at the electrodes. According to the constraints, \eqref{eq_constr}, $\bm{I}$ belongs to the mean-free subspace $\mathbb{R}^{N_{_{\tiny{\hbox{el}}}}}_{_{\tiny{\hbox{free}}}}$ of $\mathbb{R}^{N_{_{\tiny{\hbox{el}}}}}$ and $\bm{U}$ is a random element of $\mathbb{R}^{N_{_{\tiny{\hbox{el}}}}}_{_{\tiny{\hbox{free}}}}$.

\subsubsection{Finite element formulation}
We let $H^1(D)$ be the Hilbert space of $L^2$-integrable functions with $L^2$-integrable derivatives over the physical domain, $D$. Moreover, $\mathcal{H} \myeq H^1(D) \times \mathbb{R}^{N_{_{\tiny{\hbox{el}}}}}_{_{\tiny{\hbox{free}}}}$ denotes the space of the solution $(\bm{u}(\omega),\bm{U}(\omega))$ for a given random event $\omega \in \Omega$, and we introduce the Bochner space,
\begin{eqnarray*}
L^2_{\mathbb{P}} \left(\Omega;\mathcal{H} \right) \myeq \left\{ (\bm{u},\bm{U}): \Omega \rightarrow \mathcal{H} \;\;\;\hbox{s.t.}\;\; \int_{\Omega} \left\|(\bm{u}(\omega),\bm{U}(\omega))\right\|^2_{\mathcal{H}} d \mathbb{P}(\omega) < \infty \right\}.
\end{eqnarray*}
The variational form associated with the EIT problem finds $\left(\bm{u},\bm{U}\right) \in  L^2_{\mathbb{P}} \left(\Omega;\mathcal{H} \right)$ such that
\begin{eqnarray}
\mathbb{E} \left[B\left((\bm{u},\bm{U}),(\bm{v},\bm{V})\right) \right] = \bm{I} \cdot \mathbb{E} \left[ \bm{U} \right], \;\;\;\;\hbox{for all}\;\;\;(\bm{v},\bm{V}) \in L^2_{\mathbb{P}} \left(\Omega;\mathcal{H} \right),
\end{eqnarray}
where for any event $\omega \in \Omega$, the bilinear form $B:\mathcal{H}\times\mathcal{H} \rightarrow \mathbb{R}$ is 
\begin{eqnarray}
B\left((\bm{u},\bm{U}),(\bm{v},\bm{V})\right) = \int_D \bm{\jmath} \cdot \nabla \bm{v} d D + \sum_{l=1}^{N_{_{\tiny{\hbox{el}}}}} \frac{1}{z_l} \int_{E_l} \left(U_m -\bm{u} \right) \left(V_m - \bm{v} \right) d E_l.
\end{eqnarray}

We let $K$ be an element in the triangulation $\mathcal{T}_h$ and $\bar{D}\myeq\cup_{K\in\mathcal{T}_h}K$. Then we define the subspace $H^1_h(D_k) \myeq \{\bm{u}_h|_K \subset C^0(\bar{D}), \, \forall K \in \mathcal{T}_h\}$. Also, we let $\mathcal{H}_h \myeq H^1_h(D) \times \mathbb{R}^{N_{_{\tiny{\hbox{el}}}}}_{_{\tiny{\hbox{free}}}}\subset\mathcal{H}$. Next we consider the finite-dimensional problem in the Bubnov-Galerkin sense. The trial $(\bm{u}_h,\bm{U}_h)$ (test $(\bm{v}_h,\bm{V}_h)$) function pair is denoted such that $(\bm{u}_h,\bm{U}_h)\in\mathcal{H}_h$ ($(\bm{v}_h,\bm{V}_h)\in\mathcal{H}_h$). Finally, by rephrasing $L^2_{\mathbb{P}}$ on $(\Omega;\mathcal{H}_h)$, we conclude the finite element formulation.

\subsubsection{Bayesian experimental design formulation}
\begin{figure}[H]
\centering
\includegraphics[width=0.8\textwidth]{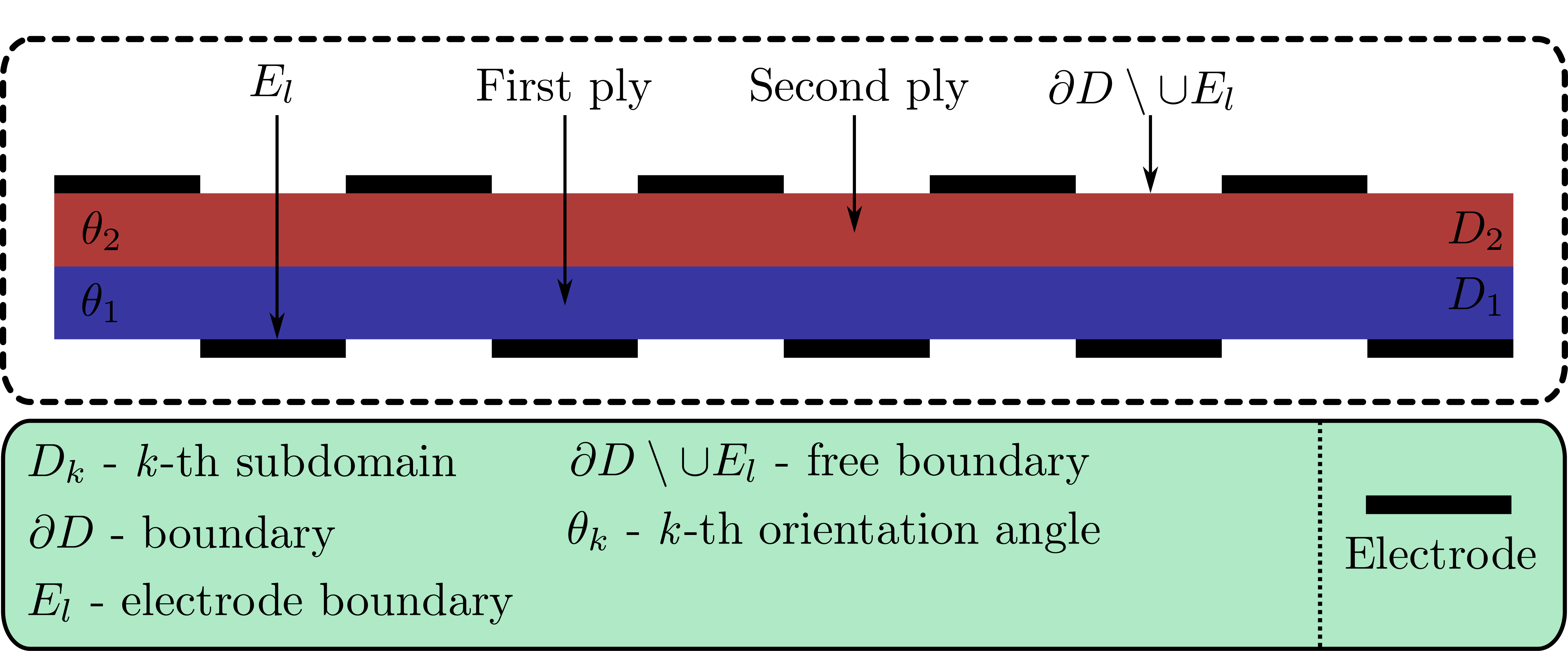}
\caption{Description of the domain for the EIT problem (Example 3).}
\label{eit_setup}
\end{figure}

We perform the experiments using the potentials measured at the electrodes. Hence, the Bayesian formulation of the EIT model is given by
\begin{eqnarray}
\bm{y}_i = \bm{g}_h(\bm{\theta}_t) + \bm{\epsilon}_i \myeq \bm{U}_h(\bm{\theta}_t) +\bm{\epsilon}_i, \quad for \quad i = 1, \cdots, N_e \,,
\end{eqnarray}
where $\bm{y}_i \in \mathbb{R}^{N_{_{\tiny{\hbox{el}}}}-1}$ (i.e., $q=N_{_{\tiny{\hbox{el}}}}-1$), $\bm{\theta}_t=(\theta_{t,1},\theta_{t,2})$, and the error distribution is Gaussian, i.e., $\bm{\epsilon} \sim \mathcal{N}(0,0.25)$.  We consider a uniform distribution to describe our prior knowledge of $\bm{\theta}_t$, i.e., 
\[
\pi(\theta_1)\sim\mathcal{U}\left(-\frac{\pi}{4}-0.05,-\frac{\pi}{4}+0.05\right),
\] 
and
\[
\pi(\theta_2)\sim\mathcal{U}\left(\frac{\pi}{4}-0.05,\frac{\pi}{4}+0.05\right).
\] 
We consider a body consisting of two plies whose parameters are $\sigma_{11}=0.05$, $\sigma_{22}=\sigma_{33}=10^{-3}$, and $z_l=0.1$. A total of 10 electrodes are placed on the surface of the plies to measure the potential at the electrodes. The orientations of the angles $\theta_1$ and $\theta_2$ of the fibers are the uncertain parameters. Figure \ref{eit_setup} presents the physical configuration of the experiment.

\begin{figure}[H]
  \centering
  \subfloat[Potential field]{\label{eit_1.a}\includegraphics[width=0.85\linewidth]{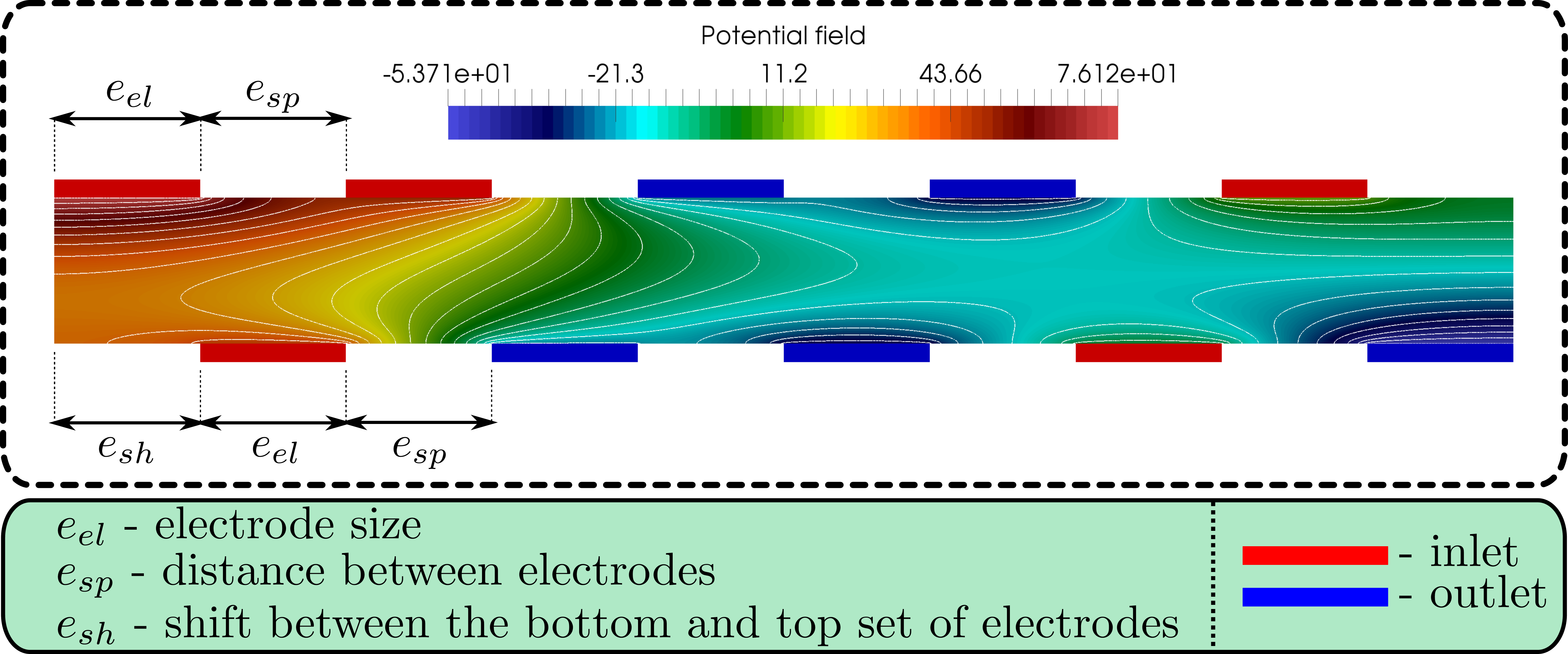}} \\
  \subfloat[Current flux]{\label{eit_1.b}\includegraphics[width=0.85\linewidth]{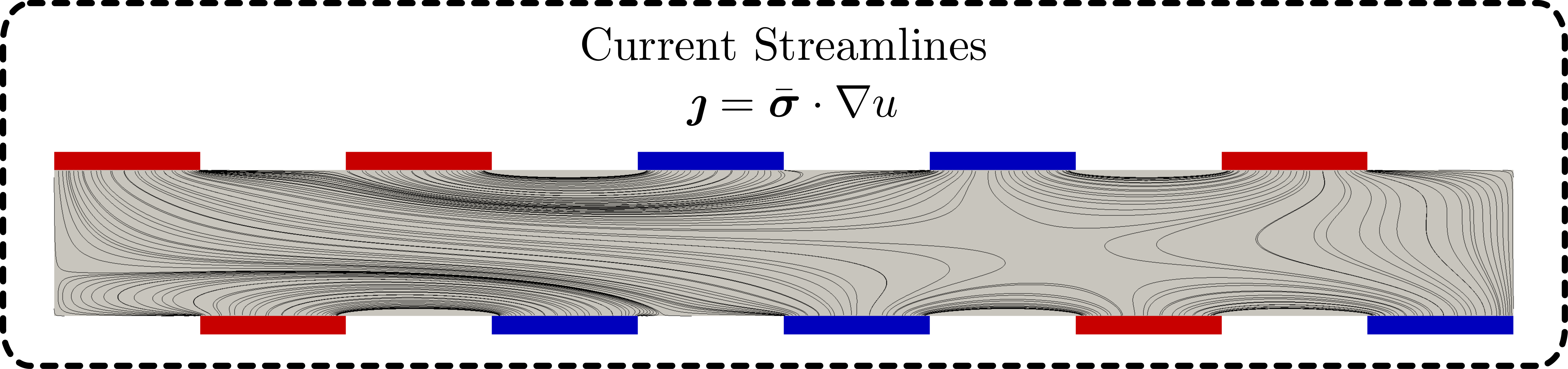}}
  \caption{Experimental setup for the EIT problem (Example 3).}
  \label{eit_1}
\end{figure}

Figure \ref{eit_1.a} presents the experiment for which we compute the respective optimal settings for MCLA, DLMC and DLMCIS. The domain is $D=[0,20]\times[0,2]$, and we set $e_{el}=e_{sh}=e_{sp}=2$. We depict the potential field with its equipotential lines. Figure \ref{eit_1.b} illustrates the current streamlines, $\bm{\jmath}$.

\begin{figure}[H]
  \centering
  \subfloat[Potential at $x_2=0$]{\label{pot_prof_2.a}\includegraphics[width=0.5\linewidth]{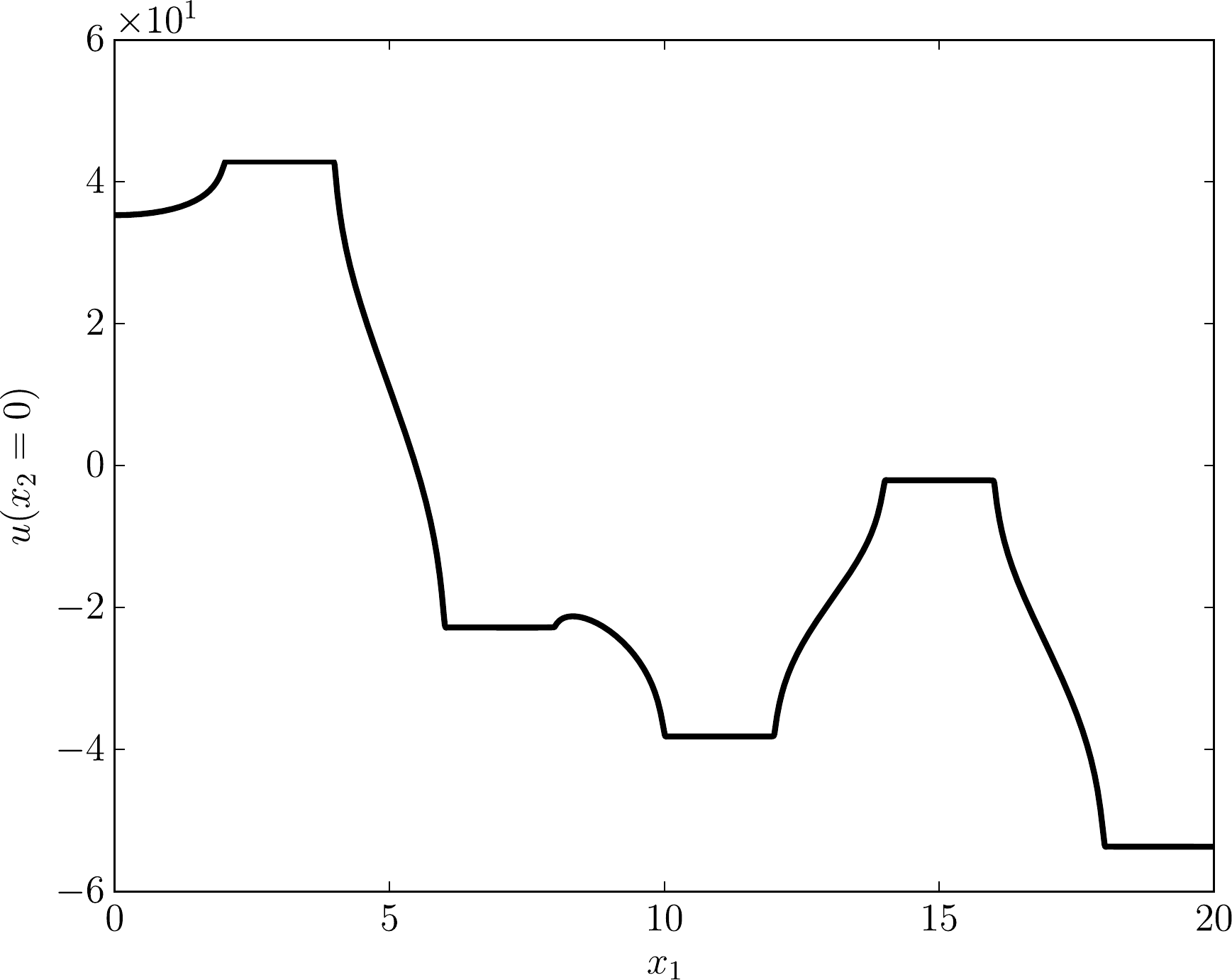}}
  \subfloat[Potential at $x_2=2$]{\label{pot_prof_2.b}\includegraphics[width=0.5\linewidth]{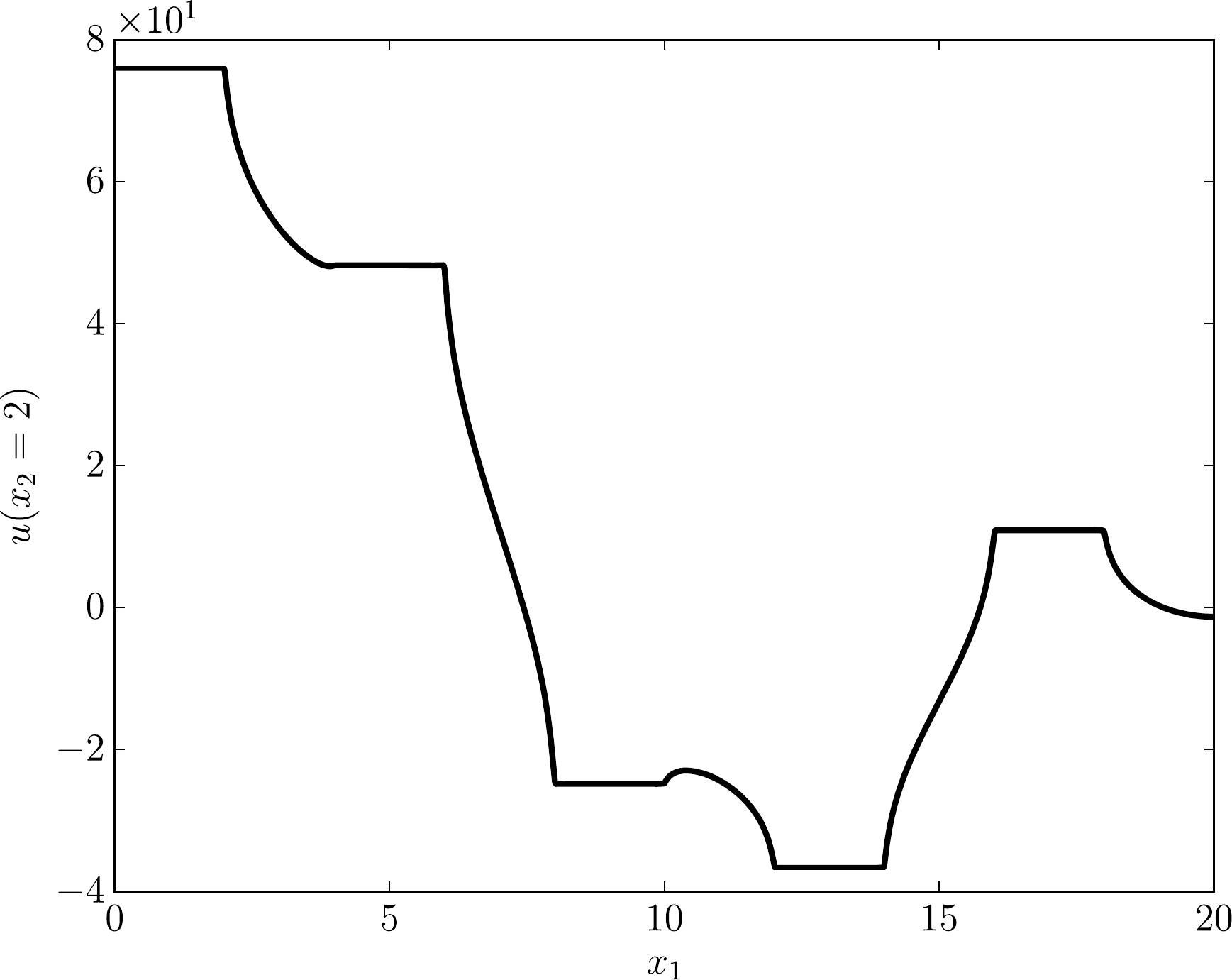}}
  \caption{Potential profiles for the EIT problem (Example 3).}
  \label{pot_prof_2}
\end{figure}

To understand how the potential field may behave at the surface boundaries, we show the potential profiles at $x_2=0$ (Figure \ref{pot_prof_2.a}) and $x_2=2$ (Figure \ref{pot_prof_2.b}). Sharp changes in the potential field are appreciated at the boundaries of the electrodes, meaning that the meshes utilized must be very fine. The solutions depicted in Figures \ref{eit_1} and \ref{eit_2} are computed using $1600 \times 160$ linear elements.

\subsubsection{Optimality and error convergence}
\begin{figure}[H]
  \centering
  \subfloat[$N^*$ and $M^*$]{\label{optimality_eit.MN}\includegraphics[width=0.75\linewidth]{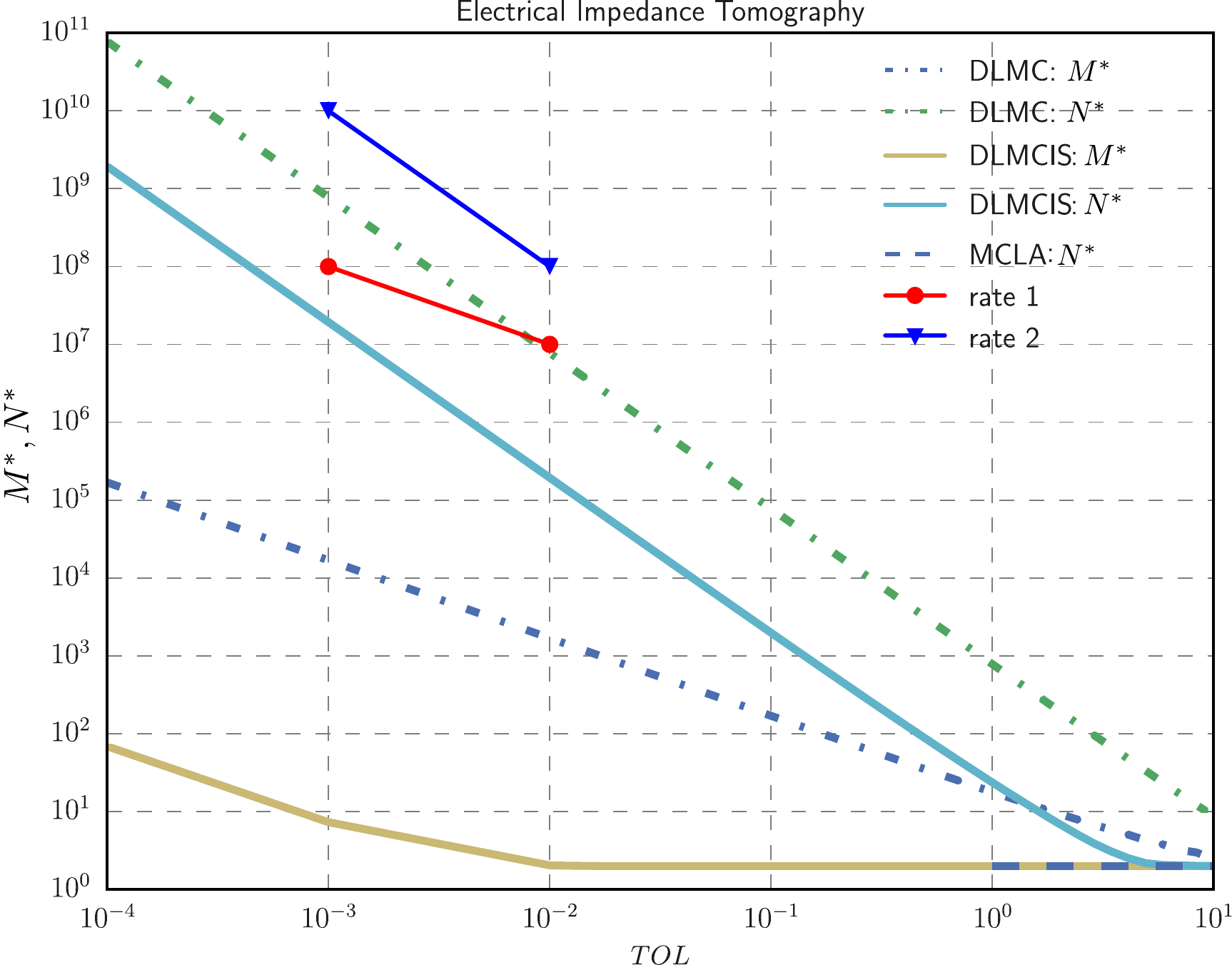}}\\
  \subfloat[$\kappa^*$]{\label{optimality_eit.nu}\includegraphics[width=0.5\linewidth]{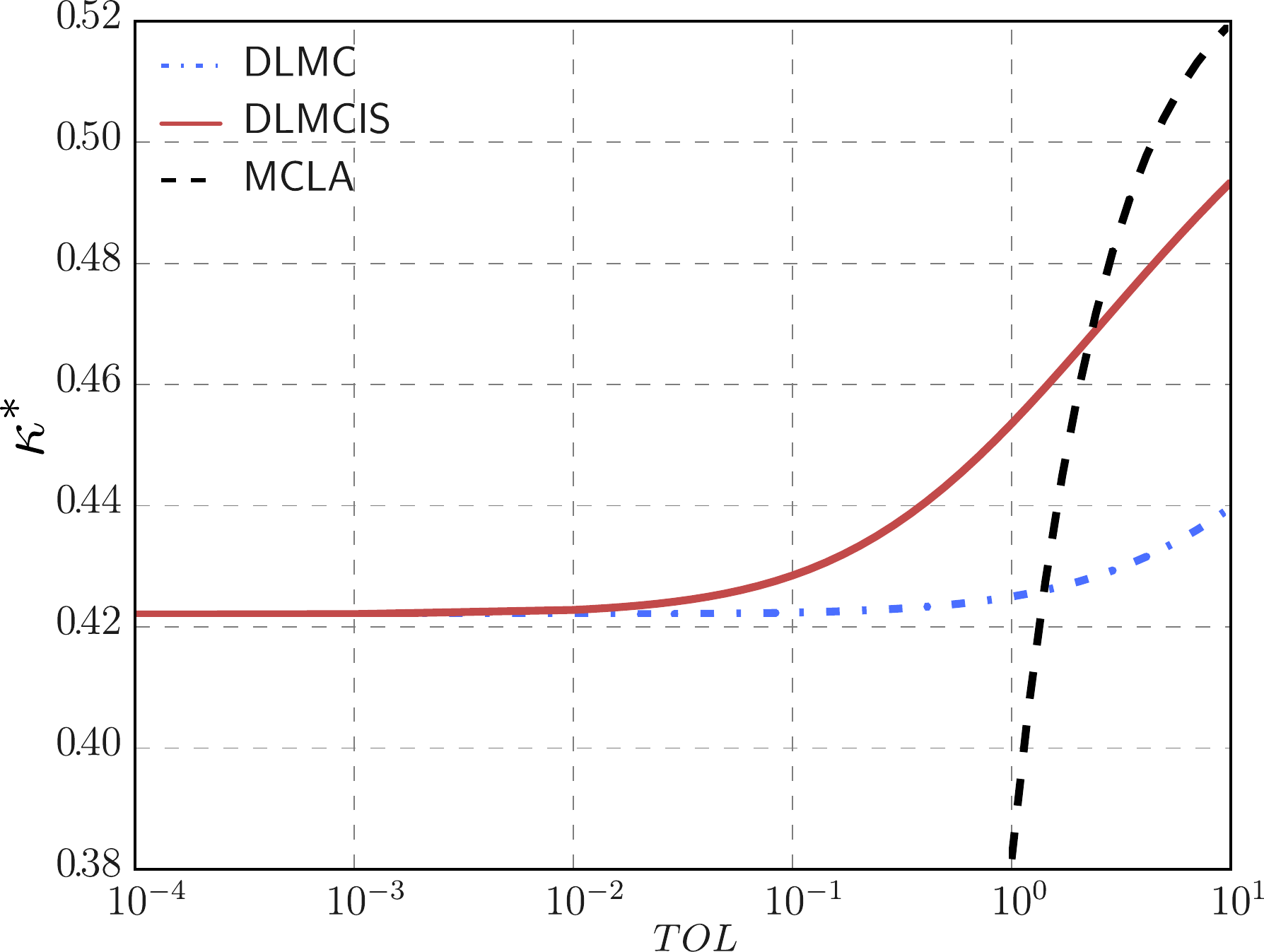}}
  \caption{Optimal setting (outer $N^*$ and inner $M^*$ number of samples, and balance $\kappa^*$) vs. tolerance for the EIT problem (Example 3).}
  \label{optimality_nu_U}
\end{figure}

\begin{figure}[H]
\centering
\subfloat[Error vs. tolerance (DLMCIS only)]{\label{vsTOL_eit.err.b}\includegraphics[width=0.5\linewidth]{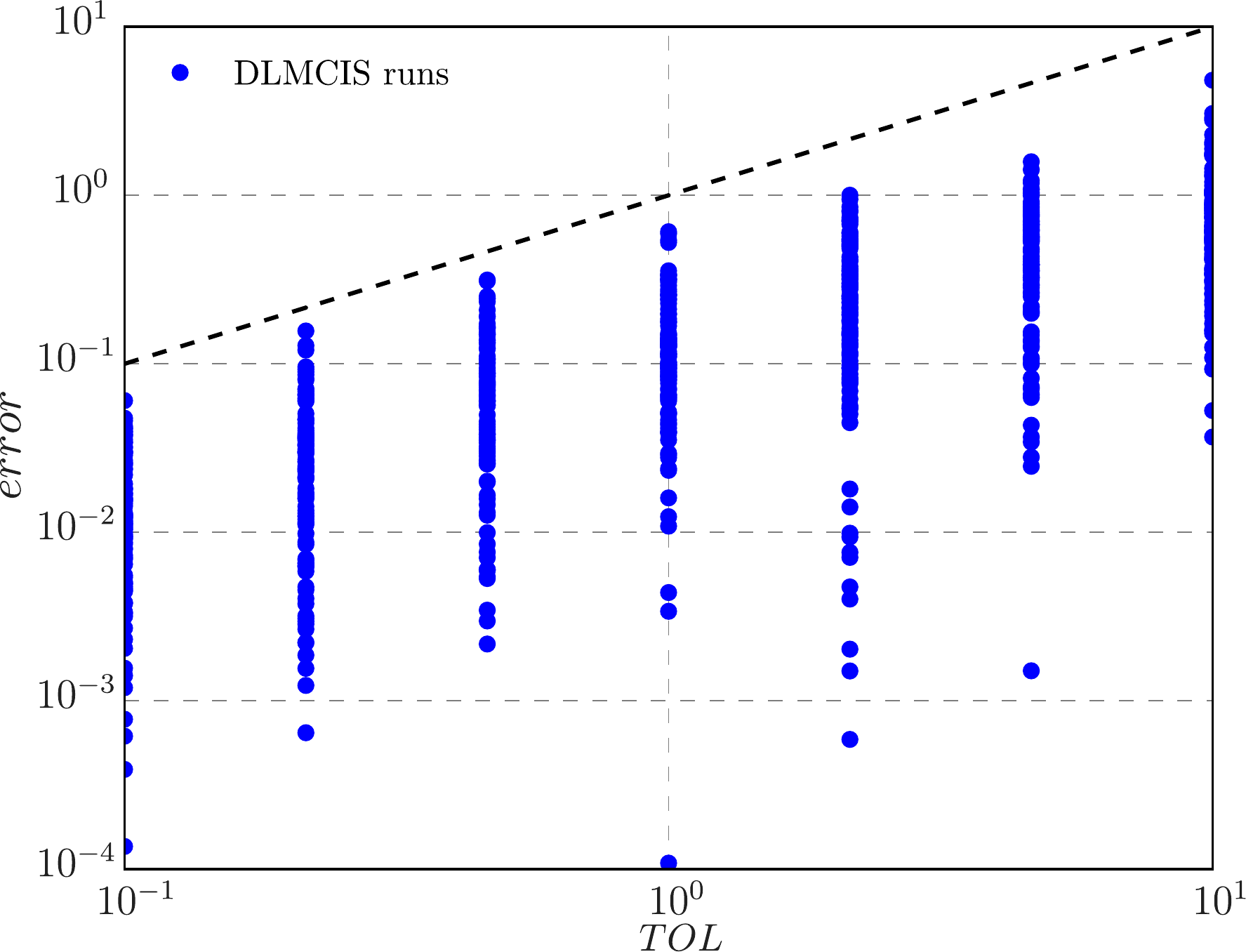}}
\subfloat[Running time vs. tolerance]{\label{vsTOL_eit}\includegraphics[width=0.5\linewidth]{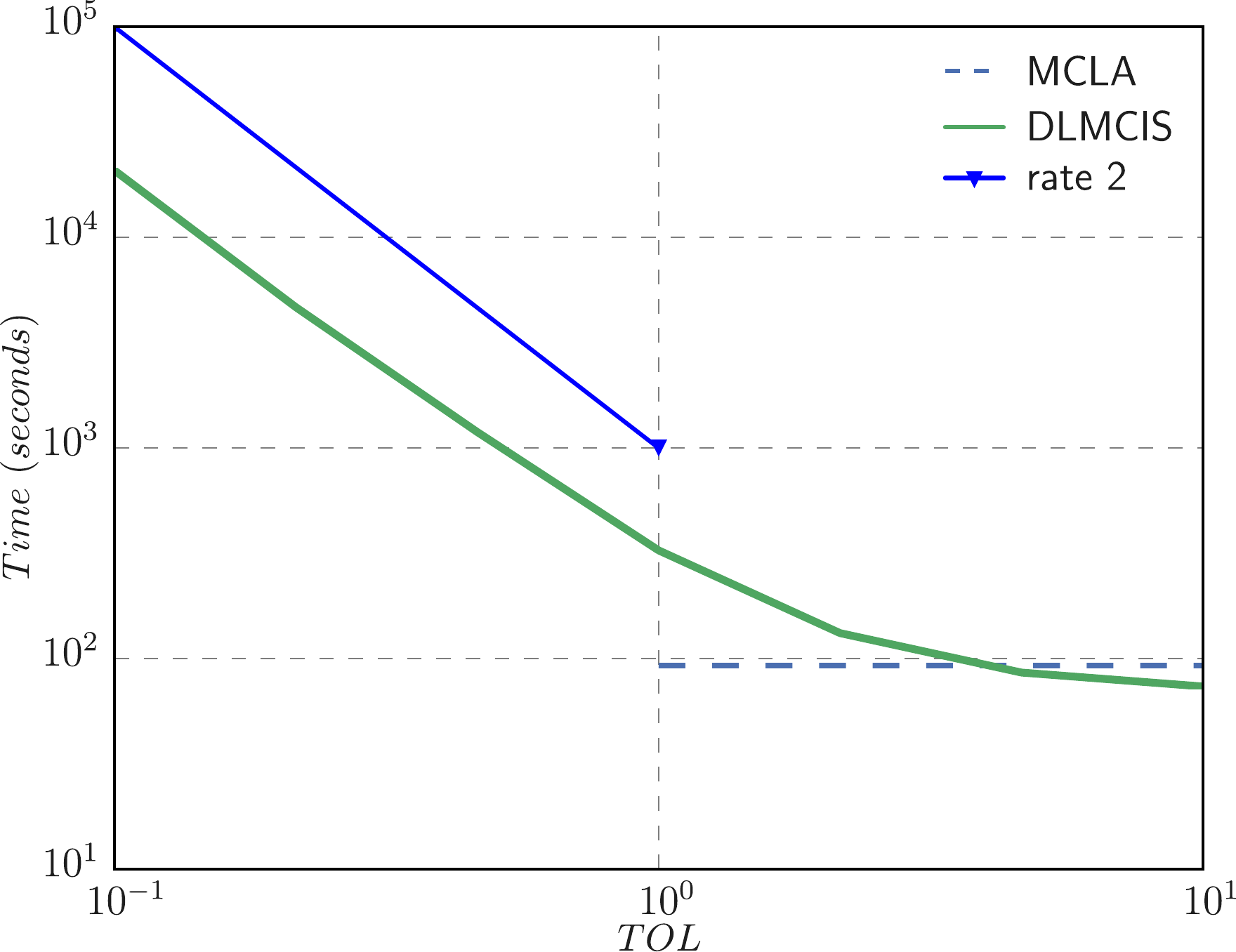}}
\caption{To the left is error vs. tolerance for DLMCIS only, and to the right is running time vs. tolerance for MCLA and DLMCIS (Example 3).}
\label{vsTOL_eit.err.xxxxxx}
\end{figure}

Figure \ref{optimality_eit.MN} shows the optimal number of outer and inner samples with respect to tolerance. Given the Laplace bias, MCLA can be performed to $\hbox{TOL}=1$, whereas DLMCIS, which is not limited by any uncontrolled bias, is carried out to $\hbox{TOL}=0.1$. The constants are estimated using $N=M=10$ for DLMCIS, instead of $N=M=100$ as in the previous examples, due to the expensive nature of this model. Even though the constants are roughly estimated, our results show a good agreement in tolerance versus absolute error for DLMCIS when using a confidence level of $97.5\%$, as seen in Figure \ref{vsTOL_eit.err.b}. 

The optimal discretization is roughly given by $h^* \approx 20\hbox{TOL}^{1/\eta}$. We constrain the maximum allowed $h$ to the mesh size that corresponds to 50 elements along the plies, and to 2 elements throughout the thickness of each ply; therefore, the same discretization size is used for $\hbox{TOL}=10$ and $\hbox{TOL}=1$. The discontinuity in the material and the discontinuity between the electrode and the no-flux boundary condition lead to discontinuity in the gradient of the potential field and, as a result, the $h$-convergence rate of $\bm{U}_h$ is $\eta \approx 1.15$ and the work rate is $\gamma \approx 3.5$.

To obtain an accurate reference solution for the computational error, we substitute the MC sampling for the outer integral by a Gauss-Legendre quadrature with $20\times20$ points, as the inner integral is low dimensional (in our case, $\text{dim}(\bm{\theta})=2$). Figure \ref{vsTOL_eit.err.b} shows the computational error in relation to tolerance. The relationship between computational work and tolerance is given in Figure \ref{vsTOL_eit}, where the theoretical asymptotic rates are numerically verified.

To reach $\hbox{TOL} \sim 0.01$, the number of samples needed for the inner loop in DLMCIS is $M^*=1$, and $M^* \sim 10^3$ for DLMC; for MCLA we are unable to reach accuracies better than $\hbox{TOL} \sim 1$. For all of the methods, the bias is larger than the statistical error; more specifically, $\kappa^*$ is less $0.5$ (see Figure \ref{optimality_eit.nu}). As discussed in Example 2, the Laplace approximation for the importance sampling requires about $30$ forward model evaluations per outer sample $\bm{\theta}_n$ to estimate $\bm{\hat{\theta}}_n$.

For MCLA, the rate of the computational work cannot be observed, as the work required is the same for $\hbox{TOL}=1$ and $10$. This is because the Laplace bias is the dominant error component and $N^*$ is equal to 1 for both tolerance levels. For DLMCIS, the observed computational work is $\mathcal{O}\left(\hbox{TOL}^{-2}\right)$, which is a pre-asymptotic rate that is two orders better than the asymptotic result. The reason for this is that the spatial discretization and the number of inner samples are constant for the range of tolerances considered.

\subsubsection{Expected information gain}
For this demonstration, we assign the design variables $\bm{\xi}=(\xi_1,\xi_2)$, where $\xi_1=e_{sh}$ and $\xi_2=e_{sp}$. We compute the expected information gain over the range $\bm{\xi} \in[0,2]\times[0.05,0.5]$. The reference solution is computed using the Gauss-Legendre quadrature for the outer integral to achieve a relative error $0.01$. The measurement error follows the distribution $\epsilon \sim \mathcal{N}(0,0.01)$.

\begin{figure}[H]
  \centering
  \subfloat[Potential field]{\label{eit_2.a}\includegraphics[width=0.85\linewidth]{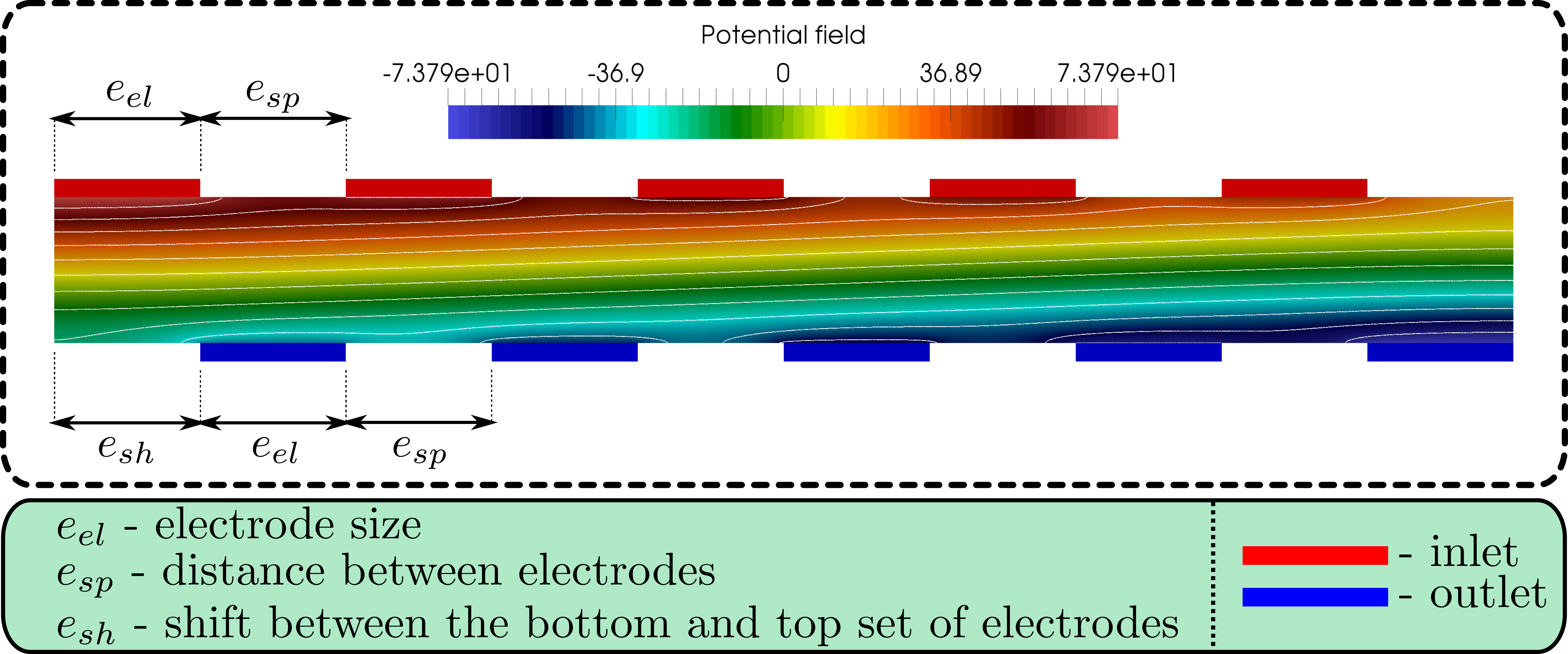}} \\
  \subfloat[Current flux]{\label{eit_2.b}\includegraphics[width=0.85\linewidth]{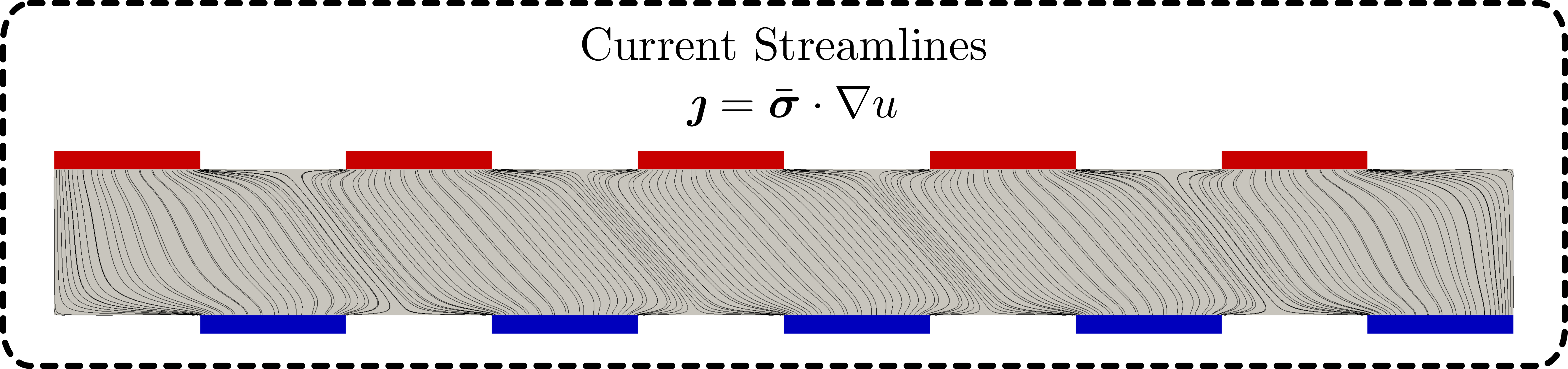}}
  \caption{The experimental setup for the expected information gain computation (Example 3).}
  \label{eit_2}
\end{figure}

The experiment is different from the previously considered experimental setup in Figure \ref{eit_1}, and is depicted in Figure \ref{eit_2.a}. All of the electrodes at the top boundary inject current, and the electrodes at the bottom act as outlets, and the domain is $D=[0,20]\times[0,2]$. The current streamlines, $\bm{\jmath}$, are shown in Figure \ref{eit_2.b}.

\begin{figure}[H]
  \centering
  \subfloat[Potential at $x_2=0$]{\label{pot_prof_1.a}\includegraphics[width=0.5\linewidth]{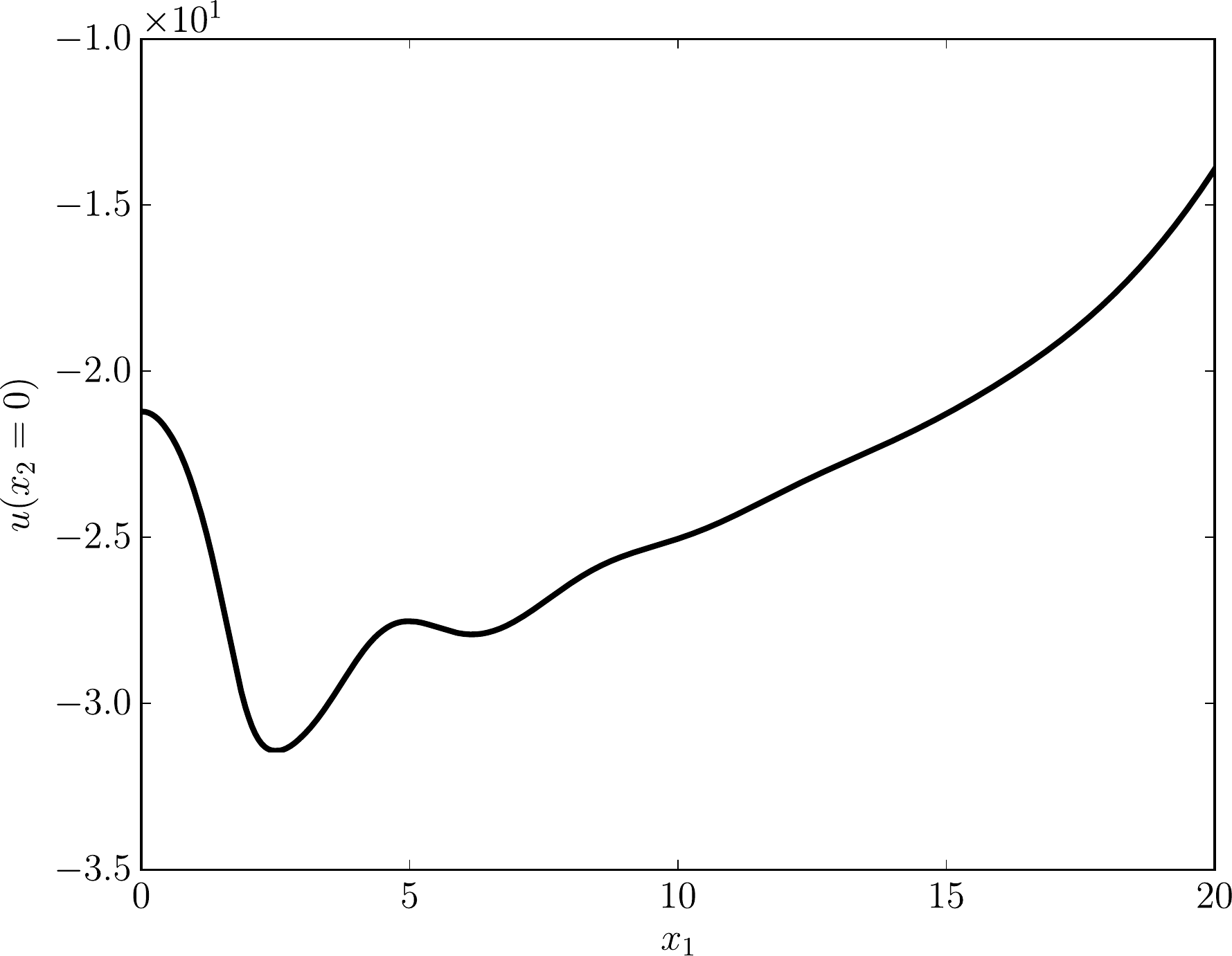}}
  \subfloat[Potential at $x_2=2$]{\label{pot_prof_1.b}\includegraphics[width=0.5\linewidth]{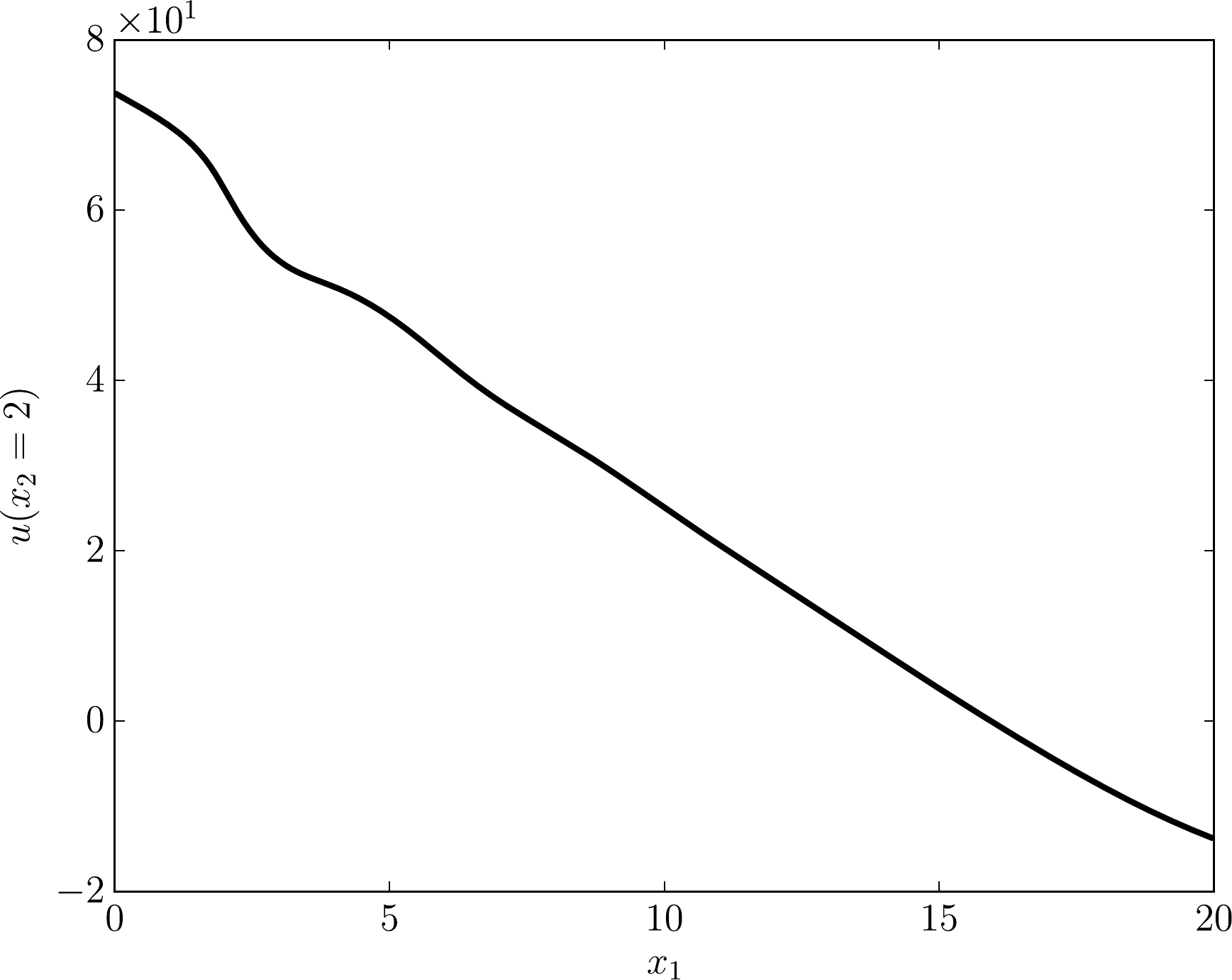}}
  \caption{The potential profile for the EIT problem (Example 3).}
  \label{pot_prof_1}
\end{figure}

Figure \ref{pot_prof_1.a} shows the potential profile at the bottom of the domain, $x_2=0$, whereas Figure \ref{pot_prof_1.b} shows the potential profile at the top of the domain $x_2=2$. The potential field is smoother than seen in the experimental setup shown in Figure \ref{eit_1}, due to all of the current being injected at the top and leaving at the bottom.

\begin{figure}[H]
\centering
\includegraphics[width=0.8\textwidth]{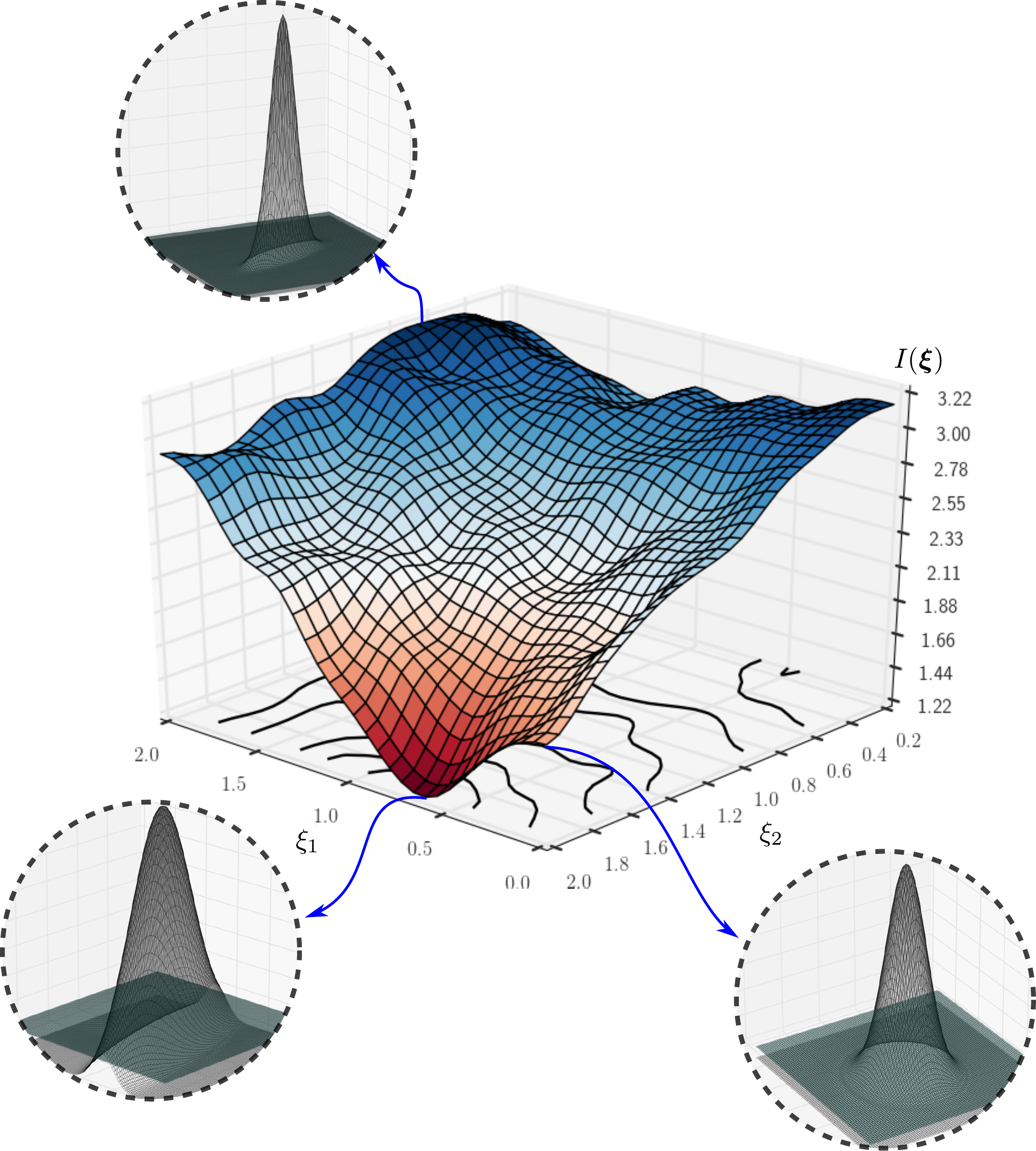}
\caption{Expected information gain for the EIT problem (Example 3).}
\label{I_eit}
\end{figure}

The expected information gain is presented in Figure \ref{I_eit} together with the posterior distributions for three experimental designs. In the figure, we see that there is oscillatory behavior at the surface; when the set of electrodes at the bottom is moved by $\xi_1$, the main current flow changes from one pair of electrodes to a different pair. Similar behavior occurs for the other design component, $\xi_2$. Moreover, the expected information gain is lower when the distance between the electrodes is larger.

\begin{figure}[H]
\centering
\includegraphics[width=0.65\textwidth]{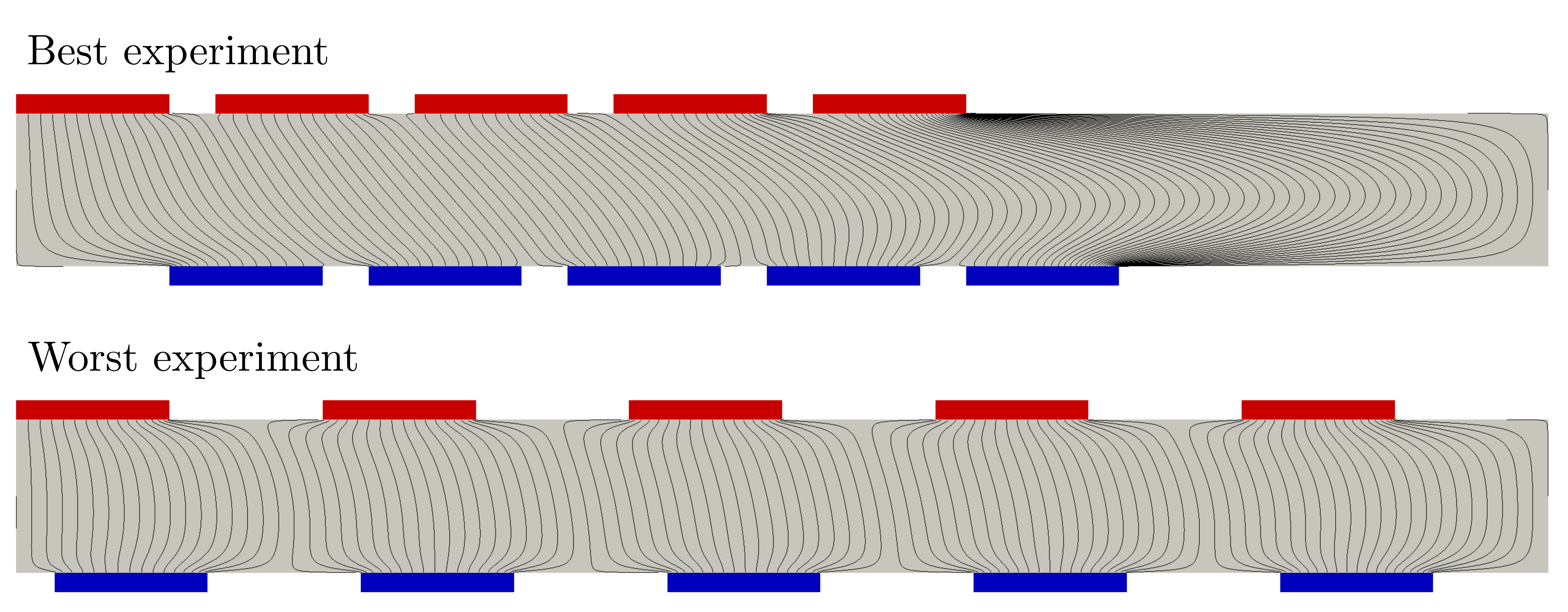}
\caption{The electrode placement and current flux for the worst (bottom) and best (top) experiments according to in accordance with the expected information gain presented in Figure \ref{I_eit} (Example 3).}
\label{I_eit_worst_best}
\end{figure}

The maximum expected information gain occurs at $\bm{\xi}=(2.0,0.6)$, and the corresponding posterior distribution is shown at the top of Figure \ref{I_eit}; the worst experiment occurs at $\bm{\xi}=(0.5,2.0)$, and its posterior distribution is depicted in the bottom-left of the figure. We show the electrode placement and current flux for the best and worst experiment in Figure \ref{I_eit_worst_best}. The best experiment suggests that even better results can be achieved by further increasing the shift of the electrodes, $e_{sh}$. From a design perspective, it is interesting that the optimal distance between the electrodes is about a quarter of the electrode size.

\section*{Conclusion}
We presented a computationally efficient method, DLMCIS, for computing the expected information gain in optimal Bayesian experimental design. The method uses importance sampling in the inner loop of the classical DLMC method through a change of measure based on the Laplace method. We demonstrated that the use of importance sampling can substantially reduce the number of samples in the inner loop, leading to an improvement in the running time of DLMC by several orders of magnitude. Further benefits include preservation of the asymptotic unbiasedness and mitigates the risk of underflow in the computation of the inner loop. We derived optimal settings for DLMCIS, DLMC and MCLA. The methods were optimized by minimizing the average computational work for a given tolerance. We showed that DLMCIS achieves the same work rate as MCLA for higher tolerances since only a few inner samples were needed, and for lower tolerances MCLA is not applicable due to its inherit bias. The efficiency of DLMCIS was demonstrated for optimal sensor placement in an EIT problem, where the composite material model was solved using high-fidelity finite elements.

\section*{Acknowledgements}
The research reported in this publication was supported by funding from King Abdullah University of Science and Technology (KAUST); KAUST CRG3 Award Ref:2281 and the KAUST CRG4 Award Ref:2584.

\section*{References}

\bibliographystyle{elsarticle-num}

\begin{thebibliography}{10}
\expandafter\ifx\csname url\endcsname\relax
  \def\url#1{\texttt{#1}}\fi
\expandafter\ifx\csname urlprefix\endcsname\relax\def\urlprefix{URL }\fi
\expandafter\ifx\csname href\endcsname\relax
  \def\href#1#2{#2} \def\path#1{#1}\fi

\bibitem{kullback1951}
S.~Kullback, R.~A. Leibler, On information and sufficiency, Annals of
  Mathematical Statistics 22 (1951) 79--86.

\bibitem{kullback1959}
S.~Kullback, Information Theory and Statistics, Wiley, 1959.

\bibitem{gin}
J.~Ginebra, On the measure of the information in a statistical experiment,
  Bayesian Analysis 2 (2007) 167--211.

\bibitem{ryan}
K.~J. Ryan, Estimating expected information gains for experimental designs with
  application to the random fatigue-limit model, Journal of Computational and
  Graphical Statistics 12 (2003) 585--603.

\bibitem{huan}
X.~Huan, Y.~M. Marzouk, Simulation-based optimal {B}ayesian experimental design
  for nonlinear systems, Journal of Computational Physics 232~(1) (2013)
  288--317.

\bibitem{Stigler1986}
S.~M. Stigler, Laplace's 1774 memoir on inverse probability, Statistical
  Science 1 (1986) 359--363.

\bibitem{tierney1986}
L.~Tierney, J.~B. Kadane, Accurate approximations for posterior moments and
  marginal densities, Journal of American Statistical Association 81 (1986)
  82--86.

\bibitem{tierney1989}
L.~Tierney, R.E.Kass, J.B.Kadane, Fully exponential {L}aplace approximations to
  expectations and variances of nonpositive functions, Journal of American
  Statistical Association 710-716~(84) (1989) 710--716.

\bibitem{kass1990}
R.E.Kass, L.Tierney, J.B.Kadane, Essays in Honor of George Barnard (eds. S.
  Geisser, J. S. Hodges, S. J. Press, and A. Zellner), North--Holland, 1990,
  Ch. The Validity of Posterior Expansions Based on {L}aplace's Method, pp.
  473--488.

\bibitem{key15}
Q.~Long, M.~Scavino, R.~Tempone, S.~Wang, Fast estimation of expected
  information gains for {B}ayesian experimental designs based on {L}aplace
  approximations, Computer Methods in Applied Mechanics and Engineering 259
  (2013) 24--39.

\bibitem{key51}
Q.~Long, M.~Motamed, R.~Tempone, Fast {B}ayesian optimal experimental design
  for seismic source inversion, Computer Methods in Applied Mechanics and
  Engineering 155 (2015) 123--145.

\bibitem{papadimitriou}
C.~Papadimitriou, Optimal sensor placement methodology for parametric
  identification of structural systems, Journal of Sound and Vibration 278
  (2004) 923--947.

\bibitem{key52}
Q.~Long, M.~Scavino, R.~Tempone, S.~Wang, A {L}aplace method for
  under-determined {B}ayesian optimal experimental design, Computer Methods in
  Applied Mechanics and Engineering​ 285 (2015) 849--876.

\bibitem{long2016}
F.~Bisetti, D.~Kim, O.~Knio, Q.~Long, R.~Tempone, Optimal {B}ayesian
  experimental design for priors of compact support with appication to
  shock-tube experiments for combustion kinetics, International Journal for
  Numerical Methods in Engineering 108~(2) (2016) 136--155.

\bibitem{alex2016}
A.~Alexanderian, N.~Petra, G.~Stadler, O.~Ghattas, A fast and scalable method
  for a-optimal design of experiments for infinite-dimensional {B}ayesian
  nonlinear inverse problems, SIAM Journal on Scientific Computing 38~(1)
  (2016) A243--A272.

\bibitem{BG2016}
J.~Beck, S.~Guillas, Sequential design with mutual information for computer
  experiments (mice): Emulation of a tsunami model, SIAM/ASA Journal on
  Uncertainty Quantification 4~(1) (2016) 739--766.

\bibitem{F2015}
C.~Feng, Optimal bayesian experimental design in the presence of model error,
  Ph.D. thesis, Massachusetts Institute of Technology (2015).

\bibitem{key13}
D.~V. Lindley, On a measure of information provided by an experiment, The
  Annals of Mathematical Statistics 27 (1956) 986--1005.

\bibitem{key12}
C.~E. Shannon, A mathematical theory of communication, Bell System Technical
  Journal 27 (1948) 379--423.

\bibitem{key14}
M.~B. Giles, Multilevel {M}onte {C}arlo path simulation, Operations Research
  56~(3) (2008) 607--617.

\bibitem{key56}
N.~Collier, A.-L. Haji-Ali, F.~Nobile, E.~von Schwerin, R.~Tempone, A
  continuation multilevel monte carlo algorithm, BIT Numerical Mathematics
  (2014) 1--34.

\bibitem{key59}
K.~J. Ryan, Estimating expected information gains for experimental designs with
  application to the random fatigue-limit model, Journal of Computational and
  Graphical Statistics 12~(3) (2003) 585--603.

\bibitem{somersalo}
E.~Somersalo, M.~Cheney, D.~Isaacson., Existence and uniqueness for electrode
  models for electric current computed tomography, SIAM J. Appl. Math, 52
  (1992) 1023--1040.

\end{thebibliography}

\appendix

\section{Estimation of the bias and variance of the DLMC estimator} \label{ap:secA}
\begin{proof}
This is the proof of Proposition \ref{prop1}. First, we show that the bias of the DLMC estimator, $\mathcal{I}_{dl,h}$, for the expected information gain can be upper-bounded by
\begin{equation*}
 \lvert I - \mathbb{E}\left[\mathcal{I}_{dl,h}\right] \rvert \leq C_{dl,3}h^\eta + \frac{C_{dl,4}}{M}+o(h^{\eta})+\mathcal{O}\left(\frac{1}{M^2}\right),
\end{equation*}
where $I$ is the expected information gain, and the subscript $h$ of the estimator shows its dependence on the forward model, $\bm{g}_h$.

We decompose the bias as follows:
\begin{eqnarray}
\label{A1totalbias}
\lvert I - \mathbb{E}\left[\mathcal{I}_{dl,h}\right]  \rvert \le \lvert I - \mathbb{E}\left[\mathcal{I}_{dl} \right] \rvert + \lvert \mathbb{E} \left[ \mathcal{I}_{dl} - \mathcal{I}_{dl,h} \right] \rvert,
\end{eqnarray}
where $\mathcal{I}_{dl}=\lim_{h \to 0^{+}} \mathcal{I}_{dl,h}$. The numerical bias of the estimator, due to the mesh discretization of $\bm{g}_h$, follows as
\begin{eqnarray}
\lvert \mathbb{E} \left[ \mathcal{I}_{dl} - \mathcal{I}_{dl,h} \right] \rvert = C h^{\eta}+o(h^{\eta}),
\end{eqnarray}
for a constant $C>0$, where $\mathcal{O}(h^{\eta})$ is of the same order as for $\bm{g}_h$. This can be seen by tracking $\mathcal{O}(h^{\eta})$ through the computation of $\log{p(\bm{Y} \vert \bm{\theta})}$ in $\mathcal{I}_{dl,h}$.

The other bias component, $\lvert I - \mathbb{E}\left[\mathcal{I}_{dl}\right] \rvert$, resulting from the inner averaging $\hat{p}_M(\bm{Y})= \frac{1}{M}\sum_{m=1}^{M}{p(\bm{Y} \vert \bm{\tilde{\theta}}_{m})}$ within the logarithmic function of the estimator, is proportional to the number of samples, $M$. This bias term can be recast as:
\begin{eqnarray}
\lvert I - \mathbb{E}\left[\mathcal{I}_{dl} \right] \rvert
&=&\mathbb{E} \left[   \frac{1}{N} \sum_{n=1}^{N} \log\left( \frac{p(\bm{Y}_n \vert \bm{\theta}_n)}{p(\bm{Y}_n)} \right) \right] - \mathbb{E} \left[  \frac{1}{N} \sum_{n=1}^{N}{\log\left( \frac{p(\bm{Y}_n \vert \bm{\theta}_n)}{\hat{p}_M(\bm{Y_n})}\right)}\right] \nonumber \\
\label{eqerror_dl1}&=& \mathbb{E} \left[ \log \left(\hat{p}_M(\bm{Y}) \right)  \right] -  \mathbb{E} \left[ \log (p(\bm{Y})) \right].
\end{eqnarray}
We apply the second-order Taylor expansion of $\log(X)$ around $\mathbb{E}[X]$ for both terms in \eqref{eqerror_dl1}:
\begin{eqnarray*}
 \log(X) =\log \left( \mathbb{E}\left[X\right] \right) + \frac{1}{\mathbb{E}\left[X\right] } \left( X- \mathbb{E}\left[X\right] \right) - \frac{1}{2}\frac{1}{\mathbb{E}\left[X\right] ^2} \left( X - \mathbb{E}\left[X\right] \right)^2 + \mathcal{O}\left(\left[ X-\mathbb{E}\left[X\right]\right] ^3 \right).
\end{eqnarray*}
That is, the expectation of $\log(\hat{p}_M(\bm{Y}))$ in \eqref{eqerror_dl1} can be written as
\begin{align*}
  \mathbb{E} \left[ \log \left(\hat{p}_M(\bm{Y}) \right)  \right] &= \mathbb{E}\left[  \log\left( \mathbb{E}\left[\hat{p}_M(\bm{Y})\right] \right) \right] 
  - \frac{1}{2} \frac{1}{ \left( \mathbb{E} \left[\hat{p}_M(\bm{Y}) \right]\right)^2 }  \mathbb{E} \left[ \left(\hat{p}_M(\bm{Y})  - \mathbb{E} \left[\hat{p}_M(\bm{Y}) \right] \right)^2 \right]  + \mathcal{O}\left(\frac{1}{M^2}\right) \nonumber\\
 &  =  \log(p(\bm{Y}))  - \frac{1}{2}\frac{1}{p^2(\bm{Y})}\mathbb{E} \left[ \left( \hat{p}_M(\bm{Y})  - p(\bm{Y}) \right)^2 \right] + \mathcal{O}\left(\frac{1}{M^2}\right).
\end{align*}
Here, the order $\mathcal{O}\left(\frac{1}{M^{2}}\right)$ comes from the fourth-order term of the Taylor expansion. Using the above results, we derive the result for \eqref{eqerror_dl1} as
\begin{eqnarray*}
\lvert I - \mathbb{E}\left[\mathcal{I}_{dl} \right] \rvert
  &=& \frac{1}{2}\mathbb{E}\left[\frac{1}{p^2(\bm{Y})}  \mathbb{E} \left[ \left( \hat{p}_M(\bm{Y})  - p(\bm{Y}) \right)^2 \bigg\vert \bm{Y} \right] \right] + \mathcal{O}\left(\frac{1}{M^{2}}\right)\\
  &=& \frac{1}{2}\mathbb{E}\left[\frac{1}{p^2(\bm{Y})}  \mathbb{V} \left[ \hat{p}_M(\bm{Y}) \bigg\vert \bm{Y} \right] \right] + \mathcal{O}\left(\frac{1}{M^{2}}\right)\\
&=& \frac{1}{2}\mathbb{E}\left[\frac{1}{p^2(\bm{Y})} \frac{1}{M^2} \mathbb{V} \left[  \sum_{m=1}^{M}{p(\bm{Y} \vert \bm{\tilde{\theta}}_m) }\bigg\vert \bm{Y} \right] \right] + \mathcal{O}\left(\frac{1}{M^{2}}\right)\\
&=& \frac{1}{2}\mathbb{E}\left[\frac{1}{p^2(\bm{Y})} \frac{1}{M^2} \sum_{m=1}^{M} \mathbb{V} \left[  p(\bm{Y} \vert \bm{\tilde{\theta}})  \bigg\vert \bm{Y} \right] \right] + \mathcal{O}\left(\frac{1}{M^{2}}\right)\\
&=& \frac{1}{2M}\mathbb{E}\left[ \frac{\mathbb{V} \left[  p(\bm{Y} \vert \bm{\tilde{\theta}}) \bigg\vert \bm{Y} \right] } {p^2(\bm{Y})} \right] + \mathcal{O}\left(\frac{1}{M^{2}}\right).
\end{eqnarray*}

The above result completes the derivation of the upper bound of the total bias \eqref{A1totalbias} with respect to $M$ and $h$. Throughout the proof, we drop the subscripts when they are not required since they are identically distributed variables.

Next, we derive the expression for the variance of the estimator,
\begin{eqnarray*}
 \mathbb{V}\left[ \mathcal{I}_{dl} \right] = \frac{C_{dl,1}}{N} + \frac{C_{dl,2}}{NM}+\mathcal{O}\left(\frac{1}{NM^2}\right).
\end{eqnarray*}

By the law of total variance, we attain the following result:
\begin{eqnarray}
  \mathbb{V}\left[\mathcal{I}_{dl} \right] 
  &=& \mathbb{V}\left[  \frac{1}{N} \sum_{n=1}^{N}{\log\left( \frac{p(\bm{Y}_n \vert \bm{\theta}_n)}{ \hat{p}_M(\bm{Y})}\right)} \right] \nonumber \\
  &=& \frac{1}{N} \mathbb{V}\left[ \mathbb{E} \left[ \log (p(\bm{Y} \vert \bm{\theta})) - \log \left(  \hat{p}_M(\bm{Y}) \right) \bigg\vert \bm{\theta},\bm{Y}\right]\right] \nonumber \\ 
  && +\frac{1}{N} \mathbb{E} \left[ \mathbb{V} \left[  \log (p(\bm{Y} \vert \bm{\theta})) - \log \left(  \hat{p}_M(\bm{Y}) \right) \bigg\vert \bm{\theta},\bm{Y}\right] \right] \nonumber \\
  &=& \frac{1}{N} \mathbb{V} \left[ \mathbb{E} \left[ \log (p(\bm{Y} \vert \bm{\theta})) - \log( \hat{p}_M(\bm{Y})) \bigg\vert \bm{\theta},\bm{Y} \right] \right]  + \frac{1}{N} \mathbb{E} \left[ \mathbb{V} \left[  \log \left(  \hat{p}_M(\bm{Y}) \right) \bigg\vert \bm{Y} \right] \right]. \nonumber \\ \label{A1var}
\end{eqnarray}
Here, we apply the second-order Taylor expansion around $\mathbb{E}\left[ \hat{p}_M(\bm{Y}) \right]$ for both terms. The first term of \eqref{A1var} can be rewritten as follows:
\begin{eqnarray}
   \quad & \quad &\frac{1}{N} \mathbb{V} \left[ \mathbb{E} \left[ \log (p(\bm{Y} \vert \bm{\theta})) - \log(\hat{p}_M(\bm{Y})) \bigg\vert \bm{\theta},\bm{Y} \right] \right] \nonumber \\
   \quad &=& \frac{1}{N} \mathbb{V}\left[ \mathbb{E} \left[\log\left(\frac{p(\bm{Y} \vert \bm{\theta})}{p(\bm{Y})}\right)-\frac{1}{2p^2(\bm{Y})} \left(\frac{1}{M} \sum_{m=1}^M p(\bm{Y} \vert \bm{\tilde{\theta}}_{m}) - p(\bm{Y}) \right)^2 \bigg\vert \bm{\theta},\bm{Y} \right] \right]+\mathcal{O}\left(\frac{1}{NM^2}\right)\nonumber\\
  \quad &=& \frac{1}{N} \mathbb{V} \left[ \log\left(\frac{p(\bm{Y} \vert \bm{\theta})}{p(\bm{Y})}\right)-\frac{1}{2}\frac{1}{M} \mathbb{V} \left[ \frac{p(\bm{Y} \vert \bm{\theta})}{p(\bm{Y})} \bigg\vert \bm{Y} \right] \right]+\mathcal{O}\left(\frac{1}{NM^2}\right)\nonumber\\
   \quad &=& \frac{1}{N} \left(\mathbb{V} \left[ \log\left(\frac{p(\bm{Y} \vert \bm{\theta})}{p(\bm{Y})}\right) \right]+\frac{1}{4M^2} \mathbb{V}\left[\mathbb{V} \left[ \frac{p(\bm{Y} \vert \bm{\theta})}{p(\bm{Y})} \bigg\vert \bm{Y} \right] \right] \right. \nonumber\\
   \quad &\quad &\left.-\frac{1}{M}Cov\left[\log\left(\frac{p(\bm{Y} \vert \bm{\theta})}{p(\bm{Y})}\right),\mathbb{V}\left[ \frac{p(\bm{Y} \vert \bm{\theta})}{p(\bm{Y})} \bigg\vert \bm{Y} \right] \right] \right)+\mathcal{O}\left(\frac{1}{NM^2}\right) \nonumber\\
   \quad &=& \frac{1}{N} \left(\mathbb{V} \left[ \log\left(\frac{p(\bm{Y} \vert \bm{\theta})}{p(\bm{Y})}\right) \right]+\frac{1}{4M^2} \mathbb{V}\left[\mathbb{V} \left[ \frac{p(\bm{Y} \vert \bm{\theta})}{p(\bm{Y})} \bigg\vert \bm{Y} \right] \right] \right. \nonumber\\
   \quad &\quad &\left.-\frac{1}{M}\left( \mathbb{E}\left[ \log\left( \frac{p(\bm{Y} \vert \bm{\theta})}{p(\bm{Y})} \right) \mathbb{V}\left[ \frac{p(\bm{Y} \vert \bm{\theta})}{p(\bm{Y})} \bigg\vert \bm{Y} \right]\right] \right.\right. \nonumber\\
   \quad &\quad &\left.\left.-\mathbb{E}\left[ \log\left( \frac{p(\bm{Y} \vert \bm{\theta})}{p(\bm{Y})} \right)\right]\mathbb{E}\left[\mathbb{V}\left[ \frac{p(\bm{Y} \vert \bm{\theta})}{p(\bm{Y})} \bigg\vert \bm{Y} \right]\right] \right) \right)+\mathcal{O}\left(\frac{1}{NM^2}\right)\nonumber\\
   \quad &=& \frac{1}{N}\mathbb{V} \left[ \log\left(\frac{p(\bm{Y} \vert \bm{\theta})}{p(\bm{Y})}\right) \right]-\frac{1}{NM}\left(\mathbb{E}\left[ \log\left( \frac{p(\bm{Y} \vert \bm{\theta})}{p(\bm{Y})} \right) \mathbb{V}\left[ \frac{p(\bm{Y} \vert \bm{\theta})}{p(\bm{Y})} \bigg\vert \bm{Y} \right]\right] \right. \nonumber \\
   \quad & \quad & \left.-\mathbb{E}\left[ \log\left( \frac{p(\bm{Y} \vert \bm{\theta})}{p(\bm{Y})} \right)\right]\mathbb{E}\left[\mathbb{V}\left[ \frac{p(\bm{Y} \vert \bm{\theta})}{p(\bm{Y})} \bigg\vert \bm{Y} \right]\right]\right) +\mathcal{O}\left(\frac{1}{NM^2}\right)
\end{eqnarray}

Similarly, for the second term of \eqref{A1var}, we obtain:
\begin{eqnarray}
   \quad & \quad & \frac{1}{N} \mathbb{E} \left[ \mathbb{V} \left[  \log \left(  \hat{p}_M(\bm{Y}) \right) \bigg\vert \bm{Y} \right] \right] \nonumber \\
   & =& \frac{1}{N} \mathbb{E}\left[ \mathbb{V} \left[- \log(p(\bm{Y})) + \frac{1}{p(\bm{Y})} \left(\frac{1}{M} \sum_{m=1}^M p(\bm{Y} \vert \bm{\tilde{\theta}}_{m}) - p(\bm{Y}) \right) \bigg\vert \bm{Y}  \right] \right] \nonumber\\
  & = & \frac{1}{N} \mathbb{E} \left[ \mathbb{V} \left[ \frac{1}{p(\bm{Y})} \left(\frac{1}{M} \sum_{m=1}^M p(\bm{Y} \vert \bm{\tilde{\theta}}_{m}) - p(\bm{Y}) \right) \bigg\vert \bm{Y}  \right] \right] \nonumber\\
   &=& \frac{1}{N} \mathbb{E}\left[  \mathbb{V} \left[ \frac{1}{M} \frac{ \sum_{m=1}^M p(\bm{Y} \vert \bm{\tilde{\theta}}_{m})}{p(\bm{Y})} \bigg\vert \bm{Y}\right]\right] \nonumber\\
 &=& \frac{1}{NM^2} \mathbb{E}  \left[ \mathbb{V} \left[ \frac{\sum_{m=1}^M p(\bm{Y} \vert \bm{\tilde{\theta}}_{m})}{p(\bm{Y})} \bigg\vert \bm{Y}\right] \right] \nonumber\\
&=& \frac{1}{NM} \mathbb{E} \left[\mathbb{V} \left[ \frac{ p(\bm{Y}\vert\bm{\theta})}{p(\bm{Y})} \bigg\vert \bm{Y}\right] \right].
\end{eqnarray}
The proof is completed by combining the two terms into the variance.
\end{proof}

\section{A leading-order expansion of the likelihood with respect to $N_e$}\label{proof_underflow}
\begin{proof}
The likelihood is defined as
\begin{eqnarray*}
 p(\bm{Y}_n \vert \bm{\theta}) = \prod_{i=1}^{N_e}  p(\bm{y}_{n,i} \vert \bm{\theta}),\;\; \forall n \leq N, \hbox{where}\;\; \bm{Y}_n=\{\bm{y}_{n,i}\}^{N_e}_{i=1},
\end{eqnarray*}
and the errors are distributed as $\bm{\epsilon}_i \sim \mathcal{N}\Big(0, \bm{\Sigma_\epsilon} \Big)$, $(\bm{\Sigma_\epsilon})_{jj} = \sigma^2_{\epsilon_j}$, and the unknown parameter value $\bm{\theta}$ follows the prior pdf, $\pi(\bm{\theta})$. Thus, the likelihood for the inner loop can be written as
\begin{eqnarray}
 p(\bm{Y}_n \vert \bm{\tilde{\theta}}_{n,m}) =  \left(2\pi \vert \bm{\Sigma_\epsilon} \vert\right)^{-\frac{N_e}{2}} \exp \left( -\frac{1}{2}  \sum_{i=1}^{N_e} \left\|\bm{y}_{n,i} - \bm{g}(\bm{\tilde{\theta}}_{n,m})\right\|^2_{\bm{\Sigma_\epsilon}^{-1}} \right),
\end{eqnarray}
where
\begin{eqnarray*}
 \left\|\bm{y}_{n,i} - \bm{g}(\bm{\tilde{\theta}}_{n,m})\right\|^2_{\bm{\Sigma_\epsilon}^{-1}}
 &=& \left( \bm{g}(\bm{\theta}_n) + \bm{\epsilon}_i - \bm{g}(\bm{\tilde{\theta}}_{n,m})\right)^T\bm{\Sigma_\epsilon}^{-1}\left( \bm{g}(\bm{\theta}_n) + \bm{\epsilon}_i- \bm{g}(\bm{\tilde{\theta}}_{n,m})\right)\\
 &=& \left( \bm{g}(\bm{\theta}_n)  - \bm{g}(\bm{\tilde{\theta}}_{n,m})\right)^T\bm{\Sigma_\epsilon}^{-1}\left( \bm{g}(\bm{\theta}_n) - \bm{g}(\bm{\tilde{\theta}}_{n,m})\right) \\
 & &  \vspace{5cm} + \bm{\epsilon}_i ^T\bm{\Sigma_\epsilon}^{-1}\left( \bm{g}(\bm{\theta}_n) - \bm{g}(\bm{\tilde{\theta}}_{n,m})\right) + \bm{\epsilon}_i ^T\bm{\Sigma_\epsilon}^{-1} \bm{\epsilon}_i \\
  &=& \left\| \bm{g}(\bm{\theta}_n)  - \bm{g}(\bm{\tilde{\theta}}_{n,m})\right\| ^2_{\bm{\Sigma_\epsilon}^{-1}}  +  \bm{\epsilon}_i ^T\bm{\Sigma_\epsilon}^{-1}\left( \bm{g}(\bm{\theta}_n) - \bm{g}(\bm{\tilde{\theta}}_{n,m})\right) + \left\| \bm{\epsilon}_i \right\|^2_{\bm{\Sigma_\epsilon}^{-1}}.
\end{eqnarray*}
We know that
\begin{eqnarray*}
  \sum_{i=1}^{N_e}  \left\| \bm{g}(\bm{\theta}_n)  - \bm{g}(\bm{\tilde{\theta}}_{n,m})\right\|^2_{\bm{\Sigma_\epsilon}^{-1}} = N_e  \left\| \bm{g}(\bm{\theta}_n)  - \bm{g}(\bm{\tilde{\theta}}_{n,m})\right\|^2_{\bm{\Sigma_\epsilon}^{-1}},
\end{eqnarray*}
which we expand to
\begin{eqnarray*}
\sum_{i=1}^{N_e}\bm{\epsilon}_i ^T\bm{\Sigma_\epsilon}^{-1}\left( \bm{g}(\bm{\theta}_n) - \bm{g}(\bm{\tilde{\theta}}_{n,m})\right)\;\;\; \hbox{and}\;\;\; \sum_{i=1}^{N_e} \left\| \bm{\epsilon}_i\right\|^2_{\bm{\Sigma_\epsilon}^{-1}}.
\end{eqnarray*}
We will use Kolmogorov's inequality for a random variable $X$:
\begin{eqnarray}\label{kolmo}
\mathbb{P}\left( \vert X - \mathbb{E}\left[ X \right] \vert \geq k \sqrt{\mathbb{V}\left[X\right]}\right) < \frac{1}{k^2},
\end{eqnarray}
where $k$ is a real number greater than $1$. 

First,
\begin{eqnarray*}
  \sum_{i=1}^{N_e} \bm{\epsilon}_i ^T\bm{\Sigma_\epsilon}^{-1}\left( \bm{g}(\bm{\theta}_n) - \bm{g}(\bm{\tilde{\theta}}_{n,m})\right) 
&=&   \sum_{i=1}^{N_e}   \sum_{j=1}^{q} \frac{\epsilon_{i,j}}{\sigma_{\epsilon_j}^2}\left(\bm{g}_j(\bm{\theta}_n) - \bm{g}_j(\bm{\tilde{\theta}}_{n,m})\right) \nonumber \\
\quad &=& \sum_{j=1}^{q} \left( \sum_{i=1}^{N_e}\epsilon_{i,j} \right)\frac{1}{\sigma_{\epsilon_j}^2}\left(\bm{g}_j(\bm{\theta}_n) - \bm{g}_j(\bm{\tilde{\theta}}_{n,m})\right).
\end{eqnarray*}
We apply \eqref{kolmo} to attain $\mathbb{P}\left( \vert \sum_{i=1}^{N_e} \epsilon_{i,j} \vert \geq k\sigma_{\epsilon_{j}}\sqrt{N_e} \right) < \frac{1}{k^2}$, which leads to $\sum_{i=1}^{N_e} \epsilon_{i,j}=\sigma_{\epsilon_j}\mathcal{O}_{\mathbb{P}}\left(\sqrt{N_e}\right)$ and, thus,
\begin{eqnarray*}
  \sum_{i=1}^{N_e} \bm{\epsilon}_i ^T\bm{\Sigma_\epsilon}^{-1}\left( \bm{g}(\bm{\theta}_n) - \bm{g}(\bm{\tilde{\theta}}_{n,m})\right)    
   &= & \sum_{j=1}^{q} \left( \sigma_{\epsilon_j} \mathcal{O}_{\mathbb{P}} \left( \sqrt{N_e} \right) \right)
    \frac{1}{\sigma_{\epsilon_k}^2}\left(\bm{g}_j(\bm{\theta}_n) - \bm{g}_j(\bm{\tilde{\theta}}_{n,m})\right)  \\
  &= & \sum_{j=1}^{q}  \frac{1}{\sigma_{\epsilon_j}}\left(g_j(\bm{\theta}_n) - g_j(\bm{\tilde{\theta}}_{n,m})\right)\mathcal{O}_{\mathbb{P}} \left( \sqrt{N_e} \right)   \\
    &= & \bm{v_\epsilon}^{T}\left(\bm{g}(\bm{\theta}_n) - \bm{g}(\bm{\tilde{\theta}}_{n,m})\right)  \mathcal{O}_{\mathbb{P}}\left( \sqrt{N_e} \right),
\end{eqnarray*}
where $\bm{v_\epsilon}$ is the vector of the diagonal elements of $\bm{\Sigma_{\epsilon}}^{-1/2}$. For the other term, we have
\begin{eqnarray*}
 \sum_{i=1}^{N_e} \left\| \bm{\epsilon}_i\right\|^2_{\bm{\Sigma_\epsilon}^{-1}} 
 &=&  \sum_{i=1}^{N_e} \bm{\epsilon}_i^T \bm{\Sigma_\epsilon}^{-1} \bm{\epsilon}_i  =  \sum_{i=1}^{N_e}\sum_{j=1}^{q} \frac{\epsilon^2_{i,j}}{\sigma_{\epsilon_j}^2} \\
  & =& \sum_{j=1}^{q} \frac{1}{\sigma_{\epsilon_j}^2} \left(\sum_{i=1}^{N_e}{\epsilon^2_{i,j}}\right).
\end{eqnarray*}

The random variable $\sum_{i=1}^{N_e}{\epsilon^2_{i,j}}$ follows the chi-squared distribution with mean $\sigma{\epsilon_{j}}N_e$ and variance $2\sigma_{\epsilon_j}^2N_e$. Similarly, we apply \eqref{kolmo} to obtain $\sum_{i=1}^{N_e}\epsilon^2_{i,j}=-\sigma_{\epsilon_{j}}N_e+\sum_{i=1}^{N_e}(\epsilon^2_{i,j}+\sigma_{\epsilon_{j}}) = -\sigma_{\epsilon_{j}}N_e+\sigma_{\epsilon_{j}}\mathcal{O}_{\mathbb{P}}\left(\sqrt{N_e}\right)$; as a result,
\begin{eqnarray*}
 \sum_{i=1}^{N_e} \left\| \bm{\epsilon}_i\right\|^2_{\bm{\Sigma_\epsilon}^{-1}} 
  & =& \sum_{j=1}^{q} \frac{1}{\sigma_{\epsilon_j}^2} \left( \sum_{i=1}^{N_e}\epsilon^2_{i,j}\right)  =  \sum_{j=1}^{q} \sigma_{\epsilon_{j}}^{-1}\left(N_e + \mathcal{O}_{\mathbb{P}}\left(\sqrt{N_e}\right) \right).
\end{eqnarray*}

Finally, we end up with
\begin{eqnarray*}
  p(\bm{Y}_n \vert \bm{\tilde{\theta}}_{n,m}) &=&  \left(2\pi \vert \bm{\Sigma_\epsilon} \vert \right)^{-\frac{N_e}{2}}
   \exp \left( -\frac{1}{2}N_e\left\|\bm{g}(\bm{\theta}_{n}) - \bm{g}(\bm{\tilde{\theta}}_{n,m})\right\|_{\bm{\Sigma_\epsilon}^{-1}}^2 \right. \nonumber \\
 && \left. - \frac{1}{2} \bm{v_\epsilon}^{T}\left(\bm{g}(\bm{\theta}_n) - \bm{g}(\bm{\tilde{\theta}}_{n,m})\right)\mathcal{O}_{\mathbb{P}}\left(\sqrt{N_e}\right) - \frac{1}{2}\sum_{j=1}^{q} \sigma_{\epsilon_{j}}^{-1}\left(N_e + \mathcal{O}_{\mathbb{P}}\left(\sqrt{N_e}\right) \right) \right).
\end{eqnarray*}

\end{proof}

\section{Laplace approximation of posterior distributions}
We derive the Laplace posterior pdf $\pi_{la}(\bm{\theta} \vert \bm{Y}) \sim \mathcal{N} \left(\bm{\hat{\theta}}, \bm{\hat{\Sigma}}(\bm{\hat{\theta}})\right)$, which is an approximation of the posterior pdf, $\pi(\bm{\theta} \vert \bm{Y})$. 
\begin{proof}
We let $F$ be the negative logarithmic of the posterior distribution:
\begin{eqnarray}
F(\bm{\theta}) &=& -\log\left( \pi(\bm{\theta} \vert \bm{Y})\right).
\end{eqnarray}

The mode $\hat{\bm{\theta}}$ of the posterior pdf is given by
\begin{eqnarray}
  \bm{\hat{\theta}} \myeq \underset{\bm{\theta} \in \Theta}{\arg\min} \ F(\bm{\theta}),
\end{eqnarray}
which implies that $ \nabla_{\bm{\theta}} F(\bm{\hat{\theta}})=0$. By matching the second-order Taylor expansion of $F$ around $\bm{\hat{\theta}}$ and $-\log\left(\pi_{la}(\bm{\theta} \vert \bm{Y})\right)$ we obtain $\bm{\hat{\Sigma}}$ as the inverse Hessian matrix of the negative logarithm of the posterior pdf evaluated at $\bm{\theta} = \bm{\hat{\theta}}$, i.e.,
\begin{eqnarray}
\label{la:thetahat1}
    \bm{\hat{\Sigma}} =\Big( \nabla_{\bm{\theta}} \nabla_{\bm{\theta}} F(\bm{\hat{\theta}})\Big)^{-1}.
\end{eqnarray}
We use $\nabla_{\bm{\theta}} F(\bm{\hat{\theta}})=0$ in the evaluation at $\bm{\theta} = \bm{\theta}_t$ of the Taylor expansion, of $\nabla_{\bm{\theta}} F$, around $\bm{\hat{\theta}}$ to write
\begin{eqnarray}
\label{la:sigma1}
\bm{\hat{\theta}} =   \bm{\theta}_t -  \nabla_{\bm{\theta}} F(\bm{\theta}_t)\left( \nabla_{\bm{\theta}} \nabla_{\bm{\theta}} F(\bm{\theta}_t)\right)^{-1}+ \mathcal{O}\left(\left\|\bm{\hat{\theta}} - \bm{\theta}_t\right\|^2\right).
\end{eqnarray}
In order to express these moments, (\ref{la:thetahat1}) and (\ref{la:sigma1}), in terms of $\bm{g}$ rather than $F$, we further expand by introducing
\begin{eqnarray*}
  \bm{E_\epsilon}(\bm{\theta})  = \sum_{i=1}^{N_e} \bm{r}_i^T(\bm{\theta}), \;\;\;
 \bm{J}(\bm{\theta}) =- \nabla_{\bm{\theta}} \bm{g}(\bm{\theta}),\;\;\; 
  \bm{H}(\bm{\theta}) =-\nabla_{\bm{\theta}}  \nabla_{\bm{\theta}} \bm{g}(\bm{\theta}),\;\;\;\hbox{and}\;\;\;
  h(\bm{\theta}) = \log (\pi(\bm{\theta})),
\end{eqnarray*}
where $\bm{J}(\bm{\theta})$ and $\bm{H}(\bm{\theta})$ are the Jacobian and the Hessian of $-\bm{g}(\bm{\theta})$, respectively. For $\bm{\theta} = \bm{\hat{\theta}}$, we have $ \bm{E_\epsilon}(\bm{\hat{\theta}}) \sim \mathcal{N}\Big(N_e(\bm{g}(\bm{\theta_t}) - \bm{g}(\bm{\hat{\theta}})), N_e\bm{\Sigma_\epsilon}\Big)$.

This allows us to recast the moments as
\begin{eqnarray}
  \hat{\bm{\theta}} = \bm{\theta}_t - \Big( \nabla_{\bm{\theta}}  h(\bm{\theta}_t) + \bm{E_\epsilon} \bm{\Sigma_\epsilon}^{-1}\bm{J}(\bm{\theta}_t)\Big) \Big( N_e \bm{J}(\bm{\theta}_t)^T\bm{\Sigma_\epsilon}^{-1}\bm{J}(\bm{\theta}_t) \nonumber \\
  - \nabla_{\bm{\theta}}  \nabla_{\bm{\theta}}  h(\bm{\theta}_t)  + \bm{H}(\bm{\theta}_t)^T\bm{\Sigma_\epsilon}^{-1}\bm{E_\epsilon} \Big)^{-1} + \mathcal{O}\left( \left\|\bm{\hat{\theta}} - \bm{\theta}_t \right\| \right),
\end{eqnarray}
and
\begin{eqnarray}
 \bm{\hat{\Sigma}} = \Big( N_e \bm{J}(\bm{\hat{\theta}})^T\bm{\Sigma_\epsilon}^{-1} \bm{J}(\bm{\hat{\theta}}) -  \nabla_{\bm{\theta}}  \nabla_{\bm{\theta}}  h(\bm{\hat{\theta}})  + \bm{H}(\bm{\hat{\theta}})^T \bm{\Sigma_\epsilon}^{-1} \bm{E_\epsilon} \Big)^{-1}.
\end{eqnarray}
Furthermore, Long et al. \cite{key15} show that
\begin{eqnarray*}
  \quad & \bm{J}(\bm{\hat{\theta}})^T \bm{\Sigma_\epsilon}^{-1} \bm{E_\epsilon} = \mathcal{O}_\mathbb{P} \left(\sqrt{N_e}\right),
\;\;
  \bm{H}(\bm{\hat{\theta}})^T \bm{\Sigma_\epsilon}^{-1} \bm{E_\epsilon} = \mathcal{O}_\mathbb{P} \left(\sqrt{N_e}\right)\;\; \\ \quad & \hbox{and}\;\;  N_e \bm{J}(\bm{\hat{\theta}})^T\bm{\Sigma_\epsilon}^{-1} \bm{J}(\bm{\hat{\theta}}) =  \mathcal{O}_\mathbb{P} \left(N_e\right).
\end{eqnarray*} 
Taking into account these probabilistic rates with respect to $N_e$, we state the approximations of the moments to complete the proof as
\begin{eqnarray}\label{eq_theta1}
\bm{\hat{\theta}} = \bm{\theta}_t -(N_e \bm{J}^T(\bm{\theta}_t) \bm{\Sigma_\epsilon}^{-1} \bm{J}(\bm{\theta}_t) + \bm{H}^T \bm{\Sigma_\epsilon}^{-1} \bm{E_\epsilon} - \nabla \nabla h(\bm{\theta}_t))^{-1} \bm{J}^T \bm{\Sigma_\epsilon}^{-1} \bm{E_\epsilon} + \mathcal{O}_\mathbb{P}\left(\frac{1}{N_e}\right)
\end{eqnarray}
and
\begin{eqnarray}\label{eq_sigma1}
\bm{\hat{\Sigma}}^{-1} =  N_e \bm{J}(\bm{\hat{\theta}})^T\bm{\Sigma_\epsilon}^{-1} \bm{J}(\bm{\hat{\theta}}) -  \nabla_{\bm{\theta}}  \nabla_{\bm{\theta}}  h(\bm{\hat{\theta}}) +  \mathcal{O}_\mathbb{P}\left(\sqrt{N_e}\right).
\end{eqnarray}
\end{proof}

\section{Expected information gain with Laplace method} \label{ap:secB}
The derivation of the expected information gain with the Laplace method.
\begin{proof}
We consider 
\begin{eqnarray}
\tilde{\pi}(\bm{\theta} \vert \bm{Y}) = (2 \pi)^{-\frac{d}{2}} \vert\bm{\hat{\Sigma}}\vert^{-\frac{1}{2}} \exp\left(-\frac{1}{2} \|\bm{\theta} - \bm{\hat{\theta}}\|^{2}_{\bm{\hat{\Sigma}}^{-1}}\right),  
\end{eqnarray}
which is the Gaussian approximation of the posterior pdf, for rewriting the Kullback-Leibler divergence,
\begin{eqnarray}
D_{kl} = \int_{\Theta} \pi(\bm{\theta} \vert \bm{Y}) \log \left( \frac{\pi(\bm{\theta} \vert \bm{Y})}{\pi(\bm{\theta})} \right)\, \text{d}\bm{\theta},
\end{eqnarray} 
into
\begin{eqnarray}
 D_{kl} &=& \int_{\Theta} \log \left( \frac{\pi(\bm{\theta} \vert \bm{Y})}{\pi(\bm{\theta})} \right)\,\tilde{\pi}(\bm{\theta} \vert \bm{Y}) \text{d}\bm{\theta} + \underbrace{\int_{\Theta}  \log \left( \frac{\pi(\bm{\theta} \vert \bm{Y})}{\pi(\bm{\theta})} \right) \left(\pi(\bm{\theta} \vert \bm{Y}) - \tilde{\pi}(\bm{\theta} \vert \bm{Y})\right)\, \text{d}\bm{\theta}}_{\varepsilon_{int}}.
\end{eqnarray}
Then, by using the decomposition
\begin{eqnarray}
  \log \left( \frac{\pi(\bm{\theta} \vert \bm{Y})}{\pi(\bm{\theta})}\right) =   \log \left( \frac{\pi(\bm{\theta} \vert \bm{Y})}{\tilde{\pi}(\bm{\theta} \vert \bm{Y})}\right) +  \log \left( \frac{\tilde{\pi}(\bm{\theta} \vert \bm{Y})}{\pi(\bm{\theta})}\right),
\end{eqnarray}
we can write
\begin{eqnarray}\label{eq_dkl}
 D_{kl} &=& \int_{\Theta} \left[\underbrace{ \log \left( \frac{\pi(\bm{\theta} \vert \bm{Y})}{\tilde{\pi}(\bm{\theta} \vert \bm{Y})}\right)}_{\varepsilon_{la}} +  \log \left( \frac{\tilde{\pi}(\bm{\theta} \vert \bm{Y})}{\pi(\bm{\theta})}\right)\right] \tilde{\pi}(\bm{\theta} \vert \bm{Y})\, \text{d}\bm{\theta} + \varepsilon_{int}.
\end{eqnarray}

The logarithmic error introduced by replacing the posterior density $\pi(\bm{\theta} \vert \bm{Y})$ by its Gaussian approximation $\tilde{\pi}(\bm{\theta} \vert \bm{Y})$ is here denoted by $\varepsilon_{la}$. Furthermore, we know that
\begin{eqnarray}
   \log \left( \frac{\tilde{\pi}(\bm{\theta} \vert \bm{Y})}{\pi(\bm{\theta})}\right) &=& \log \left( \tilde{\pi}(\bm{\theta} \vert \bm{Y})\right) -\underbrace{\log \left(\pi(\bm{\theta})\right)}_{h(\bm{\theta})}\nonumber\\
   &=& -\frac{1}{2} \log((2 \pi)^{d} \vert \bm{\hat{\Sigma}} \vert)  - \frac{1}{2} (\bm{\theta} - \bm{\hat{\theta}})^T\bm{\hat{\Sigma}}^{-1}(\bm{\theta} - \bm{\hat{\theta}}) - h(\bm{\theta}),
\end{eqnarray}
which allows us to express $D_{kl}$ as
\begin{eqnarray}
    D_{kl} = \int_{\Theta} \left[ -\frac{1}{2} \log((2 \pi)^{d} \vert \bm{\hat{\Sigma}} \vert)  - \frac{1}{2} (\bm{\theta} - \bm{\hat{\theta}})^T\bm{\hat{\Sigma}}^{-1}(\bm{\theta} - \bm{\hat{\theta}}) - h(\bm{\theta}) + \varepsilon_{la} \right] \tilde{\pi}(\bm{\theta} |y)\, \text{d}\bm{\theta} + \varepsilon_{int}. \nonumber
\end{eqnarray}
In \cite{key15}, it was shown that
\begin{eqnarray}
\int_{\Theta}  \varepsilon_{la}  \tilde{\pi}(\bm{\theta} \vert \bm{Y})\, \text{d}\bm{\theta} = \mathcal{O}_\mathbb{P}\left(\frac{1}{N_e^2}\right),\;\;\;\;\;\;\;
\varepsilon_{int} = \mathcal{O}_\mathbb{P}\left(\frac{1}{N_e^2}\right), \nonumber
\end{eqnarray}
and
\begin{eqnarray}
\int_{\Theta} h(\bm{\theta}) \tilde{\pi}(\bm{\theta} \vert \bm{Y})\, \text{d}\bm{\theta} =  - h(\bm{\hat{\theta}}) - \frac{1}{2}\hbox{tr}\left(\bm{\hat{\Sigma}}: \nabla_{\bm{\theta}} \nabla_{\bm{\theta}} h(\bm{\hat{\theta}})\right), \nonumber
\end{eqnarray}
which leads us to
\begin{eqnarray}
    D_{kl} &=& \int_{\Theta} \left[ -\frac{1}{2} \log((2 \pi)^{d} \vert \bm{\hat{\Sigma}} \vert)  - \frac{1}{2} (\bm{\theta} - \bm{\hat{\theta}})^T\bm{\hat{\Sigma}}^{-1}(\bm{\theta} - \bm{\hat{\theta}}) - h(\bm{\theta}) \right] \tilde{\pi}(\bm{\theta} \vert \bm{Y})\, \text{d}\bm{\theta} + \mathcal{O}_\mathbb{P}\left(\frac{1}{N_e^2}\right),\nonumber\\
    &=&-\frac{1}{2} \log((2 \pi)^{d} \vert \bm{\hat{\Sigma}} \vert)  -\frac{d}{2} - h(\bm{\hat{\theta}}) - \frac{\hbox{tr}\left( \bm{\hat{\Sigma}}: \nabla_{\bm{\theta}} \nabla_{\bm{\theta}} h(\bm{\hat{\theta}})\right)}{2} +  \mathcal{O}_\mathbb{P}\left(\frac{1}{N_e^2}\right).\nonumber
\end{eqnarray}
Now, the expected information gain, $I=\mathbb{E}[D_{kl}]$, can be written as:
\begin{eqnarray}
 I &=& \int_{\bm{\mathcal{Y}}} D_{kl} p(\bm{Y}) d \bm{Y}  \nonumber\\
   &=& \int_{\Theta} \int_{\bm{\mathcal{Y}}} D_{kl} p(\bm{Y} \vert \bm{\theta}_t) d \bm{Y} \pi(\bm{\theta}_t) \text{d}\bm{\theta}_t\nonumber\\
 &=& \int_{\Theta} \int_{\bm{\mathcal{Y}}}  \left[ -\frac{1}{2} \log((2 \pi)^{d} \vert\bm{\hat{\Sigma}}\vert)  -\frac{d}{2} - h(\bm{\hat{\theta}}) - \frac{1}{2}\hbox{tr}\left( \bm{\hat{\Sigma}}: \nabla_{\bm{\theta}} \nabla_{\bm{\theta}} h(\bm{\hat{\theta}})\right) \right]  p(\bm{Y} \vert \bm{\theta}_t) d \bm{Y} \pi(\bm{\theta}_t) \text{d}\bm{\theta}_t +  \mathcal{O}\left(\frac{1}{N_e^2}\right).\nonumber
\end{eqnarray}
Finally, to conclude the proof we use approximation \eqref{eq_theta1} to obtain
\begin{eqnarray}
I = \int_{\Theta}{\left[-\frac{1}{2} \log((2 \pi)^{d} \vert \bm{\hat{\Sigma}} \vert) -\frac{d}{2} - h(\bm{{\theta}_t}) \right] p(\bm{\theta}_t) d\bm{\theta}_t} + \mathcal{O}\left(\frac{1}{N_e}\right).  
\end{eqnarray}
\end{proof}

\section{Effect of the change of measure on the arithmetic underflow}\label{correct_underflow}
We show the effect of the proposed importance sampling approach on the inner loop computation of DLMCIS with respect to the number of repetitive experiments, $N_e$. The derivation below follows closely the derivation in \ref{proof_underflow}. 
\begin{proof}
By using the importance sampling based on the Laplace approximation, we obtain
\begin{eqnarray}
 L(\bm{Y}_n;\bm{\tilde{\theta}}_{n,m}) =  \frac{p (\bm{Y}_n \vert \bm{\tilde{\theta}}_{n,m} ) \pi(\bm{\tilde{\theta}}_{n,m})}{\tilde{\pi}_n (\bm{\tilde{\theta}}_{n,m})},
\end{eqnarray}
where
\begin{eqnarray*}
  p(\bm{Y}_n \vert \bm{\tilde{\theta}}_{n,m}) &=&  
   \left(2\pi \vert \bm{\Sigma_\epsilon} \vert \right)^{-\frac{N_e}{2}} \exp \left( -\frac{1}{2}N_e\left\|\bm{g}(\bm{\theta}_n) - \bm{g}(\bm{\tilde{\theta}}_{n,m})\right\|_{\bm{\Sigma_\epsilon}^{-1}}^2  \right. \nonumber \\
 && \left.  - \frac{1}{2}\bm{v_\epsilon}^{T}\left(\bm{g}(\bm{\theta}_n) - \bm{g}(\bm{\tilde{\theta}}_{n,m})\right)\mathcal{O}_{\mathbb{P}}\left(\sqrt{N_e}\right) - \frac{1}{2}\sum_{j=1}^{q} \sigma_{\epsilon_{j}}^{-1}\left(N_e + \mathcal{O}_{\mathbb{P}}\left(\sqrt{N_e}\right) \right) \right). \nonumber
\end{eqnarray*}
The distribution $\tilde{\pi}$ is the Gaussian distribution $\tilde{\pi}_n \sim \mathcal{N} (\bm{\hat{\theta}}_n, \bm{\hat{\Sigma}}(\bm{\hat{\theta}}_n))$:
\begin{eqnarray*}
  \tilde{\pi}_n(\bm{\tilde{\theta}}_
  {n,m}) &=&  \left(2\pi \vert \bm{\hat{\Sigma}}(\bm{\hat{\theta}}_n) \vert \right)^{-\frac{1}{2}} \exp\left(-\frac{1}{2} \left(\bm{\tilde{\theta}}_{n,m} - \bm{\hat{\theta}}_n\right)^T\bm{\hat{\Sigma}}^{-1}(\bm{\hat{\theta}}_n) \left(\bm{\tilde{\theta}}_{n,m} - \bm{\hat{\theta}}_n\right)\right),
\end{eqnarray*}
where the inverse of the covariance matrix is given by
\begin{eqnarray*}
 \bm{\hat{\Sigma}}^{-1}(\bm{\hat{\theta}}_n) = N_e \bm{J}(\bm{\hat{\theta}}_n)^T \bm{\Sigma_\epsilon}^{-1}
                                                     \bm{J}(\bm{\hat{\theta}}_n) - \nabla_{\bm{\theta}} \nabla_{\bm{\theta}} h(\bm{\hat{\theta}}_n) + \mathcal{O}_\mathbb{P}\left( \sqrt{N_e}\right).
\end{eqnarray*}
This allows us to write
\begin{eqnarray*}
  \tilde{\pi}_n(\bm{\tilde{\theta}}_{n,m}) 
  &=&  
   \left(2\pi \vert \bm{\hat{\Sigma}}(\bm{\hat{\theta}}_n)\vert \right)^{-\frac{1}{2}} \exp\left( -\frac{N_e}{2} \left(\bm{\tilde{\theta}}_{n,m} - \bm{\hat{\theta}}_n\right)^T
    \bm{J}(\bm{\hat{\theta}}_n)^T \bm{\Sigma_\epsilon}^{-1} \bm{J}(\bm{\hat{\theta}}_n) \left(\bm{\tilde{\theta}}_{n,m} - \bm{\hat{\theta}}_n\right) \right. \\
    && \hspace{2cm} \left.  +\;  \frac{1}{2} \left(\bm{\tilde{\theta}}_{n,m} - \bm{\hat{\theta}}_n\right)^T  \nabla_{\bm{\theta}} \nabla_{\bm{\theta}} h(\bm{\hat{\theta}}_n) \left(\bm{\tilde{\theta}}_{n,m} - \bm{\hat{\theta}}_n\right) + 
    \frac{1}{2}\mathcal{O}_\mathbb{P}\left( \sqrt{N_e}\right) \right)\\
  &=&  \left(2\pi \vert \bm{\hat{\Sigma}}(\bm{\hat{\theta}}_n) \vert \right)^{-\frac{1}{2}} \exp\left( -\frac{N_e}{2} \left(\bm{J}(\bm{\hat{\theta}}_n)\left(\bm{\tilde{\theta}}_{n,m} - \bm{\hat{\theta}}_n\right)\right)^T
 \bm{\Sigma_\epsilon}^{-1} \left(\bm{J}(\bm{\hat{\theta}}_n)\left(\bm{\tilde{\theta}}_{n,m} - \bm{\hat{\theta}}_n\right)\right) \right. \\
    && \hspace{2cm} \left. +\; \frac{1}{2} \left(\bm{\tilde{\theta}}_{n,m} - \bm{\hat{\theta}}_n\right)^T  \nabla_{\bm{\theta}} \nabla_{\bm{\theta}} h(\bm{\hat{\theta}}_n) \left(\bm{\tilde{\theta}}_{n,m} - \bm{\hat{\theta}}_n\right) + \mathcal{O}_\mathbb{P}\left( \sqrt{N_e}\right) \right)\\
      &=&  \left(2\pi \vert \bm{\hat{\Sigma}}(\bm{\hat{\theta}}_n) \vert \right)^{-\frac{1}{2}}  \exp\left(-\frac{N_e}{2} \left\|\bm{J}(\bm{\hat{\theta}}_n)\left(\bm{\tilde{\theta}}_{n,m} - \bm{\hat{\theta}}_n\right)\right\|^2_{ \bm{\Sigma_\epsilon}^{-1}} \right.\\
   && \hspace{2cm} \left.  + \frac{1}{2} \left\|\bm{\tilde{\theta}}_{n,m} - \bm{\hat{\theta}}_n \right\|^2_{ \nabla_{\bm{\theta}} \nabla_{\bm{\theta}} h(\bm{\hat{\theta}}_n)} +\;\; \mathcal{O}_\mathbb{P}\left( \sqrt{N_e}\right) \right).
\end{eqnarray*}
Dividing $p(\bm{Y}_n \vert \bm{\tilde{\theta}}_{n,m}) \pi(\bm{\tilde{\theta}}_{n,m}) $ by $\tilde{\pi}_n(\bm{\tilde{\theta}}_{n,m})$ yields
\begin{eqnarray*}
 L(\bm{Y}_n;\bm{\tilde{\theta}}_{n,m}) &=&  \frac{\left(2\pi \vert \bm{\Sigma_\epsilon} \vert \right)^{-\frac{N_e}{2}}}{\left(2\pi \vert \bm{\hat{\Sigma}}(\bm{\hat{\theta}}_n) \vert \right)^{-\frac{1}{2}}}  \pi \left( \bm{\tilde{\theta}}_{n,m}\right)  \exp \left( - \frac{N_e}{2} \left\|\bm{g}(\bm{\theta}_n) - \bm{g}(\bm{\tilde{\theta}}_{n,m})\right\|^2_{\bm{\Sigma_\epsilon}^{-1}}
\right.\nonumber \\
 && \left. +\frac{N_e}{2}  \left\|\bm{J}(\bm{\hat{\theta}}_n)\left(\bm{\tilde{\theta}}_{n,m} - \bm{\hat{\theta}}_n\right)\right\|^2_{ \bm{\Sigma_\epsilon}^{-1}} 
   - \frac{1}{2}  \left\| \bm{\tilde{\theta}}_{n,m} - \bm{\hat{\theta}}_n \right\|^2_{ \nabla_{\bm{\theta}} \nabla_{\bm{\theta}} h(\bm{\hat{\theta}}_n)} + \mathcal{O}_\mathbb{P}\left(\sqrt{N_e}\right) \right. \nonumber\\
  &&\left.  - \frac{1}{2}\bm{v_\epsilon}^{T}\left(\bm{g}(\bm{\theta}_n) - \bm{g}(\bm{\tilde{\theta}}_{n,m})\right)\mathcal{O}_{\mathbb{P}}\left(\sqrt{N_e}\right) - \frac{1}{2}\sum_{j=1}^{q} \sigma_{\epsilon_{j}}^{-1}\left(N_e + \mathcal{O}_{\mathbb{P}}\left(\sqrt{N_e}\right) \right) \right).
\end{eqnarray*}

Then we apply the logarithm to obtain the log-likelihood expression.

\end{proof}

\end{document}